\documentclass[10pt]{amsart}
\usepackage[usenames,dvipsnames]{xcolor}
\usepackage{etex}
\usepackage{mathtools,tensor}
\usepackage{latexsym,amssymb,amsthm,amsmath,enumerate,epsfig,setspace,bbding,color,graphicx,caption}
\usepackage[noadjust]{cite}
\usepackage{stmaryrd}
\usepackage[all]{xy}
\usepackage[retainorgcmds]{IEEEtrantools}
\usepackage{hyperref}
\usepackage{subfiles}
\usepackage[margin=10pt,font=small,labelfont=bf]{caption}
\usepackage{tikz-cd}
\usepackage{stmaryrd}
\usetikzlibrary{arrows}
\usepackage[normalem]{ulem}

\newtheorem{thm}{Theorem}[section]
\newtheorem{Thm}{Theorem}
\newtheorem{cor}[thm]{Corollary}
\newtheorem{Cor}{Corollary}
\newtheorem{Prop}{Proposition}
\newtheorem{prop}[thm]{Proposition}
\newtheorem{lemma}[thm]{Lemma}
\newtheorem*{Q}{Question}

\theoremstyle{definition}
\newtheorem{defi}[thm]{Definition}

\newtheorem{Ex}{Example}
\newtheorem{ex}[thm]{Example}
\newtheorem{rem}[thm]{Remark}

\pgfmathsetmacro{\shift}{0.65ex}
\renewcommand{\L}{\mathcal{L}}

\newcommand{\G}{\mathcal{G}}

\newcommand{\U}{\text{$\mathcal{U}$}}
\newcommand{\g}{\mathfrak{g}}
\newcommand{\h}{\mathfrak{h}}
\renewcommand{\t}{\mathfrak{t}}

\renewcommand{\H}{\mathcal{H}}

\newcommand{\T}{\mathbb{T}}
\newcommand{\A}{\mathcal{A}}
\newcommand{\No}{\mathcal{N}}
\newcommand{\N}{\mathbb{N}}
\newcommand{\R}{\mathbb{R}}

\newcommand{\Z}{\mathbb{Z}}
\newcommand{\D}{\mathcal{D}}

\renewcommand{\dim}{\text{dim}}
\renewcommand{\ker}{\mathrm{Ker}}

\renewcommand{\O}{\mathcal{O}}

\newcommand\sS{\mathcal{S}}

\renewcommand{\d}{\mathrm{d}}

\renewcommand{\P}{\mathcal{P}}

\renewcommand{\Q}{\mathcal{Q}}

\newcommand{\p}{\mathfrak{p}}
\renewcommand{\c}{\mathfrak{c}}

\newcommand{\Af}{\emph{Af}}

\begin{document}


\newcommand{\Addresses}{{
  \bigskip
  \footnotesize

  \textsc{University of Toronto, Canada}\par\nopagebreak
  \textit{E-mail address}: \texttt{maarten.mol.math@gmail.com}

}}

\title[\resizebox{5.5in}{!}{Constructibility of momentum maps and linear variation for singular symplectic reduced spaces}]{Constructibility of momentum maps and linear variation for singular symplectic reduced spaces}
\author{Maarten Mol\\
University of Toronto}
\date{}
\maketitle
\begin{abstract} In this paper we show that the transverse image of the momentum map of a Hamiltonian Lie group action admits a natural integral affine stratification with the property that over each stratum the momentum map is an equivariantly locally trivial fibration, provided the group is compact and the momentum map is proper. Using this we extend the linear variation theorem of Duistermaat and Heckman to singular values of the momentum map by showing that the cohomology classes of the symplectic forms on the reduced spaces at values within a stratum vary linearly. We also point out an instance of an invariant cycle theorem for momentum maps. Finally, we extend all of the above to Hamiltonian actions of proper quasi-symplectic groupoids. 
\end{abstract}
\tableofcontents

\section*{Introduction} Proper maps in various categories studied in singularity theory (such as the semi-algebraic, subanalytic, complex algebraic and complex analytic categories) are known to be constructible, in the sense that the image of the map can be stratified (with strata in the given category) in such a way that the map is a topologically locally trivial fibration over each stratum (see, e.g., \cite[\S1.5-1.7]{GoMac}). This is one of the fundamental facts on how families of such spaces (encoded as the fibers of a map in the corresponding category) behave. 
With the aim of gaining a better understanding of how the symplectic reduced spaces of a Hamiltonian action of a compact Lie group (or, more generally, of a proper quasi-symplectic groupoid) vary when reducing at different values of the momentum map, in this paper we show that the momentum map of such an action is constructible in a similar equivariant sense, whenever it is proper. That is, we show that there is a natural stratification of the projection of the momentum map image to the coadjoint orbit space (or, more generally, the orbit space of the quasi-symplectic groupoid) with the property that the momentum map is equivariantly locally trivial over each stratum. This stratification is determined entirely by the data of the isotropy groups of the action and those of the coadjoint action (or, more generally, those of the quasi-symplectic groupoid). Moreover, its strata are integral affine manifolds. Since the underlying topological space of the symplectic reduced spaces at values within a given stratum does not change, for such values we can compare the cohomology classes of the symplectic forms and the corresponding symplectic volumes of these spaces. As we will show, these cohomology classes vary linearly and the symplectic volumes form a polynomial function with respect to the integral affine structure on such a stratum, extending the Duistermaat-Heckman linear variation theorem \cite{DuHe} and the linear variation theorem for regular Poisson manifolds of compact types \cite{CrFeTo}. \\
 
In the rest of this introduction we explain our results in more detail for Hamiltonian actions of compact and connected Lie groups, postponing the details of their extension to actions of possibly disconnected compact Lie groups and proper quasi-symplectic groupoids to \S\ref{sec:extension:disconnectedgps} and \S\ref{sec:extension:quasisympgpoids}. We focus on this case first for the sake of simplicity and because (to the best of our knowledge) already in this case our results have not appeared in the literature before. \\

So, let $G$ be a compact connected Lie group and let $J:(S,\omega)\to \g^*$ be a Hamiltonian $G$-space with the property that the momentum map is proper as map into its image $\Delta:=J(S)$. Fix a maximal torus $T$ in $G$ and a closed Weyl chamber $\t^*_+\subset \t^*$. We view $\t^*_+$ as subspace of $\g^*$ via the canonical inclusion of $\t^*$ into $\g^*$ (recalled in \S\ref{sec:background:Lie}). As such, $\t^*_+$ is a fundamental domain for the coadjoint $G$-action, in the sense that this inclusion induces a homeomorphism between $\t^*_+$ and $\underline{\g^*}:=\g^*/G$. We denote subsets of $\t^*_+$ with a subscript $+$ and subsets of $\underline{\g}^*$ with an underline. Given a subset $X_+$ of $\t_+^*$, we denote the corresponding subset of $\underline{\g}^*$ by $\underline{X}$ and the corresponding invariant subset of $\g^*$ by $X$, and similarly when instead we are given a subset of $\underline{\g}^*$ or an invariant subset of $\g^*$. In particular, we denote the subset of $\t^*_+$ corresponding to the image of $J$ by $\Delta_+$ (which by the non-abelian convexity theorem is often a convex polytope \cite{At,GS1,Kir} or, more generally, a convex locally polyhedral set \cite{HiNePl}). Using these notational conventions, our theorem on the equivariant constructibility of the momentum map can be stated as follows.
\begin{Thm}[Constructibility]\label{thm:stratification:image:momentum:map} The space $\Delta_+$ admits a natural integral affine 
stratification ${\sS}_\mathrm{Ham}(\Delta_+)$ with the property that, for each stratum $\sigma_+\in {\sS}_\mathrm{Ham}(\Delta_+)$, the restriction
\begin{equation}\label{eqn:restrictionstratum:mommap} J:J^{-1}(\sigma)\to \sigma
\end{equation} is a $G$-equivariantly locally trivial fibration (in the sense made precise in Theorem \ref{thm:localtriviality}). It is uniquely 
determined by the property that each stratum is an integral affine submanifold of $(\t^*,\Lambda^*_T)$ and that at any point $x\in \Delta_+$ the tangent space to the stratum $\sigma_+$ through $x$ is given by
\begin{equation}\label{eqn:definingproperty:stratofimage}
T_x\sigma_+=\bigcap_{p\in J^{-1}(x)} (\g_p^0)^{G_x}\subset \t^*, 
\end{equation}  where 
\begin{itemize} 
\item $(\t^*,\Lambda^*_T)$ is the integral affine vector space with $\Lambda^*_T$ the character lattice of $T$, 
\item $\g_p$ is the isotropy Lie algebra at $p$ of the $G$-action on $S$,
\item $G_x$ is the isotropy group at $x$ of the coadjoint $G$-action and $\g_x$ is its Lie algebra,
\item $\g_p^0\subset \g_x^*$ is the annihilator of $\g_p$ in $\g_x$, 
\item $(\g_p^0)^{G_x}:=\g_p^0\cap (\g_x^*)^{G_x}$ is the intersection of $\g_p^0$ with the fixed point set of the coadjoint $G_x$-action (viewed as linear subspace of $\t^*$ via the identifications recalled in \S\ref{sec:background:Lie}).
\end{itemize} Moreover, each stratum of ${\sS}_\mathrm{Ham}(\Delta_+)$ is contained in an open face of $\t^*_+$.
\end{Thm}
This stratification is closely related to the canonical Hamiltonian stratification $\sS_\textrm{Ham}(S)$ studied in \cite{Mol1}. The strata of $\sS_\textrm{Ham}(S)$ are the connected components of the members of the partition given by the equivalence relation
\begin{equation*} p\sim q \iff (G_p,G_{J(p)})\text{ and }(G_q,G_{J(q)})\text{ are conjugate pairs of subgroups of }G.
\end{equation*} 
For each stratum $\Sigma\in \sS_\textrm{Ham}(S)$, the set $J(\Sigma)_+\subset \t^*_+$ corresponding to the image of $\Sigma$ under $J$ is an integral affine submanifold of $(\t^*,\Lambda_T)$ (by Proposition \ref{prop:images:canhamstrata:affineopen}). The defining property of $\sS_\mathrm{Ham}(\Delta_+)$ in Theorem \ref{thm:stratification:image:momentum:map} says that, near a point $x\in \Delta_+$, the stratum of $\sS_\mathrm{Ham}(\Delta_+)$ through $x$ is given by the intersection of those members of $\{J(\Sigma)_+\mid\Sigma\in \sS_\textrm{Ham}(S)\}$ that contain $x$. This is because any connected affine submanifold of $\t^*$ is an open subset of its affine hull and, for any $\Sigma$ intersecting $J^{-1}(x)$ and any $p\in \Sigma\cap J^{-1}(x) $, the affine hull of $J(\Sigma)_+$ is parallel to $(\g_p^0)^{G_x}$.\\

All of this simplifies when $G=T$ is a torus. In that case, 
\begin{itemize}\item $\sS_\textrm{Ham}(S)$ is simply the stratification of $S$ by orbit types of the $T$-action, meaning that its strata are the connected components of the members of the partition given by
\begin{equation*} p\sim q\iff T_p=T_q;
\end{equation*}
\item $\t^*_+=\t^*=\g^*$, $\Delta_+=\Delta$ and $J(\Sigma)_+=J(\Sigma)$ for each $\Sigma\in \sS_\textrm{Ham}(S)$;
\item $G_x=T$ each $x$, so that $(\g_p^0)^{G_x}=\g_p^0$ is just the annihilator of $\g_p$ in $\t^*$.
\end{itemize}
\begin{Ex} To illustrate the stratification $\sS_\mathrm{Ham}(\Delta)$, consider the Hamiltonian $\T^2$-action on $(S,\omega):=(\mathbb{C}P^1,\omega_\textrm{st})\times (\mathbb{C}P^2,3\,\omega_\textrm{st})$ given by
\begin{equation*} (\lambda_1,\lambda_2)\cdot ([w_0:w_1],[z_0:z_1:z_2])=([w_0:\lambda_1w_1],[z_0,\lambda_1z_1:\lambda_2z_2]),
\end{equation*} with the momentum map $J:\mathbb{C}P^1\times\mathbb{C}P^2\to \R^2$ given by 
\begin{equation*} J([w_0:w_1],[z_0:z_1:z_2])=\left(\frac{|w_1|^2}{||w||^2}+3\frac{|z_1|^2}{||z||^2},3\frac{|z_2|^2}{||z||^2}\right).
\end{equation*}
The image $\Delta$ of $J$ with the stratification $\mathcal{S}_\mathrm{Ham}(\Delta)$ is as depicted below. 
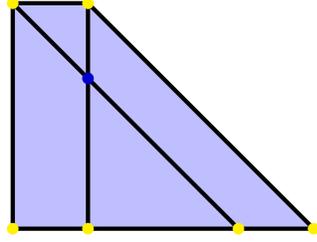
\begin{center}
\begin{tikzpicture} 
\filldraw[blue!25!white] (0,0) -- (0,3) -- (1,3) -- (4,0) -- cycle ;

\draw[ultra thick,black] (0,0) -- (0,3) ;
\draw[ultra thick,black] (0,3) -- (1,3) ;
\draw[ultra thick,black] (1,3) -- (4,0) ;
\draw[ultra thick,black] (0,0) -- (4,0) ;
\draw[ultra thick, black] (1,3) -- (1,0) ;
\draw[ultra thick, black] (0,3) -- (3,0) ;
\filldraw[yellow] (0,0) circle (2pt) ;
\filldraw[yellow] (0,3) circle (2pt) ;
\filldraw[yellow] (1,3) circle (2pt) ;
\filldraw[yellow] (4,0) circle (2pt) ;
\filldraw[yellow] (1,0) circle (2pt) ;
\filldraw[yellow] (3,0) circle (2pt) ;
\filldraw[blue!80!black] (1,2) circle (2pt);
\end{tikzpicture}
\end{center} To see this, note first that the action has six fixed points. Their images (the yellow dots) are zero-dimensional strata and $\Delta$ is their convex hull. To compute the rest of the stratification we ought to find the images of all the orbit type strata with isotropy groups of dimension one in $\T^2$. These turn out to be the subgroups $\mathbb{S}^1\times \{\ast\}$, $\{\ast\}\times\mathbb{S}^1$ and the diagonally embedded $\mathbb{S}^1$. For each of these there are two strata with that isotropy group and their images are the open line segments
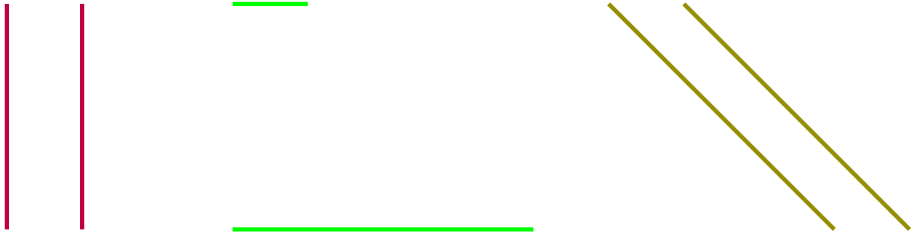
\begin{center}
\begin{tikzpicture} 

\draw[ultra thick,purple] (0,0) -- (0,3) 
;
\draw[ultra thick,purple] (1,0) -- (1,3) ;

\draw[ultra thick,green] (3,3) -- (4,3) ;
\draw[ultra thick,green] (3,0) -- (7,0);

\draw[ultra thick,olive] (9,3) -- (12,0) ;
\draw[ultra thick,olive] (8,3) -- (11,0) ;

\end{tikzpicture}
\end{center}

Following the recipe of (locally) intersecting these we find that there is one more zero-dimensional stratum (the blue dot) and there are ten one-dimensional strata (the black open line segments). The four connected open regions that remain (coloured light blue) are the two-dimensional strata.
\end{Ex}

It is a general fact that, if $G=T$ acts effectively, then the open strata of $\sS_\mathrm{Ham}(\Delta)$ are the connected components of the set of regular values of $J$ in $\Delta$. Folklore says that, if we further suppose that $S$ is compact and connected, then the closures of these components are convex polytopes, any two of which intersect in a common face. It turns out that, in that case, the strata of $\sS_\mathrm{Ham}(\Delta)$ are exactly the open faces of these polytopes. So, $\sS_\mathrm{Ham}(\Delta)$ is then entirely determined by the convex geometry of the set of regular values of $J$ in $\Delta$. A proof of this will be not be given here, but in a separate note \cite{MoSe}. We are hopeful that such a convex geometric description can also be given for actions of arbitrary compact connected Lie groups, but a proof of this is still work in progress. \\

We now return to the more general setting in which $G$ is any compact connected Lie group. 
Recall that the reduced spaces are the topological quotient spaces 
\begin{equation*} \underline{S}_\L:=J^{-1}(\L)/G,\quad \L\in \underline{\g^*}. 
\end{equation*} These are the fibers of the transverse momentum map
\begin{equation*} \underline{J}:\underline{S}\to \underline{\g^*},
\end{equation*} by which we mean the map between the orbit spaces $\underline{S}:=S/G$ and $\underline{\g^*}$ induced by $J$. Denote by 
\begin{equation*} \underline{J}_+:\underline{S}\to \t^*_+
\end{equation*} its composition with the canonical homeomorphism between $\underline{\g^*}$ and $\t^*_+$. Theorem \ref{thm:stratification:image:momentum:map} has the following immediate corollary, which in particular says that the symplectic reduced spaces $\underline{S}_{\L_x}$ do not change topologically as we vary $x$ within a stratum of $\sS_\mathrm{Ham}(\Delta_+)$. 
\begin{Cor}\label{cor:thm:stratification:image:transverse:momentum:map} For each $\sigma_+\in {\sS}_\mathrm{Ham}(\Delta_+)$, the restricted map
\begin{equation}\label{eqn:restrictionstratum:transversemommap} \underline{J}_+:\underline{J}_+^{-1}(\sigma_+)\to \sigma_+
\end{equation}  is a topologically locally trivial fibration.
\end{Cor}
Actually, it follows from Theorem \ref{thm:stratification:image:momentum:map} that the map \eqref{eqn:restrictionstratum:transversemommap} is not only locally trivial topologically, but also in a stratified sense (see Remark \ref{rem:trivstratpreserving}). Moreover, the local trivializations that we construct are not only smooth stratum-wise, but are in fact induced by smooth maps on the level of $S$ (see Theorem \ref{thm:localtriviality}). \\

Theorem \ref{thm:stratification:image:momentum:map} and the above corollary further imply the following.
\begin{Cor}\label{cor:constructibility:pushforwards} Let $i\in \N$.
\begin{itemize}
\item[(a)] The derived push-forward $R^i(\underline{J}_+^{EG})_*(\underline{\R})$ of the constant sheaf $\underline{\R}$ along 
\begin{equation*} \underline{J}_+^{EG}:=\underline{J}_+\circ \mathrm{pr}_{\underline{S}}:EG\times_GS\to \Delta_+ 
\end{equation*} is constructible with respect to ${\sS}_\mathrm{Ham}(\Delta_+)$. 
\item[(b)] The same holds for the derived push-forward $R^i(\underline{J}_+)_*(\underline{\R})$ of the constant sheaf $\underline{\R}$ along 
\begin{equation*} 
\underline{J}_+:\underline{S}\to \Delta_+.
\end{equation*} 
\end{itemize}
\end{Cor}
Here, by constructibility of a sheaf of vector spaces with respect to a stratification we mean that its restriction to each stratum is locally constant and that its stalks are finite-dimensional. For $x\in \underline{\Delta}$, the stalk at $x$ of the sheaf in (a) is the equivariant cohomology $H^i_{G}(S_{\L_x},\R)$ of the $G$-space $S_{\L_x}:=J^{-1}(\L_x)$, whilst the stalk at $x$ of that in (b) is the singular cohomology $H^i(\underline{S}_{\L_x},\R)$. So, the above corollary is a statement about how the cohomology and the equivariant cohomology of the reduced spaces vary with $x\in \underline{\Delta}$. In particular, it implies that these are isomorphic for all $x$ within a stratum $\underline{\sigma}$. Such an isomorphism between two points $x,y\in \underline{\sigma}$ is canonically associated to any path $\gamma$ starting at $x$ and ending at $y$. Constructibility further implies that, for any path $\gamma$ that starts at $x\in \underline{\sigma}$ and remains in $\underline{\sigma}$ for all time except for at its end-point $y$, there are associated maps
\begin{equation}\label{eqn:specialization:maps} H^i(\underline{S}_{\L_y},\R)\to H^i(\underline{S}_{\L_x},\R), \quad H^i_G({S}_{\L_y},\R)\to H^i_G({S}_{\L_x},\R).
\end{equation}  
Whilst we leave a more extensive study of these maps for future work, we do make the following observation, which is reminiscent of the local invariant cycle theorem for complex algebraic maps.



\begin{Prop}\label{prop:specializationmaps} In the setting of Theorem \ref{thm:stratification:image:momentum:map}, let $x,y\in \Delta_+$ be regular values of $J$, connected by a path $\gamma$ from $x$ to $y$ as above, with $x$ in the interior of $\t^*_+$. Then $H^\bullet(\underline{S}_{\L_x},\R)$ is naturally acted on by the Weyl group $\mathbb{W}_y$ of $G_y$ and the  maps \eqref{eqn:specialization:maps} associated to $\gamma$ (which coincide in this case) are isomorphisms onto the set of $\mathbb{W}_y$-fixed points of $H^\bullet(\underline{S}_{\L_x},\R)$. 
\end{Prop}
Next, we turn to our second main result, which is an extension of the so-called Duistermaat-Heckman linear variation theorem \cite{DuHe} to arbitrary (possibly singular) reduced spaces. 
\begin{Thm}[Linear variation]\label{thm:linvar:main} Let $\sigma_+$ be a stratum of the stratification ${\sS}_\mathrm{Ham}(\Delta_+)$ in Theorem \ref{thm:stratification:image:momentum:map}. When compared via the parallel transport of the Gauss-Manin connection of the fibration (\ref{eqn:restrictionstratum:transversemommap}), the cohomology class $\varpi_x\in H^2(\underline{S}_{\L_x})$ of the symplectic structure on the reduced space $\underline{S}_{\L_x}$ at $\L_x$ (defined as in \cite{Sj2}; see Section \ref{subsec:deRhammodel:1} as well) varies linearly with $x\in \sigma_+$. 
\end{Thm}
Duistermaat and Heckman proved linear variation for the cohomology class of the symplectic form on the reduced spaces at regular values of $J$ in the interior of $\t^*_+$. When $J$ is so that the set $\t^*_{\mu-\mathrm{reg}}$ of such values is non-empty, it is a dense open subset of $\t^*$ and it coincides with the union of all open strata of ${\sS}_\mathrm{Ham}(\Delta_+)$. Therefore, Duistermaat and Heckman's theorem is recovered by applying Theorem \ref{thm:linvar:main} to the open strata of ${\sS}_\mathrm{Ham}(\Delta_+)$. \\

An important consequence of the linear variation theorem in \cite{DuHe} is that the symplectic volumes $\mathrm{Vol}(\underline{S}_{\L_x}):=\mathrm{Vol}(\underline{S}_x,\varpi_x)$ of the reduced spaces at values $x\in \t^*_{\mu-\mathrm{reg}}$ are polynomial as function of $x$ on each connected component of $\t^*_{\mu-\mathrm{reg}}$, of degree at most half the dimension of the reduced spaces at regular values. Theorem \ref{thm:linvar:main} leads to the following extension of this for symplectic volumes of reduced spaces at arbitrary values (which were shown to be finite in \cite{Sj2}).
\begin{Cor}\label{cor:sympvolpolynomial}
For each stratum $\sigma_+$ of the stratification ${\sS}_\mathrm{Ham}(\Delta_+)$ in Theorem \ref{thm:stratification:image:momentum:map}, the symplectic volume $\mathrm{Vol}(\underline{S}_{\L_x})$ of the reduced space at $x$ is polynomial as function of $x\in \sigma_+$, of degree at most half the dimension of these reduced spaces.
\end{Cor}

\textbf{\underline{Outline:}} The majority of \S 1 is devoted to the proof of Theorem \ref{thm:stratification:image:momentum:map}. First, in \S 1.1 we give a construction (which could be of independent interest) with input a certain type of collection of (not necessarily disjoint) affine submanifolds of a (finite-dimensional real) vector space and output an affine stratification of the subset covered by the given collection. In \S 1.2 we apply this to show that there is a (necessarily unique) integral affine stratification of $\Delta_+$ with the defining property in Theorem \ref{thm:stratification:image:momentum:map} and in \S 1.3 we complete the proof of Theorem \ref{thm:stratification:image:momentum:map} by proving that $J$ has the desired local triviality property with respect to this stratification. After that, we prove Corollary \ref{cor:constructibility:pushforwards} and Proposition \ref{prop:specializationmaps} in \S 1.4.\\

In \S 2.1 we recall Sjamaar's de Rham model for the cohomology of the symplectic reduced spaces. In addition to this, we give a model for their homology build out of ``smooth" simplices, which allows us to express the canonical pairing between Sjamaar's de Rham cohomology and singular homology by integrating differential forms over such simplices. This is used in the proof of Theorem \ref{thm:linvar:main}, which is given in \S 2.2 along with that of Corollary \ref{cor:sympvolpolynomial}. \\

Finally, the extensions of our main results to actions of arbitrary (possibly disconnected) compact Lie groups are stated and proved in \S 3.1 and those for actions of proper quasi-symplectic groupoids are stated and proved in \S 3.2.\\


\textbf{\underline{Acknowledgements:}} I am grateful to Ana B\u{a}libanu for interesting discussions related to this project. This work was supported partly by the Max Planck Institute for Mathematics and partly by NSERC Discovery Grant RGPIN-2024-05764. \\

\textbf{\underline{Conventions:}} Throughout, we require smooth manifolds to be both Hausdorff and second countable and we require the same for both the base and the space of arrows of a Lie groupoid. By a submanifold we mean an embedded submanifold. Our conventions and terminology on stratifications are as in \cite{Mol1}. 

\section{Constructibility}
\subsection{Affine partitions and stratifications} 
\subsubsection{A construction of affine stratifications of subsets of vector spaces} The construction of the affine stratification in Theorem \ref{thm:stratification:image:momentum:map} will be based on a more general mechanism, given by the following proposition. 
\begin{prop}\label{prop:piecewise-affine:stratification} 
Let $V$ be a real finite-dimensional vector space and $X$ a subset of $V$. Given a piecewise-affine cover $\A=\{P_i\mid i\in I\}$ of a subset $X$ of $V$, there is a unique partition $\sS(\A)$ of $X$ into connected affine submanifolds of $V$ with the property that for each $\Sigma\in \sS(\A)$ and each $x\in \Sigma$:
\begin{equation}\label{eqn:tngtsp:cond} T_x\Sigma=\bigcap_{i\in I_x}T_xP_i, \quad I_x:=\{i\in I\mid x\in P_i\}.
\end{equation}
Moreover, if $X$ is locally closed in $V$, then $\sS(\A)$ is a Whitney $(b)$-regular stratification of $X$.
\end{prop}
Here, piecewise-affine cover is meant in the following sense. 
\begin{defi}\label{def:affineopen:affinecover}
We call a non-empty subset $P$ of $V$ \textbf{affine-open} if it is open in its affine hull. If in addition $P\subset X$, then we say that it is \textbf{affine-open in $X$}. By a \textbf{piecewise-affine cover} of $X$ we mean a locally finite cover $\A$ of $X$ with the following properties:
\begin{itemize}
\item[(i)] each member of $\A$ is affine-open in $X$,
\item[(ii)] the closure in $X$ of any member of $\A$ is a union of members of $\A$. 
\end{itemize} 
\end{defi}
\begin{rem}\label{rem:affopenvsaffsubm} Affine-opens are particular types of affine submanifolds of $V$. Whilst general affine submanifolds of $V$ need not be affine-open, all connected affine submanifolds of $V$ are affine-open. To see this, note that if $P$ is an affine submanifold of $V$, then for each affine subspace $A$ of $V$ the subset $P_A:=\{x\in P\mid x\in A\textrm{ and }T_xP=T_xA\}$ is open in $P$. So, since $P$ is partitioned by such subsets, when connected it equals $P_A$ for one $A$ and, hence, it is affine-open with affine hull $A$. 
\end{rem}
The main point of this subsection will be to prove Proposition \ref{prop:piecewise-affine:stratification}. First, let us give some examples.
\begin{ex} To illustrate the relevance of condition (ii) in Definition \ref{def:affineopen:affinecover} for Proposition \ref{prop:piecewise-affine:stratification} to hold, note that the cover of $X:=\R^2$ consisting of the three affine-open subspaces $\{x<0,y=0\}$, $\{x>0,y=0\}$ and $\R^2$ does not admit a partition as in the proposition.
\end{ex}
\begin{ex} If $X$ is not locally closed, then $\sS(\A)$ need not be locally finite (and, hence, it need not be a stratification). Indeed, if $V:=\R^2$ and
\begin{equation*} X:=[0,\infty[\times]-1,1[\,-\bigsqcup_{n\in \N}\{\tfrac{1}{n}\}\times ]0,1[,
\end{equation*} then the piecewise-affine cover $\A$ of $X$ consisting of the affine-opens
\begin{equation*} 
\{(0,0)\},\quad \{0\}\times]0,1[,\quad \{0\}\times]-1,0[,\quad ]0,\infty[\times]-1,1[\,-\bigsqcup_{n\in \N}\{\tfrac{1}{n}\}\times [0,1[,\quad ]0,\infty[\times \{0\}
\end{equation*}
is so that any open around a point in $\{0\}\times]0,1[$ intersects infinitely many members of $\sS(\A)$. 
\end{ex}
We now turn to the proof of Proposition \ref{prop:piecewise-affine:stratification}. Given a partition $\P$ of a subset of $V$, let $\D(\P)$ denote the discrete subbundle of $TV$ with fiber at $x$ the tangent space $T_xP$ to the member $P$ of $\P$ containing $x$. 
\begin{lemma}\label{lemma:affinestrat:refinementcrittangentspaces} Suppose that $X\subset Y\subset V$ and let $\P$ and $\Q$ be partitions of $X$ and $Y$, respectively, into connected affine submanifolds of $V$. If $\D(\P)\subset \D(\Q)$, then each member of $\P$ is contained in a member of $\Q$. 
\end{lemma}
\begin{proof} If $P\in \P$ and $Q\in \Q$, then since $P$ and $Q$ are affine submanifolds, the germs of $P$ and $Q$ in $V$ at any $x\in P\cap Q$ coincide with those of $x+T_xP$ and $x+T_xQ$, respectively. So, since in addition $T_xP=\D(\P)_x\subset \D(\Q)_x=T_xQ$ for any such $x$, $P\cap Q$ is open in $P$. Hence, by connectedness, each $P\in \P$ intersects only one $Q\in \Q$, which therefore contains $P$. 
\end{proof}
It follows from this that a partition $\P$ of $X$ into connected affine submanifolds is entirely determined by the discrete subbundle $\D(\P)$ of $TV$. Lemma \ref{lem:affine-integrability} below provides a useful converse to this.
\begin{lemma}\label{lem:affine-integrability} Let $X$ be a subset of $V$ and let $\D$ be a discrete wide subbundle of $TV\vert_X$ which is affine-integrable, in the sense that for every $x\in X$ there is an affine-open $P$ in $X$ such that $x\in P$ and $T_yP=\D_y$ for all $y\in P$. Then there is a (necessarily unique) partition $\P$ of $X$ into connected affine-opens with the property that $T_xP=\D_x$ for all $P\in \P$ and $x\in P$. 
\end{lemma}
\begin{proof} For $x\in X$, let $P_x$ denote the union of all connected affine-opens $P$ in $X$ that contain $x$ and are such that $T_yP=\D_y$ for all $y\in P$. This is itself an affine-open in $X$ with these properties and, in fact, it is the greatest such affine-open. Using this it is readily seen that $\{P_x\mid x\in X\}$ is the desired partition.  
\end{proof}
For the proof of Proposition \ref{prop:piecewise-affine:stratification} we will use one more lemma, which is straightforward to prove. 
\begin{lemma}\label{lem:intersection:aff-opens} Let $P_1,...,P_n$ be affine-opens in $V$. If non-empty, then the intersection $P_1\cap ... \cap P_n$ is also affine-open in $V$, with affine hull equal to the intersection $\Af(P_1)\cap...\cap \Af(P_n)$ of the individual affine-hulls.  
\end{lemma}
\begin{proof}[Proof of Proposition \ref{prop:piecewise-affine:stratification}] Given a non-empty finite subset $J\subset I$, consider:
\begin{equation}\label{eqn:affopen:subindex} P_J:=\bigcap_{j\in J} P_j,\quad P^{-}_{J}:=P_J-\bigcup_{i\in I(J)} P_i,
\end{equation} with $I(J)\subset I$ the subset of those $i\in I$ for which $P_i\cap P_J$ is not open in $P_J$. We claim that:
\begin{itemize}\item[(a)] $P^{-}_J$ is open in $P_J$,
\item[(b)] for every $x\in P^{-}_J$ it holds that $T_xP^{-}_J=\bigcap_{i\in I_x}T_xP_i$, with $I_x$ as in (\ref{eqn:tngtsp:cond}).
\end{itemize} Since $I_x$ is non-empty ($\A$ being a cover of $X$) and finite (by local finiteness of $\A$) and $x\in P^{-}_{I_x}$, for each $x\in X$, this would show that the assumptions of Lemma \ref{lem:affine-integrability} are satisfied for the subbundle $\D:=\D(\A)$ with fiber at $x$ given by the right-hand side of (\ref{eqn:tngtsp:cond}). So, to prove the existence and uniqueness of the  partition $\sS(\A)$, it suffices to prove claims (a) and (b). For (a), we will show that the complement:
\begin{equation}\label{eqn:complPJ} P_J-P^{-}_J=\bigcup_{i\in I(J)}P_i\cap P_J
\end{equation} is closed in $P_J$. By local finiteness of $\A$, the closure of this complement is the union of closures:
\begin{equation}\label{eqn:complPJclos} \mathrm{Cl}_{P_J}(P_J-P^{-}_J)=\bigcup_{i\in I(J)}\mathrm{Cl}_{P_J}(P_i\cap P_J).
\end{equation} Let $i\in I(J)$. By condition (ii), there is a subset $J_i\subset I$ such that:
\begin{equation*} \mathrm{Cl}_X(P_i)=\bigcup_{k\in J_i} P_k. 
\end{equation*} Since the closure $\mathrm{Cl}_X(P_i)$ is contained in the affine hull $\mathrm{Af}(P_i)$, and $\mathrm{Af}(P_i)\cap \mathrm{Af}(P_J)=\mathrm{Af}(P_i\cap P_J)$ (Lemma \ref{lem:intersection:aff-opens}), it holds that $\mathrm{Af}(\mathrm{Cl}_X(P_i)\cap P_J)=\mathrm{Af}(P_i\cap P_J)$. So, if $k\in J_i$ is such that $P_k\cap P_J\neq \emptyset$, then: 
\begin{equation*} \mathrm{Af}(P_k\cap P_J)\subset \mathrm{Af}(P_i\cap P_J)\subset \mathrm{Af}(P_J).
\end{equation*} If moreover $P_k\cap P_J$ were open in $P_J$, then it would hold that $\mathrm{Af}(P_k\cap P_J)=\mathrm{Af}(P_J)$, so that $\mathrm{Af}(P_i\cap P_J)=\mathrm{Af}(P_J)$, which would force $P_i\cap P_J$ to be open in $P_J$ as well. Therefore, $P_k\cap P_J$ is not open in $P_J$. From this we conclude that:
\begin{equation*} \mathrm{Cl}_{P_J}(P_i\cap P_J)\subset \mathrm{Cl}_X(P_i)\cap P_J\subset P_J-P_J^-.
\end{equation*} This being true for all $i\in I(J)$, it follows that (\ref{eqn:complPJclos}) is contained in (\ref{eqn:complPJ}), so that (\ref{eqn:complPJ}) is closed in $P_J$. This proves claim (a). For (b), let $x\in P^{-}_J$. Then $P_i\cap P_J$ is open in $P_J$ for each $i\in I_x$. Moreover, $J\subset I_x$ and so:
\begin{equation*} P_{I_x}=\bigcap_{i\in I_x}P_i\cap P_J.
\end{equation*} Since $I_x$ is finite it follows that $P_{I_x}$ is open in $P_J$ as well. By first using claim (a), then using that $P_{I_x}$ is open in $P_J$ and finally using Lemma \ref{lem:intersection:aff-opens}, we conclude that:
\begin{equation*} T_xP^{-}_J=T_xP_J=T_xP_{I_x}=\bigcap_{i\in I_x} T_xP_i,
\end{equation*} which proves claim (b). \\

To prove that $\sS(\A)$ is a Whitney (b)-regular stratification of $X$ when $X$ is locally closed in $V$, we will use the following fact (see, e.g., \cite{Mat1,GWPL}): given a locally closed subset $X$ of a smooth manifold $M$ and a locally finite partition $\P$ of $X$ into smooth submanifolds of $M$ that is Whitney (b)-regular (in the sense that each any member of $\P$ is Whitney (b)-regular over any other member), the partition $\P^c$ consisting of connected components of the members of $\P$ is a Whitney (b)-regular stratification of $X$ (meaning, in particular, that $\P^c$ is locally finite and satisfies the frontier condition). In view of this, it suffices to construct a partition $\P:=\P(\A)$ with these properties, such that $\P(\A)^c=\sS(\A)$. For this, consider as before the affine-integrable discrete subbundle $\D:=\D(\A)$ of $TV\vert_X$ with fiber at $x$ given by (\ref{eqn:tngtsp:cond}). Given $x\in X$, let $\widehat{P}_x$ denote the union of all (not necessarily connected) affine-opens $P$ in $X$ that contain $x$ and are such that $T_yP=\D_y$ for all $y\in P$. This is itself an affine-open with these properties. In fact, it is the greatest such affine-open. Using this it is readily seen that the collection $\P(\A):=\{\widehat{P}_x\mid x\in X\}$ is a partition of $X$. It is clear that $\sS(\A)=\P(\A)^c$. So, we are left to prove that $\P(\A)$ is locally finite and Whitney regular. \\

For local finiteness, let $x\in X$. Since $\A$ is locally finite, there is an open neighbourhood $U$ of $x$ that intersects $P_i$ for only finitely many $i\in I$. For such $U$, there are only finitely many $J\subset I$ that are of the form $J=I_x$ for some $x\in U$. If $J=I_x$ for some $x\in U$, then $x\in P_J^-$ and (due to claim (b) above) $P_J^-$ is contained in a single member of $\P(\A)$. Therefore, any $P\in \P(\A)$ that intersects $U$ contains $P^-_J$ for some such $J$ (take $J=I_x$ for any $x\in P\cap U$), and so $U$ intersects only finitely members of $\P(\A)$.\\

For Whitney regularity, let $P,Q\in \P(\A)$ and suppose that $(x_n)$ and $(y_n)$ are sequences in $P$ and $Q$ respectively, both converging to some $x\in P$, such that $T_{y_n}Q$ converges in the Grassmannian of $\dim(Q)$-dimensional linear subspaces of $V$, say with limit $\tau$, and such that the sequence of lines $[x_n-y_n]$ converges in the projectivization of $V$, say with limit $\ell$. Then $x\in \mathrm{Af}(Q)$, since $\mathrm{Af}(Q)$ is closed in $V$. Since $Q$ is affine-open, its tangent space at any two points is the same subspace of $V$, and hence equal to $\tau$. So, $\mathrm{Af}(Q)=x+\tau$ and hence $y_n-x\in \tau$ for all $n$. We will further show that $T_xP\subset \tau$, so that $\mathrm{Af}(P)\subset \mathrm{Af}(Q)$ and hence $x_n-x\in \tau$ for all $n$. This would imply that $x_n-y_n=(x_n-x)-(y_n-x)\in \tau$ for all $n$, from which it is clear that $\ell\subset \tau$, as is to be shown. So, it remains to prove that $T_xP\subset \tau$. Note that, since $\A$ is locally finite, there is an $n$ such that $x\in \mathrm{Cl}_X(P_i)$ for all $i\in I_{y_n}$. It follows from property (ii) in Definition \ref{def:affineopen:affinecover} that, for any $i\in I_{y_n}$, there is a $j\in I_x$ such that $P_j\subset \mathrm{Cl}_X(P_i)$. Since $\mathrm{Cl}_X(P_i)\subset \mathrm{Af}(P_i)$, we conclude that 
\begin{equation*} T_xP=\D_x:=\bigcap_{j\in I_x}T_xP_j \subset \bigcap_{i\in I_{y_n}}T_{y_n}P_i=:\D_{y_n}=T_{y_n}Q=\tau, 
\end{equation*} which completes the proof.
\end{proof}

\subsubsection{From local to global for abstract affine stratifications} The stratification in Theorem \ref{thm:stratification:image:momentum:map} will be constructed by means of Proposition \ref{prop:piecewise-affine:stratification}. Extending Theorem \ref{thm:stratification:image:momentum:map} to actions of proper quasi-symplectic groupoids will involve affine stratifications of more abstract spaces than subsets of a vector space. For the purpose of that extension, the following will be used.

\begin{lemma}\label{lemma:localtoglobal:abstractaffstrat} Let $Y$ be a topological space with a partition $\P$ of $Y$ into affine manifolds equipped with the subspace topology. Further, let $X$ be a subspace of $Y$ and $\D\subset (\bigsqcup_{P\in \P}TP)\vert_X$ a discrete subbundle. Suppose that every point in $X$ admits an open neighbourhood $U$ in $X$ and a stratification ${\sS}_U$ of $U$ into affine submanifolds of the members of $\P$ such that $T_x\Sigma=\D_x$ for each $\Sigma\in {\sS}_U$ and $x\in \Sigma$. Then there is a (necessarily unique) stratification $\sS$ of $X$ into affine submanifolds of the members of $\P$ with the property that $T_x\Sigma=\D_x$ for each $\Sigma\in \sS$ and $x\in \Sigma$.
\end{lemma}
\begin{proof} Uniqueness follows from the observation that the germ of an affine submanifold at a given point in the submanifold is determined by its tangent space at that point, via an argument similar to that in the proof of Lemma \ref{lemma:affinestrat:refinementcrittangentspaces}. \\

For existence, suppose that $x\in X$ and let $P\in \P$ denote the member through $x$. Define $\Sigma_x$ to be the union of all connected affine submanifolds $A$ of $P$ containing $x$ and integrating $\D$, in the sense that $T_yA=\D_y$ for all $y\in A$. Then $\Sigma_x$ is itself such a submanifold of $P$. To see this, let $y\in \Sigma_x$ and consider an open neighbourhood $U$ of $x$ and a stratification ${\sS}_U$ as in the statement of the lemma. Since ${\sS}_U$ is locally finite, there is an open $V$ in $U$ around $y$ such that $y$ is contained in the closure in $U$ of each stratum of ${\sS}_U$ that intersects $V$. Each such stratum different from  the stratum $\Sigma_{U,y}$ through $y$ then has dimension greater than that of $\Sigma_{U,y}$, which equals $\dim(\D_y)$. Since $y\in \Sigma_x$, there is an $A_y$ as above containing both $x$ and $y$, and so $\dim(\D_y)=\dim(A_y)=\dim(\D_x)$. It follows from this that $\Sigma_x\cap V$ is an open in $\Sigma_{U,y}$, for if $A$ is a submanifold of $P$ as above, then $A\cap V$ must be contained in and open in $\Sigma_{U,y}$ since for each $w\in A\cap V$ the stratum of ${\sS}_U$ through $w$ has dimension $\dim(\D_w)=\dim(A)=\dim(\D_x)=\dim(\D_y)$. This shows that $\Sigma_x$ is indeed an affine submanifold of $P$ integrating $\D$. Moreover, it is connected, being a union of connected subspaces with non-empty intersection. \\

To conclude the proof of existence, note that the manifolds $\Sigma_x$ define a partition $\sS$ of $X$ (since for any $x,y\in X$ either $\Sigma_x=\Sigma_y$ or $\Sigma_x\cap \Sigma_y=\emptyset$) and that, to see that this is the desired stratification, it remains to show that it is locally finite and satisfies the frontier condition. By assumption, any point in $X$ admits an open neighbourhood $U$ and a stratification ${\sS}_U$ of $U$ into affine submanifolds of the members of $\P$ that integrate $\D$. The same argument as for the uniqueness part shows that the partition of connected components of the members of $\{\Sigma\cap U\mid \Sigma\in \sS\}$ coincides with ${\sS}_U$. This implies that $\sS$ is itself a locally finite and satisfies the frontier condition, in view of the fact that these properties can be verified locally (see, e.g., \cite[Lemma 2.2]{Mol1}).
\end{proof}

\subsection{The stratification of the transverse momentum map image} 
\subsubsection{Background on structure theory of compact and connected Lie groups}\label{sec:background:Lie} Here we collect some background and set notation on Lie theory (we refer the reader to \cite{DuKo,GuLeSt} for more details). Let $G$ be a compact Lie group with Lie algebra $\g$. Fix a maximal torus $T$ in $G$ with Lie algebra $\t$ and a choice of closed Weyl chamber $\t^*_+$ in $\t^*$. Given $x \in \g^*$, we denote by $G_x$ the isotropy group of the coadjoint $G$-action and by $\g_x$ its Lie algebra. Throughout, as per usual: 
\begin{itemize}
\item we identify $\t^*$ (via the decomposition $\g=\t\oplus [\g,\t]$) with the $T$-fixed point set $(\g^*)^T$ in $\g^*$ (which is the annihilator of $[\g,\t]$),
\item for each $x\in \t^*$, we identify $\g_x^*$ (via the $G_x$-invariant decomposition $\g^*=(G_x \cdot \t^*)\oplus \g_x^0$) with the linear subspace $G_x\cdot \t^*$ of $\g^*$ (which is indeed linear, for it is identified with the linear subspace $\g_x$ of $\g$ via any $G$-invariant inner product on $\g$). Here $\g_x^0$ denotes the annihilator of $\g_x$ in $\g^*$, which is equal to the tangent space $T_x\L$ to the coadjoint $G$-orbit through $x$. 
\end{itemize}
Suppose now that $G$ is connected. Then the following hold. 
\begin{itemize}\item Each coadjoint orbit intersects $\t^*_+$ in exactly one point, so that the inclusion $\t^*_+\hookrightarrow \g^*$ induces a homeomorphism:
\begin{equation}\label{eqn:sweeping:map} \t^*_+\xrightarrow{\sim} \g^*/G=:\underline{\g}^*.
\end{equation} 
\item The polyhedral cone $\t^*_+$ is naturally stratified by its open faces, whilst the orbit space $\g^*/G$ comes with a natural stratification $\sS_G(\underline{\g}^*)$ by connected components of $G$-orbit types. The homeomorphism (\ref{eqn:sweeping:map}) is compatible with these stratifications, in the sense that it restricts to a diffeomorphism between strata. 
\item Each open face $F$ of $\t^*_+$ is affine-open with affine hull the linear subspace $(\g_x^*)^{G_x}$ of $\t^*$, for any choice of $x\in F$. The isotropy group $G_x$ is the same for all $x\in F$. Therefore, it also denoted as $G_F$. The subspace $\Sigma_F:=G\cdot F$ of $\g^*$ is a submanifold of $\g^*$ 
and there is a $G$-equivariant diffeomorphism:
\begin{equation}\label{eqn:saturationopenface:diffeo} \varphi_F:G/G_F\times F\to \Sigma_F,\quad ([g],x)\mapsto g\cdot x. 
\end{equation} Note that $\underline{\Sigma}_F$ is a stratum of $\sS_G(\underline{\g}^*)$.
\end{itemize} 
\subsubsection{The symplectic cross-section theorem}\label{sec:sympcrossthm} In various parts of the proof of Theorem \ref{thm:stratification:image:momentum:map} we will reduce the proof of statements about symplectic reduced spaces at arbitrary values to the case of the symplectic reduced space at value zero. This will be done using the symplectic cross-section theorem, which we now recall. 
\begin{thm}\label{thm:symplectic:cross-section} Let $G$ be a compact Lie group, $T$ a maximal torus in $G$, $\t^*_+$ a closed Weyl chamber in $\t^*$ and $J:(S,\omega)\to \g^*$  a Hamiltonian $G$-space. Further, let $x\in \t^*_+$ and let $\mathfrak{S}$ be a slice through $x$ for the coadjoint $G$-action that is contained and open in $x+\g^*_x$. Then
\begin{itemize}
\item[(a)] $J$ is transverse to $\mathfrak{S}$ and $J^{-1}(\mathfrak{S})$ is a $G_x$-invariant symplectic submanifold of $(S,\omega)$,
\item[(b)] the restriction of $J$ to $(J^{-1}(\mathfrak{S}),\omega)$ is a momentum map for the restricted $G_x$-action,
\item[(c)] the map 
\begin{equation*} G\times_{G_x}J^{-1}(\mathfrak{S})\to S,\quad [g,p]\mapsto g\cdot p,
\end{equation*} is a $G$-equivariant embedding onto the open $J^{-1}(G\cdot \mathfrak{S})$ in $S$.
\end{itemize}
\end{thm}
The Hamiltonian $G_x$-space obtained by restricting to a slice $\mathfrak{S}$ as in this theorem will be called the \textbf{symplectic cross-section} of the Hamiltonian $G$-space at $\mathfrak{S}$. 
\begin{rem} A slice $\mathfrak{S}$ as in Theorem \ref{thm:symplectic:cross-section} always exists, because $G$ is compact and $\g^*=\g^*_x\oplus T_x\O$ for each $x\in \g^*$ (cf. the proof of \cite[Theorem 2.3.3]{DuKo}).
\end{rem}
\begin{rem} The symplectic cross-section at $\mathfrak{S}$ is Morita equivalent (in the sense of \cite{Mol1}) to the Hamiltonian $G$-space obtained by restricting $J$ to $G\cdot\mathfrak{S}$ (a $G$-invariant open in $\g^*$).
\end{rem}
\subsubsection{The canonical Hamiltonian stratification}\label{sec:canhamstrat} To construct the stratification in Theorem \ref{thm:stratification:image:momentum:map} we will use the \textbf{canonical Hamiltonian stratification} studied in \cite{Mol1}, which we denote as $\sS_\mathrm{Ham}(\underline{S})$. This is a stratification of the $G$-orbit space $\underline{S}$ that refines the orbit type stratification, with the property that the transverse momentum map $\underline{J}$ maps strata of $\sS_\mathrm{Ham}(\underline{S})$ into strata of the $G$-orbit type stratification of $\underline{\g}^*$. Moreover, $\underline{J}$ restricts to a smooth map of constant rank between each such pair of strata. As discussed in \cite[Section 2.2]{Mol1}, there are various partitions $\P$ of $\underline{S}$ (with possibly disconnected members) for which the induced partition $\P^\mathrm{c}$ formed by the connected components of the members of $\P$ coincides with $\sS_\mathrm{Ham}(\underline{S})$. One such partition is defined by the equivalence relation that declares the orbits through two points $p,q\in S$ to be equivalent when there is a $g\in G$ such that $gG_pg^{-1}=G_q$ and $gG_{J(p)}g^{-1}=G_{J(q)}$. Another such partition is defined by the equivalence relation that declares the orbits through two points $p,q\in S$ to be equivalent when $G_p$ is isomorphic to $G_q$ and $G_{J(p)}$ is isomorphic to $G_{J(q)}$. One way to characterize the smooth structures on the strata of $\sS_\mathrm{Ham}(\underline{S})$ is as follows: for any such stratum $\underline{\Sigma}$ the corresponding invariant subset $\Sigma$ of $S$ is a submanifold and the quotient map $q_S:S\to \underline{S}$ restricts to a submersion $\Sigma\to \underline{\Sigma}$.

\subsubsection{Construction of the stratification in Theorem \ref{thm:stratification:image:momentum:map}} In this section we show that, when $\underline{J}$ is proper as map onto its image, the images of the strata of $\underline{S}_\mathrm{Ham}(\underline{S})$ form a piecewise-affine cover of the image $\underline{\Delta}$ of $\underline{J}$. Applying Proposition \ref{prop:piecewise-affine:stratification} to this, we then obtain the stratification in Theorem \ref{thm:stratification:image:momentum:map}. \\

To begin with, we show the following. 
\begin{prop}\label{prop:images:canhamstrata:affineopen} Let $G$ be a compact and connected Lie group and $J:(S,\omega)\to \g^*$ a Hamiltonian $G$-space. Fix a maximal torus $T$ in $G$. For each stratum $\underline{\Sigma}\in {\sS}_\mathrm{Ham}(\underline{S})$, the subset $J(\Sigma)_+=J(\Sigma)\cap \t^*_+$ is affine-open in $\t^*$, with affine hull given by:
\begin{equation}\label{eqn:affhull:imofcanhamstrat} \emph{Af}(J(\Sigma)_+)=x+(\g_p^0)^{G_x}\subset \t^*
\end{equation} for any $x\in J(\Sigma)_+$ and $p\in \Sigma$ such that $J(p)=x$. Moreover, the restriction of $\underline{J}$ to $\underline{\Sigma}$, viewed as map into \eqref{eqn:affhull:imofcanhamstrat} via \eqref{eqn:sweeping:map}, is a submersion.
\end{prop}
\begin{proof} Let $\underline{\Sigma}\in {\sS}_\mathrm{Ham}(\underline{S})$ and let $\underline{J}_+:\underline{\Sigma}\to \t^*$ denote the restriction of $\underline{J}$ to $\underline{\Sigma}$, viewed as map into $\t^*$ via \eqref{eqn:sweeping:map}. Note that $J(\Sigma)_+=\underline{J}_+(\underline{\Sigma})$. We will show the following.
\begin{itemize}
\item[(i)] For each affine subspace $A$ of $\t^*$, the set 
\begin{equation*} \underline{\Sigma}_A:=\{\O\in \underline{\Sigma}\mid \underline{J}_+(\O)\in A\textrm{ and }\d(\underline{J}_+)_\O(T_\O\underline{\Sigma})=T_{\underline{J}_+(\O)}A\}
\end{equation*} is open in $\underline{\Sigma}$.
\item[(ii)] For each $\O\in \underline{\Sigma}$ it holds that $\d(\underline{J}_+)_{\O}(T_{\O}\underline{\Sigma})=(\g_p^0)^{G_x}$ as subspaces of $\t^*$, for any choice of $p\in \O$ such that $x:=J(p)\in \t^*_+$.
\end{itemize} It would follow from (i), by connectedness of $\underline{\Sigma}$, that there is a unique affine subspace $A$ of $\t^*$ containing $J(\Sigma)_+$ and such that $\underline{J}_+:\underline{\Sigma}\to A$ is a submersion. Since the image of a submersion is open, this means in particular that $J(\Sigma)_+$ is affine-open with affine hull $A$, which by (ii) can be described as in the statement of the proposition. So, it remains to prove (i) and (ii). For this, let $\O\in \underline{\Sigma}$ and fix a $p\in \O$ such that $x:=J(p)\in \t^*_+$. We ought to show that (ii) holds at $\O$ and that, if $A$ is such that $\O\in \underline{\Sigma}_A$, then there is an open neighbourhood of $\O$ in $\underline{\Sigma}$ that is contained in $\underline{\Sigma}_A$. Consider a slice $\Sigma$ at $x$ as in Theorem \ref{thm:symplectic:cross-section}. The inclusion $J^{-1}(\Sigma)$ into $S$ induces a homeomorphism of $J^{-1}(\Sigma)/G_x$ onto the open $\underline{J}^{-1}(\underline{G\cdot \mathfrak{S}})$ in $\underline{S}$. Moreover, as is readily verified using the second description of $\sS_\mathrm{Ham}$ recalled at the beginning of this subsection, this homeomorphism is an isomorphism of stratified spaces with respect to the canonical Hamiltonian stratification of the symplectic cross-section at $\mathfrak{S}$ and that of the restricted Hamiltonian $G$-action on $J^{-1}(G\cdot \mathfrak{S})$ (more conceptually, this is an instance of the invariance of $\sS_\mathrm{Ham}(\underline{S})$ under Hamiltonian Morita equivalence \cite[Remark 2.49]{Mol1}). Therefore, after passing to such a symplectic cross-section, we can assume that $J$ maps $\O$ to the origin in $\g^*$. Using the Marle-Guillemin-Sternberg normal form theorem around $\O_p$, the proof further reduces to the case in which the Hamiltonian $G$-space is the Marle-Guillemin-Sternberg local model. So, we can assume that $p=[1,0,0]\in G\times_H(\h^0\oplus V)$ and that $J$ is the restriction to a $G$-invariant open around $p$ of the map \begin{equation}\label{eqn:map:MGSlocalmodel}
G\times_H(\h^0\oplus V)\to \g^*,\quad [g,\alpha,v]\mapsto g\cdot \big(\alpha+(\p\circ J_V)(v)\big),
\end{equation} where
\begin{itemize}
\item $H$ is a closed subgroup of $G$, 
\item $\p:\h^*\to \g^*$ is an $H$-equivariant splitting of the canonical projection $\g^*\to \h^*$,
\item $J_V:V\to \h^*$ is the quadratic momentum map of a symplectic $H$-representation $(V,\omega_V)$, meaning that $\langle J_V(v), \xi\rangle=\tfrac{1}{2}\omega_V(\xi\cdot v,v)$ for $v\in V$ and $\xi \in \h$,
\item the action of $G$ on $G\times_H(\h^0\oplus V)$ is that induced by left translation of the $G$-factor. 
\end{itemize} 
Since $J$ maps $\O$ to the origin in $\g^*$, the stratum of $\sS_G(\underline{\g}^*)$ containing $\underline{J}(\underline{\Sigma})$ is $(\g^*)^G$ and the derivative of $\underline{J}_+$ at $\O$ is given by the composition
\begin{equation}\label{eqn:comp:pfprop:canhamstrata:affineopen}
T_{\O_p}\underline{\Sigma}\xrightarrow{\d(\underline{J}\vert_{\underline{\Sigma}})} (\g^*)^G\hookrightarrow \t^*.
\end{equation} By \cite[Proposition 2.56 and Lemma 2.57]{Mol1}, in an open neighourhood of $\O$ the stratum $\underline{\Sigma}$ coincides (as smooth manifold) with the vector space $(\h^0)^G\oplus V^H$ (which is a subspace of $(\h^0\oplus V)/H$). Moreover, $\underline{J}$ restricts to the map $\underline{\Sigma}\to (\g^*)^G$ given by the restriction of $(\h^0)^G\oplus V^H\to (\g^*)^G$, $(\alpha,v)\mapsto \alpha$. Since $H=G_p$, this shows that the derivative of $\underline{J}_+$ at $\O$ has image $(\g_p^0)^G$ and that, if $A$ is such that $\O\in \underline{\Sigma}_A$, then $A=(\g_p^0)^G$ and an open neighbourhood of $\O$ in $\underline{\Sigma}$ is contained in $\underline{\Sigma}_A$. This completes the proof. 
\end{proof} 
The affine-opens in Proposition \ref{prop:images:canhamstrata:affineopen} are, in fact, integral affine submanifolds of $(\t^*,\Lambda_T^*)$. This follows from the following lemma. 
\begin{lemma}\label{lemma:rationality:Lietheoreticsubspaces} Let $G$ be a compact connected Lie group and $T$ a maximal torus in $G$. 
\begin{itemize}
\item[(a)] $(\g^*)^G$ is a rational linear subspace of $(\t^*,\Lambda^*_T)$,
\item[(b)] for any closed subgroup $H$ of $G$, $\h^0\cap (\g^*)^G$ is a rational linear subspace of $(\t^*,\Lambda^*_T)$ as well. 
\end{itemize}
\end{lemma}
\begin{proof} Recall that $\g=Z(\g)\oplus \g^\mathrm{ss}$, with $Z(\g)=\g^G$ the center of $\g$ and $\g^\mathrm{ss}:=[\g,\g]$ the semisimple part of $\g$. Since $(\g^*)^G$ is the annihilator of $\g^\mathrm{ss}$, this composition gives an identification between $(\g^*)^G$ and $Z(\g)^*$. Via this, the dual of the inclusion map $(\g^*)^G\hookrightarrow \t^*$ is identified with the projection $\t=Z(\g)\oplus (\g^\mathrm{ss}\cap \t)\to Z(\g)$. The fact that $G^\mathrm{ss}\cap T$ is a torus (in fact, it is a maximal torus in $G^\mathrm{ss}$), implies that the image of $\Lambda_T$ under this projection is $\exp_{Z(G)}^{-1}(G^\mathrm{ss})$, or in other words, it is the kernel of the exponential map of the torus $G/G^\mathrm{ss}$. This image being discrete, we conclude that part (a) holds and that $(\g^*)^G\cap \Lambda^*_T$ is dual to the lattice $\ker(\exp_{G/G^{\mathrm{ss}}})$. Turning to part (b), suppose that $H$ is a closed subgroup of $G$. Denote by $H^0$ the identity component of $H$. Since 
\begin{equation*}
\ker(\exp_{G/G^{\mathrm{ss}}})\cap \big((\h+\g^{\mathrm{ss}})/\g^{\mathrm{ss}}\big)
\end{equation*} is the kernel of the exponential map of the torus $H^0/H^0\cap G^{\mathrm{ss}}$, it is a full rank lattice in $(\h+\g^{\mathrm{ss}})/\g^{\mathrm{ss}}$. So, it also holds that 
\begin{equation*}
\ker(\exp_{G/G^{\mathrm{ss}}})^*\cap \big((\h+\g^{\mathrm{ss}})/\g^{\mathrm{ss}}\big)^0
\end{equation*} is of full rank in the annihilator $\big((\h+\g^{\mathrm{ss}})/\g^{\mathrm{ss}}\big)^0\subset (\g/\g^{\mathrm{ss}})^*$. Since the canonical isomorphism between $(\g^*)^G=(\g^\mathrm{ss})^0$ and $(\g/\g^{\mathrm{ss}})^*$ identifies $(\g^*)^G\cap \Lambda^*_T$ with $\ker(\exp_{G/G^{\mathrm{ss}}})^*$ and $\h^0\cap (\g^*)^G$ with the annihilator of $(\h+\g^{\mathrm{ss}})/\g^{\mathrm{ss}}$,  part (b) follows.
\end{proof}
Next, we turn to the construction of the stratification in Theorem \ref{thm:stratification:image:momentum:map} that we outlined above. 
\begin{proof}[Proof of Theorem \ref{thm:stratification:image:momentum:map}, except local triviality] Consider the cover $\A_\mathrm{Ham}$ of $\underline{\Delta}$ given by
\begin{equation*} \A_\mathrm{Ham}:=\{\underline{J}_+(\underline{\Sigma})\mid \underline{\Sigma}\in {\sS}_\mathrm{Ham}(\underline{S})\}.
\end{equation*} To obtain the desired stratification $\sS_\mathrm{Ham}(\Delta_+)$, we will apply Proposition \ref{prop:piecewise-affine:stratification} to $\A_\mathrm{Ham}$. To verify the relevant conditions, note first that, being continuous and proper as map into the Hausdorff and first countable space $\Delta_+$, $\underline{J}_+$ is closed as map into $\Delta_+$ (by \cite{Pa}) and its fibers are compact. Since $\sS_{\mathrm{Ham}}(\underline{S})$ is locally finite and the fibers of $\underline{J}_+$ are compact, for every $x\in \Delta_+$ there is an open neighbourhood $\underline{U}$ of $\underline{J}_+^{-1}(x)$ in $\underline{S}$ that intersects finitely many strata of ${\sS}_\mathrm{Ham}(\underline{S})$. Because $\underline{J}_+$ is closed as map into $\Delta_+$, there is an open $V$ around $x$ in $\Delta_+$ such that $\underline{J}^{-1}_+(V)\subset \underline{U}$. By construction, $V$ intersects only finitely many members of $\A_\mathrm{Ham}$, which shows that $\A_\mathrm{Ham}$ is locally finite. Secondly, because $\underline{J}_+$ is closed and continuous as map into $\Delta_+$, the image under $\underline{J}_+$ of the closure of $\underline{\Sigma}$ in $\underline{S}$ coincides with the closure in $\Delta_+$ of the image $\underline{J}_+(\underline{\Sigma})$, for any subset $\underline{\Sigma}$ of $\underline{S}$. So, it follows from the frontier condition for ${\sS}_\mathrm{Ham}(\underline{S})$ that the closure in $\Delta_+$ of each member of $\A_\mathrm{Ham}$ is a union of members of $\A_\mathrm{Ham}$. So, $\A_\mathrm{Ham}$ satisfies condition (ii) in Definition \ref{def:affineopen:affinecover}. From the fact that $\underline{J}_+$ is closed as map into its image and has compact fibers, it follows that its image $\Delta_+$ inherits local compactness from $\underline{S}$. Therefore, $\Delta_+$ is locally closed in $\t^*$. In view of this and Proposition \ref{prop:images:canhamstrata:affineopen}, $\A_\mathrm{Ham}$ indeed meets the conditions in Proposition \ref{prop:piecewise-affine:stratification} and we obtain a Whitney (b)-regular affine stratification of $\Delta_+$ that (by Lemma \ref{lemma:affinestrat:refinementcrittangentspaces}) refines the stratification of $\t^*_+$ by open faces and is uniquely determined by the property that its strata satisfy \eqref{eqn:definingproperty:stratofimage}. \\

It remains to verify that the strata of the resulting stratification are not only affine, but integral affine, meaning that the tangent spaces to the strata are rational linear subspaces of $(\t^*,\Lambda^*_T)$. This follows from Lemma \ref{lemma:rationality:Lietheoreticsubspaces}, since any intersection of rational linear subspaces is again rational. 
\end{proof}

\subsection{Local triviality} Having constructed the stratification $\sS_\mathrm{Ham}(\Delta_+)$, to prove Theorem \ref{thm:stratification:image:momentum:map} it remains to address the local triviality. More precisely, it remains to prove the following. 
\begin{thm}\label{thm:localtriviality}
In the setting of Theorem \ref{thm:stratification:image:momentum:map}, let $\sigma_+\in {\sS}_\mathrm{Ham}(\Delta_+)$ and let $\sigma$ be the corresponding $G$-invariant submanifold of $\g^*$. Any $x\in \sigma$ admits an invariant open neighbourhood $W$ in $\sigma$ with a $G$-equivariant homeomorphisms $\Phi$ and a $G$-equivariant diffeomorphism $\varphi$ for which 
\begin{center} 
\begin{tikzcd} 
J^{-1}(W)\arrow[rr,"\Phi", "\sim"'] \arrow[d,"J"'] & & J^{-1}(\L_x)\times W_+ \arrow[d,"J\times\mathrm{id}_{W_+}"] \\
W \arrow[rd,"q_+"']\arrow[rr,"\varphi", "\sim"'] & & \L_x\times W_+\arrow[ld,"\mathrm{pr}_2"] \\
  & W_+ & 
\end{tikzcd} 
\end{center}
commutes. These can be chosen such that $\Phi$ restricts to the identity map on $J^{-1}(\L_x)$ and extends to a smooth map into $S\times W_+$ on an open neighbourhood of $J^{-1}(W)$ in $S$ and such that its inverse extends to a smooth map into $S$ on an open neighbourhood of $J^{-1}(\L_x)\times W_+$ in $S\times W_+$.
\end{thm}
Here, the $G$-action on $J^{-1}(\L_x)\times W_+$ is given by $g\cdot (p,-)=(g\cdot p,-)$. Similarly for that on $\g^*\times W_+$.
\begin{rem}\label{rem:trivstratpreserving} Both $J^{-1}(W)$ and $J^{-1}(\L_x)$ are naturally stratified by connected components of the subspaces obtained by intersecting $J^{-1}(W)$, respectively $J^{-1}(\L_x)$, with $G$-orbit types in $S$. Taking products of the strata of $J^{-1}(\L_x)$ with $W_+$ we obtain a natural stratification of $J^{-1}(\L_x)\times W_+$ (assuming that $W_+$ is connected). It follows from $G$-equivariance of $\Phi$ that it respects these stratifications. The same goes for the induced map between $J^{-1}(W)/G$ and $J^{-1}(\L_x)/G\times W_+$. 
\end{rem}
\begin{rem} The proof of Theorem \ref{thm:localtriviality} will show that, if $x\in \sigma_+$, then we can in fact always choose $(\Phi,\varphi)$ so that $\varphi$ restricts to the map $y\mapsto (x,y)$ on $W_+$, or in other words, so that $\varphi^{-1}$ is the map given by $(g\cdot x,\alpha)\mapsto g\cdot \alpha$. If $G$ is a torus, this means that each $x\in \sigma$ admits an open neighbourhood $W$ in $\sigma$ and a $G$-equivariant homeomorphism $\Phi$ for which 
\begin{center} 
\begin{tikzcd} 
J^{-1}(W)\arrow[rr,"\Phi", "\sim"'] \arrow[rd,"J"'] & & J^{-1}(x)\times W \arrow[ld,"\mathrm{pr}_{2}"] \\
& W & 
\end{tikzcd} 
\end{center} commutes and which can be chosen so that it restricts to the identity on $J^{-1}(x)$ and so that it and its inverse have smooth extensions as in Theorem \ref{thm:localtriviality}. 
\end{rem}

\subsubsection{Proof of local triviality}
Although more technical, the essence of our proof of Theorem \ref{thm:localtriviality} is the same as one of the classical proofs of Ehresmann's fibration theorem for proper submersions: we lift coordinate vector fields on $\sigma_+$ along $J$ and use the flows of these lifts to obtain the desired trivializations. The following lemma ensures the existence of the types of lifts that we will use. 
\begin{lemma}\label{lemma:liftingvecfields:stratumwise}
In the setting of Theorem \ref{thm:localtriviality}, for any constant vector field $X$ on $\sigma_+$ there is a pair of $G$-invariant vector fields ($\widehat{X}$,$\widehat{X}^\sigma)$ defined on a $G$-invariant open $U$ in $S$ around $J^{-1}(\sigma)$ and on $\sigma$, respectively, that lift $X$ in the sense that 
\begin{equation}\label{eqn:liftcondition} \mathrm{d}J(\widehat{X}_p)=\widehat{X}^\sigma_{J(p)},\quad \d(q_+\vert_\sigma)(\widehat{X}^\sigma_x)=X_{q_+(x)},
\end{equation} for all $p\in J^{-1}(\sigma)$ and $x\in \sigma$. 
\end{lemma}
Here and throughout, we call a vector field $X$ on $\sigma_+$ constant if there is a vector $\beta$ in the linear subspace of $\g^*$ that the affine hull of $\sigma_+$ is modelled on, such that $X$ is given by
\begin{equation}\label{eqn:dfnconstantvectorfield} X_x=\left.\frac{\d}{\d t}\right|_{t=0} x+t\beta, \quad x\in \sigma_+.
\end{equation}
\begin{proof} It suffices to show that for each orbit in $\sigma$ there is a pair of $G$-invariant vector fields $(\widetilde{X},\widetilde{X}^\sigma)$, the second defined on some invariant open $V$ in $\g^*$ around that orbit and the first on an invariant open in $S$ around $J^{-1}(V\cap \sigma)$, such that \eqref{eqn:liftcondition} holds for all $p\in J^{-1}(V\cap\sigma)$ and $x\in V\cap \sigma$. Indeed, since $J$ is proper we can, after shrinking each such open $V$, arrange for the corresponding $\widehat{X}$ to be defined on $J^{-1}(V)$ and then, using a $G$-invariant partition of unity on the open covered by the opens $V$, we can patch the pairs of vector fields $(\widetilde{X},\widetilde{X}^\sigma)$ together to obtain a pair of vector fields as in the statement of the lemma. To prove the existence of such a pair of lifts near a given orbit $\L$ in $\sigma$, consider $x\in \L\cap\sigma_+$ and fix a slice $\mathfrak{S}$ through $x$ that is open in $x+\g^*_x$. In view of part $(c)$ of Theorem \ref{thm:symplectic:cross-section} it is enough to construct, after possibly shrinking $\mathfrak{S}$, a pair of $G_x$-invariant vector fields $(\overline{X},\overline{X}^\sigma)$, the second defined on $\mathfrak{S}$ and the first on an invariant open in $J^{-1}(\mathfrak{S})$ around $J^{-1}(\mathfrak{S}\cap\sigma)$, such that the respective conditions \eqref{eqn:liftcondition} hold at all points in $J^{-1}(\mathfrak{S}\cap\sigma)$ and in $\mathfrak{S}\cap \sigma$. To do so, note that $\mathfrak{S}\cap \sigma$ is an open subset of
\begin{equation*} x+\bigcap_{p\in J^{-1}(x)} (\g_p^0)^{G_x}.
\end{equation*}
So, after shrinking $\mathfrak{S}$, we can arrange for $\mathfrak{S}\cap \sigma$ to be an open subset of $\sigma_+$ and define $\overline{X}^\sigma:=X\vert_{\mathfrak{S}\cap \sigma}$. Then $\overline{X}^\sigma$ satisfies the second condition in \eqref{eqn:liftcondition} at all points in $\mathfrak{S}\cap \sigma$, because $q_+:\sigma\to \sigma_+$ restricts to the inclusion map $ \mathfrak{S}\cap \sigma\hookrightarrow\sigma_+$. To construct the desired lift on $J^{-1}(V)$ we will first construct, for each $G_x$-orbit in $J^{-1}(\mathfrak{S})$, a $G_x$-invariant vector field $\widehat{X}$ on an invariant open $U$ in $J^{-1}(\mathfrak{S})$ around that orbit, such that the first condition in \eqref{eqn:liftcondition} holds for all $p\in U\cap J^{-1}(\mathfrak{S}\cap\sigma)$. To construct such a local lift near a given $G_x$-orbit $\O$ in $J^{-1}(\mathfrak{S}\cap\sigma)$, we will use the MGS normal form theorem for the symplectic cross-section of $J$ at $\mathfrak{S}$. Consider $p\in \O$ such that $J(p)\in \sigma_+$. The corresponding local model near $\O$ for the symplectic cross-section of $J$ at $\mathfrak{S}$ is the map \begin{equation}\label{eqn:mommapMGSlocalmodel}
G_x\times_{G_p}(\g_p^0\oplus V)\to \g_x^*,\quad [g,\alpha,v]\mapsto g\cdot \big(\alpha+(\p\circ J_V)(v)\big),
\end{equation} where $\p$, $V$, $J_V$ and the $G_x$-action are as in the proof of Proposition \ref{prop:images:canhamstrata:affineopen}. Let $\beta$ be the vector defining the constant vector field $X$, as in \eqref{eqn:dfnconstantvectorfield}. Note that $G_x=G_{J(p)}$ since both $x$ and $J(p)$ lie in the same open face of $\t^*_+$ (for both belong to $\sigma_+$). Therefore $\beta$ belongs to $(\g_p^0)^{G_x}$ and, hence, we can define an explicit $G_x$-invariant lift of $X$ on $G_x\times_{G_p}(\g_p^0\oplus V)$ given by
\begin{equation*} [g,\alpha,v]\mapsto \left.\frac{\d}{\d t}\right|_{t=0}[g,\alpha+t\beta,v].
\end{equation*} Via the symplectomorphism obtained from the MGS normal form theorem, this transports to the desired lift on an invariant open $U$ around $\O$ in $J^{-1}(\mathfrak{S})$. Having constructed such a local lift around each $\O\in J^{-1}(\mathfrak{S}\cap \sigma)$, by patching these together using a $G$-invariant partition of unity we obtain the desired lift $\overline{X}$ on an invariant open in $J^{-1}(\mathfrak{S})$ around $J^{-1}(\mathfrak{S}\cap \sigma)$.  \end{proof}

Besides this, we will use the next lemmas to ensure the completeness of the flows of the vector fields on a domain as in Theorem \ref{thm:localtriviality}. Given a topological space $X$, by a \textbf{local flow} we will mean a continuous map $\Phi:U\to X$, $(t,p)\mapsto \Phi^t(p)$, with $U$ open in $\R\times X$, for which the following hold.
\begin{itemize}
\item[(i)] For each $p\in X$, the set $I_p:=\{t\in \R\mid (t,p)\in U\}$ is an open interval containing $0$.
\item[(ii)] For any $p\in X$, $t\in I_p$ and $s\in I_{\Phi^t(p)}$, it holds that $t+s\in I_p$ and $\Phi^s(\Phi^t(p))=\Phi^{t+s}(p)$.
\end{itemize}
\begin{lemma}\label{lemma:relatedlocalflows} Let $f:X\to Y$ be a proper continuous map between topological spaces and suppose $\Phi_X:U\to X$ and $\Phi_Y:V\to Y$ are $f$-related local flows, meaning that for all $(t,p)\in U$ it holds that $\big(t,f(p)\big)\in V$ and $f\big(\Phi_X^t(p)\big)=\Phi_Y^t\big(f(p)\big)$. Then $I_{f(p)}\subset I_p$ for all $p\in X$.
\end{lemma}
\begin{proof} Let $p\in X$, let $t\in I_{f(p)}$ and write $I_p=]a,b[$ with $a<0<b$. If $t=0$ then $t\in I_p$ by (i) above. Next, suppose $t>0$. Since $K:=\Phi_Y([0,t]\times \{f(p)\})$ is compact, its pre-image under $f$ is compact as well, by properness of $f$. So, by the tube lemma there is a $\delta\in ]0,b[$ such that 
\begin{equation}\label{eqn:pflemma:relatedlocalflows} ]-\delta,\delta[\times f^{-1}(K)\subset U.
\end{equation} Note that, if $t\notin I_p$, then $t\geq b$, so $s:=b-\tfrac{\delta}{2}\in]0,b[\subset ]0,t[$ and, hence, $f\big(\Phi_X^s(p)\big)=\Phi_Y^s\big(f(p)\big)\in K$, so that $]-\delta,\delta[\subset I_{\Phi_X^s(p)}\subset I_p-s$ by \eqref{eqn:pflemma:relatedlocalflows} and (ii) above, which is a contradiction since $\tfrac{\delta}{2}+s=b$. Therefore, $t\in I_p$, as was to be shown. By analogous reasoning one sees that $t\in I_p$ when $t<0$.
\end{proof}
\begin{lemma}\label{lemma:flowcontainedcriterion} Let $M$ be a smooth manifold, $X$ a vector field on $M$ and $C$ a closed subset of $M$. If there is a partition $\P$ of $C$ into submanifolds of $M$ such that $X_p\in T_pP$ for each $P\in \P$ and each $p\in P$, then the flow of $X$ preserves $C$, meaning that $\Phi_X^t(p)\in C$ for all $p\in C$ and $t\in I_p$.  
\end{lemma}
\begin{proof} Let $p\in C$. By continuity of the maximal integral curve $\Phi_X^{(\cdot)}:I_p\to M$ of $X$, the subset
\begin{equation*} I_C:=\{t\in I_p\mid \Phi_X^t(p)\in C\}
\end{equation*} is closed in $I_p$. We will show that it is open in $I_p$ as well, so that it must be all of $I_p$ since it contains $0$ and $I_p$ is connected. Let $t\in I_C$. Then $q:=\Phi_X^t(p)\in C$, so there is a unique $P\in \P$ containing $q$. Since $X$ is tangent to the submanifold $Q$, it restricts to a vector field $X_Q$ on $Q$. Moreover, the maximal integral curve $I_{q,Q}\to Q$ of $X_Q$ starting at $q$ is an integral curve for $X$ as well, so $I_{q,Q}+t\subset I_p$ and this maximal integral curve of $X_Q$ takes the value $\Phi_X^{t+s}(p)$ for $s\in I_{q,Q}$. Since integral curves of $X_Q$ take values in $Q$ and, hence, in $C$, this implies that $I_C$ contains the open interval $I_{q,Q}+t$ around $t$, which shows that $I_C$ is indeed open in $I_p$. 
\end{proof}
\begin{lemma}\label{lemma:integralcurvepreservation} In the setting of Theorem \ref{thm:localtriviality}, suppose $X$ is a vector field on $\sigma_+$ and $(\widehat{X},\widehat{X}^\sigma)$ a pair lifting $X$ as in Lemma \ref{lemma:liftingvecfields:stratumwise}, with domain $U$ of $\widehat{X}$ chosen so that $J^{-1}(\sigma)$ is closed in $U$. Then each integral curve of $\widehat{X}$ that starts in $J^{-1}(\sigma)$ remains in $J^{-1}(\sigma)$ and is mapped to an integral curve of $\widehat{X}^\sigma$ by $J$.
\end{lemma}
\begin{proof} Once we have shown the first claim (integral curves of $\widehat{X}$ that start in $J^{-1}(\sigma)$ remain in there), the second claim (such integral curves are mapped to integral curve of $X$ by $J$) readily follows from the first condition in \eqref{eqn:liftcondition}. To prove the first claim, note that $J:\Sigma\to \g^*$ intersects $\sigma$ cleanly for each stratum $\underline{\Sigma}\in \sS_\mathrm{Ham}(\underline{S})$, because $\Sigma\cap J^{-1}(\sigma)$ is the pre-image of the submanifold $\sigma$ of $\mathrm{Af}(\underline{J}_+(\underline{\Sigma}))$ under the composition of $\underline{J}_+:\underline{\Sigma}\to \mathrm{Af}(\underline{J}_+(\underline{\Sigma}))$ with the quotient map $q_S:\Sigma\to \underline{\Sigma}$, which are submersions by Proposition \ref{prop:images:canhamstrata:affineopen} and \cite[Proposition 2.56]{Mol1}. This means that $\Sigma\cap J^{-1}(\sigma)$ is a submanifold of $S$ with tangent space given by
\begin{equation}\label{eqn:tangentsp:integralcurvepreservation} T_p(\Sigma\cap J^{-1}(\sigma))=T_p\Sigma\cap \d J_p^{-1}(T_x\sigma),\quad p\in \Sigma\cap J^{-1}(\sigma),\quad x:=J(p). 
\end{equation}  Since the submanifolds $\Sigma\cap J^{-1}(\sigma)$ partition $J^{-1}(\sigma)$, by Lemma \ref{lemma:flowcontainedcriterion} it is enough to show that $\widehat{X}$ is tangent to $\Sigma\cap J^{-1}(\sigma)$ for each such $\underline{\Sigma}$. 
In view of the proof of \cite[Proposition 2.56]{Mol1}, $T_p\Sigma$ consists of those $v\in T_pS$ for which $v$ is tangent to the $G$-orbit type stratum through $p$ and for which $\d J_p(v)$ is tangent to the $G$-orbit type stratum through $J(p)$. So, since $\sigma$ is a submanifold of the $G$-orbit type stratum through $J(p)$, \eqref{eqn:tangentsp:integralcurvepreservation} consists of those $v\in T_pS$ that are tangent to the $G$-orbit type stratum through $p$ and for which $\d J_p(v)$ is tangent to $\sigma$. Since $\widehat{X}$ is $G$-invariant and lifts $\widehat{X}^\sigma$, $\widehat{X}_p$ is indeed such a vector for each $p\in \Sigma\cap J^{-1}(\sigma)$. Therefore, the lemma follows.  
\end{proof}
\begin{proof}[Proof of Theorem \ref{thm:localtriviality}] Consider a basis $\beta_1,...,\beta_n$ of the linear subspace of $\g^*$ that the affine hull of $\sigma_+$ is modelled on. For each $i$, let $X_i$ be the constant vector field on $\sigma$ given by $\beta_i$ as in \eqref{eqn:dfnconstantvectorfield}. Fix lifts $(\widehat{X}_i,\widehat{X}_i^\sigma)$ as in Lemma \ref{lemma:liftingvecfields:stratumwise}, all defined on one $G$-invariant open $U$ in $S$ chosen so that $J^{-1}(\sigma)$ is closed in $U$ (which is possible since $\sigma$ is locally closed in $\g^*$). Choose $\varepsilon >0$ such that
\begin{equation*} C_+:=x+[-\varepsilon,\varepsilon]\cdot\beta_1+...+[-\varepsilon,\varepsilon]\cdot\beta_n\subset \sigma_+.
\end{equation*} Consider
\begin{equation*} W_+:=x+]-\varepsilon,\varepsilon[\cdot\beta_1+...+]-\varepsilon,\varepsilon[\cdot\beta_n, \quad C^i_+:=x+[-\varepsilon,\varepsilon]\cdot\beta_1+...+[-\varepsilon,\varepsilon]\cdot\beta_{i-1},
\end{equation*} and the diffeomorphism
\begin{equation*}  \tau:W_+\to ]-\varepsilon,\varepsilon[^n,\quad x+t_1\beta_1+...+t_n\beta_n\mapsto (t_1,...,t_n). 
\end{equation*}
Since the time $t$ flow of $X_i$ is translation by $t\beta_i$, each such flow is defined on $C^i_+$ for all $t\in [-\varepsilon,\varepsilon]$. 
By Lemma \ref{lemma:integralcurvepreservation} the local flow of $\widehat{X}_i$ restricts to a local flow on $J^{-1}(\sigma)$ which is related to that of $\widehat{X}^\sigma_i$ by the map $J:J^{-1}(\sigma)\to \sigma$. Moreover, since $\widehat{X}_i^\sigma$ lifts $X^\sigma_i$, their local flows are related by $q_+$. In view of this and Lemma \ref{lemma:relatedlocalflows}, we can consider the maps
\begin{align*} 
&\Psi:J^{-1}(\L_x)\times W_+\to J^{-1}(W),\quad (p,x+t_1\beta_1+...+t_n\beta_n)\mapsto (\Phi_{\widehat{X}_n}^{t_n}\circ ...\circ \Phi_{\widehat{X}_1}^{t_1})(p),\\
 &\Phi:J^{-1}(W)\to J^{-1}(\L_x)\times W_+,\quad p\mapsto \big((\Phi_{\widehat{X}_1}^{-(\tau_1\circ J_+)(p)}\circ ...\circ \Phi_{\widehat{X}_n}^{-(\tau_n\circ J_+)(p)})(p),J_+(p)\big),
\end{align*} and
\begin{align*} 
&\psi:\L_x\times W_+\to W,\quad (y,x+t_1\beta_1+...+t_n\beta_n)\mapsto (\Phi_{\widehat{X}^\sigma_n}^{t_n}\circ ...\circ \Phi_{\widehat{X}^\sigma_1}^{t_1})(y),\\
& \varphi:W\to \L_x\times W_+,\quad y\mapsto \big((\Phi_{\widehat{X}^\sigma_1}^{-(\tau_1\circ q_+)(y)}\circ ...\circ \Phi_{\widehat{X}^\sigma_n}^{-(\tau_n\circ q_+)(y)})(y),q_+(y)\big).
\end{align*} These first two and the last two of these are mutually inverse, so the first two are homeomorphisms and the last two are diffeomorphisms. Moreover, $\Phi$ and $\varphi$ fit into the commutative diagram in the statement of the theorem. They are $G$-equivariant, because the vector fields $\widehat{X}_i$ and $\widehat{X}_i^\sigma$ are $G$-invariant, and by construction they restrict to the respective identity maps on $J^{-1}(\L_x)$ and $
\L_x$. Since the set of $(p,t)\in U\times \R^n$ for which 
\begin{equation*} (\Phi_{\widehat{X}_n}^{t_n}\circ ...\circ \Phi_{\widehat{X}_1}^{t_1})(p)
\end{equation*} is defined is an open around $J^{-1}_+(x)\times [-\varepsilon,\varepsilon]^n$, by the tube lemma there is an invariant open $V$ in $S$ around $J_+^{-1}(x)$ such that the same formula as for $\Psi$ can be used to define a smooth map from $V\times W$ into $S$. So, $\Psi$ admits the desired smooth extension. To show the same for $\Phi$, consider a smooth map $\widehat{J}_+:V\to W_+$, with $V$ an open in $S$ around $J^{-1}(W)$, that coincides with $J_+$ on $J^{-1}(W)$. Such a map can be constructed, for instance, by choosing a smooth retract $r$ from an open $O$ in $\g^*$ onto the submanifold $W$, and defining $\widehat{J}_+$ as the composition
\begin{equation*} V:=J^{-1}(O)\xrightarrow{J} O\xrightarrow{r}W\to W_+,
\end{equation*} where the last map is the restriction of $\g^*\to \underline{\g^*}\xrightarrow{\eqref{eqn:sweeping:map}^{-1}} \t^*_+$. Since the set of $p\in J^{-1}(O)$ for which 
\begin{equation*} (\Phi_{\widehat{X}_1}^{-(\tau_1\circ \widehat{J}_+)(p)}\circ ...\circ \Phi_{\widehat{X}_n}^{-(\tau_n\circ \widehat{J}_+)(p)})(p)
\end{equation*} is defined is an open around $J^{-1}(W)$, we can define a smooth extension of $\Phi$ on this set using the same formula as for $\Phi$, but with $J_+$ replaced by $\widehat{J}_+$.
\end{proof}
\subsection{Push-forwards of the constant sheaf}\label{sec:push-forward-sheaves} 
\subsubsection{Constructibility} We now turn to Corollary \ref{cor:constructibility:pushforwards}. Throughout this section, we consider cohomology with $\R$-coefficients and omit this from the notation. We will use the following standard lemma.
\begin{lemma}\label{lemma:constructibilitypushforward:criterion} Suppose that $f:X\to Y$ is a proper, continuous, surjective map between locally compact Hausdorff topological spaces and that $\sS$ a stratification of $Y$ such that $f:f^{-1}(\sigma)\to \sigma$ is a topologically locally trivial fibration for each $\sigma\in \sS$. Then $R^i(f_*)(\underline{\R})$ restricts to a locally constant sheaf on each stratum of $\sS$, for each $i\in \mathbb{N}$. If in addition $H^i(f^{-1}(y))$ is finite-dimensional for each $y\in Y$, then $R^i(f_*)(\underline{\R})$ is constructible. 
\end{lemma}
\begin{proof} This is a direct consequence of the proper base change theorem. 
\end{proof}
We will further use the lemma below. 
\begin{lemma}\label{lemma:finitegencohom:totalspacefibbun} Let $\pi:P\to M$ be a topological fiber bundle over a compact manifold, with fiber $F$ such that $H^i(F)$ is finite-dimensional for all $i$. Then $H^i(P)$ is also finite-dimensional for all $i$.
\end{lemma}
\begin{proof} Since $M$ is a compact manifold, it admits a finite good cover $U_1,...,U_n$ with the property that over each open in the cover the fiber bundle $\pi:P\to M$ is trivializable. By inspecting the Mayer-Vietoris sequence of the pair of opens $\big(\pi^{-1}(U_1),\pi^{-1}(U_2)\cup ... \cup \pi^{-1}(U_n)\big)$ and argueing by induction on the size of such good covers, the lemma follows.
\end{proof}
\begin{proof}[Proof of Corollary \ref{cor:constructibility:pushforwards}] Since the fibers of $\underline{J}_+$ (the reduced spaces) are compact Whitney (b)-regular differentiable stratified spaces (in the terminology used in \cite{Mol1}), they admit finite good covers (see, e.g., \cite[Section 7]{PfPoTa}) and, hence, their cohomology is finite-dimensional in each degree. 
Because of this and Corollary \ref{cor:thm:stratification:image:transverse:momentum:map}, Lemma \ref{lemma:constructibilitypushforward:criterion} applies to $\underline{J}_+:\underline{S}\to \Delta_+$ and $\sS_\mathrm{Ham}(\Delta_+)$, so part (a) follows. Part (b) concerns the push-forward along a map that is not proper. However, since $G$ is compact, we can find a compact, connected $G$-invariant subspace $EG_i$ such that $H^n(EG_i)=0$ for all $0<n\leq i$ and derive part (b) by instead applying Lemma  \ref{lemma:constructibilitypushforward:criterion} to the map 
\begin{equation}\label{eqn:pushforward:approximation} 
\underline{J}_+^{EG_i}:=\underline{J}_+\circ \mathrm{pr}_{\underline{S}}:EG_i\times_GS\to \Delta_+.
\end{equation} To see this note that, given an open $W$ in $\Delta_+$, the fibers of both of the canonical fiber bundles
\begin{center} 
\begin{tikzcd} & (EG_i\times EG)\times_GJ_+^{-1}(W)\arrow[ld]\arrow[rd] & \\
EG_i\times_GJ_+^{-1}(W)& & EG\times_GJ_+^{-1}(W)
\end{tikzcd} 
\end{center} are connected and have trivial cohomology in all degrees $0<n\leq i$, so the bundle maps induce isomorphisms in cohomology in all degrees $n\leq i$ (see, e.g., \cite[Corollary A.7]{Qui}) and, in particular, they induce an isomorphism 
\begin{equation*} H^i(EG_i\times_GJ_+^{-1}(W))\cong H^i(EG\times_GJ_+^{-1}(W)).
\end{equation*} Because these isomorphisms are compatible with restriction to smaller opens in $\Delta_+$, they induce an isomorphism of sheaves
\begin{equation*} R^i(\underline{J}^{EG_i}_+)_*(\underline{\R})\cong R^i(\underline{J}^{EG}_+)_*(\underline{\R}).
\end{equation*}
In contrast to $\underline{J}_+^{EG}$, \eqref{eqn:pushforward:approximation} is proper. Moreover, its fibers have finite-dimensional cohomology in each degree by Lemma \ref{lemma:finitegencohom:totalspacefibbun}, since each of its fibers is a topological fiber bundle over $EG_i/G$ with fiber $J_+^{-1}(x)$, for some $x\in \Delta_+$, which has finite-dimensional cohomology in each degree because, like $\underline{J}^{-1}_+(x)$, it is a compact Whitney (b)-regular differentiable stratified space (stratified by connected components of $G$-orbit types -- cf. \cite{LeSj}). Finally, it follows from Theorem \ref{thm:localtriviality} that \eqref{eqn:pushforward:approximation} restricts to a topologically locally trivial fibration over each stratum of $\sS_\mathrm{Ham}(\Delta_+)$. So, Lemma \ref{lemma:constructibilitypushforward:criterion} applies to \eqref{eqn:pushforward:approximation} and part (b) follows. 
\end{proof} 

\subsubsection{A local invariant cycle type result} Next, we address Proposition \ref{prop:specializationmaps}. Before turning to its proof, let us elaborate on the definition of the maps \eqref{eqn:specialization:maps}. Suppose that $\gamma$ is a path in $\Delta_+$ that starts at $x$, ends at $y$ and remains in a single stratum $\sigma_+\in \sS_\mathrm{Ham}(\Delta_+)$ except for at its end-point. The associated map 
\begin{equation*} s^\gamma:H^i(\underline{S}_{\L_y})\to H^i(\underline{S}_{\L_x}) 
\end{equation*} is defined as follows. Given $\mathrm{c}\in H^i(\underline{S}_{\L_y})$, since $H^i(\underline{S}_{\L_y})$ is the stalk of $R^i(\underline{J}_+)_*(\underline{\R})$ at $y$, there is an open $U_+$ around $y$ and a $\widehat{\mathrm{c}}\in H^i(\underline{J}_+^{-1}(U_+))$ such that $\widehat{c}\,\vert_{\underline{S}_{\L_y}}=\mathrm{c}$. The map $s^\gamma$ sends $\mathrm{c}$ to 
\begin{equation*} s^\gamma(\mathrm{c}):=m^{-1}_{\gamma_\varepsilon}(\,\widehat{\mathrm{c}}\,\vert_{\underline{S}_{\gamma(\varepsilon)}}),
\end{equation*} where $\varepsilon\in [0,1[$ is chosen such that $\gamma(\varepsilon)\in U_+$, $\gamma_\varepsilon$ is the path in $\sigma_+$ given by $\gamma_\varepsilon(t)=\gamma(\varepsilon t)$ for $t\in [0,1]$, and the map 
\begin{equation*} m_{\gamma_\varepsilon}:H^i(\underline{S}_{\L_x})\xrightarrow{\sim} H^i(\underline{S}_{\L_{\gamma(\varepsilon)}})
\end{equation*} is the associated monodromy isomorphism of the locally constant sheaf $R^i(\underline{J}_+)_*(\underline{\R})\vert_{\sigma_+}$. This does not depend on the choice of $\widehat{\mathrm{c}}$ or $\varepsilon$, as follows from the fact that the monodromy isomorphism associated to a path only depends on its homotopy class and the fact that the monodromy isomorphism associated to a composition of paths is the composition of those associated to the individual paths. The associated maps $H^i_G(S_{\L_y})\to H^i_G(S_{\L_x})$ are defined similarly. \\

For the proof of Proposition \ref{prop:specializationmaps} we will use the following extension of Ehresmann's fibration theorem (for which we could not find a reference, although it appears implicitly in \cite{DuHe}).
\begin{lemma}\label{lemma:Ehresmannfib:equiv} Let $G$ be a Lie group, $S$ be a manifold equipped with a smooth proper $G$-action and let $f:S\to M$ be a $G$-invariant, proper and surjective submersion. Then the map $\underline{f}:\underline{S}\to M$ induced by $f$ is a topologically locally trivial fibration.
\end{lemma}
\begin{proof} Let $x\in M$. Since the $G$-action is proper, there is a $G$-invariant Ehresmann connection for $f$, which can be used to lift a set of coordinate vector fields $X_1,...,X_n$ defined around $x$ to $G$-invariant vector fields $\widehat{X}_1,...,\widehat{X}_n$. Like in the proof of Theorem \ref{thm:localtriviality}, the flows of such lifts compose to give a smooth trivialization
\begin{center} \begin{tikzcd} f^{-1}(U)\arrow[rr,"\Phi", "\sim"'] \arrow[rd,"f"'] & & f^{-1}(x)\times U \arrow[ld,"\mathrm{pr}_2"] \\
  & U &
\end{tikzcd} 
\end{center} over an open $U$ around $x$, which is $G$-equivariant with respect to the $G$-action on $J^{-1}(x)\times U$ given by $g\cdot (p,y)=(g\cdot p,y)$. This descends to a trivialization of $\underline{f}$ over $U$.
\end{proof}
\begin{proof}[Proof of Proposition \ref{prop:specializationmaps}] Throughout this proof, for each $x\in \Delta_+$ we will canonically identify $\underline{S}_{\L_x}$ with $\underline{S}_x:=J^{-1}(x)/G_x$. Turning to the proof, let $x,y$ and $\gamma$ be as in the statement of the proposition. To define the $\mathbb{W}_y$-action, note that $y$ admits an open neighbourhood $U$ in $\Delta$ such that $J:J^{-1}(U)\to U$ is a submersion, since $y$ is a regular value of $J$ and $J$ is proper as map into its image. Consider $\varepsilon\in [0,1[$ such that $\gamma([\varepsilon,1])\subset U$. Let $\gamma_\varepsilon$ be the path from $x$ to $\gamma(\varepsilon)$ given by $\gamma_\varepsilon(t)=\gamma(\varepsilon t)$ and let $\gamma^\varepsilon$ be the path in $U$ from $\gamma(\varepsilon)$ to $y$ given by $\gamma^\varepsilon(t)=\gamma((1-t)\varepsilon+t)$. By Lemma \ref{lemma:Ehresmannfib:equiv}, $\underline{J}_T:J^{-1}(U_{\t^*})/T\to U_{\t^*}$ is a topologically locally trivial fibration, where $U_{\t^*}:=U\cap \t^*$. So, for each $i$, the sheaf $R^i({\underline{J}_T}_*)(\underline{\R})$ is locally constant. Therefore, we can transfer the $\mathbb{W}_y$-action on $H^i(\underline{J}_T^{-1}(y))$ induced by the natural $\mathbb{W}_y$-action on $\underline{J}_T^{-1}(y)=J^{-1}(y)/T$ to an action on $H^i(\underline{S}_x)$ via the composition
\begin{equation}\label{eqn:comp:monodromy:defweylgpaction} H^i\big(\underline{S}_x\big) \xrightarrow{m_{\gamma_\varepsilon}} H^i\big(\underline{S}_{\gamma(\varepsilon)}\big)\xrightarrow{m_{\gamma^\varepsilon}}H^i\big(\underline{J}_T^{-1}(y)\big)
\end{equation} of the monodromy isomorphism of $R^i(\underline{J}_*)(\underline{\R})$ along $\gamma_\varepsilon$ with that of $R^i({\underline{J}_T}_*)(\underline{\R})$ along $\gamma^\varepsilon$. Note here that $\underline{S}_{\gamma(\varepsilon)}=\underline{J}_T^{-1}(\gamma(\varepsilon))$, because $\gamma(\varepsilon)$ belongs to the interior of $\t^*_+$ (since $x$ does so and $\gamma([0,1[)$ is contained in the stratum of $\sS_\mathrm{Ham}(\underline{\Delta})$ through $x$). The isomorphism \eqref{eqn:comp:monodromy:defweylgpaction} and, hence, the resulting $\mathbb{W}_y$-action on $H^i(\underline{S}_x)$ do not depend on the choice of $U$ or $\varepsilon$. Having defined this action, it remains to show that the image of $s^\gamma:H^i(\underline{S}_y)\to H^i(\underline{S}_x)$ is the $\mathbb{W}_y$-fixed point set. For this, note that $s^\gamma$ is equal to the composition
\begin{equation*} H^i\big(\underline{S}_y\big)\xrightarrow{\pi^*} H^i\big(\underline{J}_T^{-1}(y)\big)\xrightarrow{\sim} H^i\big(\underline{S}_x\big),
\end{equation*} where $\pi$ is the canonical map from $\underline{J}_T^{-1}(y)=J^{-1}(y)/T$ to $\underline{S}_y=J^{-1}(y)/G_y$ and the second map is given by the inverse of \eqref{eqn:comp:monodromy:defweylgpaction}. Therefore, it is enough to show that $\pi^*$ is an isomorphism onto the $\mathbb{W}_y$-fixed point set in $H^i\big(\underline{J}_T^{-1}(y)\big)$. Since $\mathbb{W}_y$ is finite, it follows for instance from \cite[Corollary to Proposition 5.2.3, Theorem 5.3.1]{Gro} that the quotient map from $\underline{J}_T^{-1}(y)$ to $J^{-1}(y)/N_{G_y}(T)$ for the $\mathbb{W}_y$-action on $\underline{J}_T^{-1}(y)$ induces an isomorphism
\begin{equation}\label{eqn:pflocinvprop} 
H^i\big(J^{-1}(y)/N_{G_y}(T)\big)\cong H^i\big(\underline{J}_T^{-1}(y)\big)^{\mathbb{W}_y},
\end{equation} where $N_{G_y}(T)$ denotes the normalizer of $T$ in $G_y$. Since $G_y/N_{G_y}(T)$ is connected and has trivial cohomology in all strictly positive degrees (see, e.g., \cite[Chapter III, Lemma 1.1]{Hsi}) it follows from the Vietoris-Begle theorem for sheaf cohomology (see, e.g., \cite[Corollary A.7]{Qui}) that the upper horizontal map in the canonical commutative diagram
\begin{center} 
\begin{tikzcd} EG_y\times_{N_{G_y}(T)}J^{-1}(y)\arrow[r]\arrow[d] & EG_y\times_{G_y}J^{-1}(y)\arrow[d]\\
J^{-1}(y)/N_{G_y}(T) \arrow[r] & \underline{S}_y
\end{tikzcd}
\end{center}
induces an isomorphism in cohomology,
\begin{equation*} H^i_{N_{G_y}(T)}(J^{-1}(y))\cong H^i_{G_y}(J^{-1}(y)),
\end{equation*} for it is a fiber bundle with fiber $G_y/N_{G_y}(T)$. Since the $G_y$-action on $J^{-1}(y)$ has only finite isotropy groups ($y$ being a regular value) and we use  $\mathbb{R}$-coefficients, the vertical maps also induce isomorphisms in cohomology. So, the same goes for the lower horizontal map. The proposition now follows by noting that $\pi^*$ is the composition of this isomorphism with \eqref{eqn:pflocinvprop}. 
\end{proof} 

\section{Linear variation}
\subsection{On Sjamaar's de Rham model for symplectic reduced spaces}
\subsubsection{The de Rham model}\label{subsec:deRhammodel:1} In \cite{Sj2} Sjamaar proved a de Rham theorem for the cohomology of singular symplectic reduced spaces, as well as some related results, which will be used when studying the variation of the reduced spaces. He mostly focuses on reduced spaces at level zero and indicates how the theory in \cite[Section 1-7]{Sj2} extends to reduced spaces at arbitrary coadjoint orbits in an addendum. Since we will use the latter, in this section we explicitly give the relevant definitions and statements from \cite{Sj2} at arbitrary coadjoint orbits. Let $G$ be a compact Lie group and $J:(S,\omega)\to \g^*$ a Hamiltonian $G$-space. Fix a coadjoint $G$-orbit $\L\subset \g^*$. We use the notation 
\begin{equation*} S_\L:=J^{-1}(\L),\quad \underline{S}_\L:=J^{-1}(\L)/G,
\end{equation*} and we let 
\begin{equation*} S_\L^\mathrm{prin},\quad  \underline{S}_\L^\mathrm{prin}
\end{equation*} denote the principal parts (the unions of all open strata) of the respective Lerman-Sjamaar stratifications \cite{LeSj} of $S_\L$ and $\underline{S}_\L^\mathrm{prin}$. Further, we let 
\begin{equation*} i_\L:S^\mathrm{prin}_\L\hookrightarrow S, \quad q_\L:S_\L^\mathrm{prin}\to \underline{S}_\L^\mathrm{prin},
\end{equation*} denote the inclusion and quotient map. 
\begin{defi} A differential form $\alpha$ on an open $\underline{U}$ in $\underline{S}_\L$ is a differential form on $\underline{U}\cap \underline{S}_\L^\mathrm{prin}$ with the property that $q_\L^*\alpha$ extends smoothly to an open in $S$ around the invariant subset $U$ in $S_\L$ corresponding to $\underline{U}$. 
\end{defi} 
\begin{rem} In the above definition, it is equivalent to require that $q_\L^*\alpha$ extends smoothly to a \textit{$G$-invariant} differential form on an invariant open in $S$ around $U$.
\end{rem}
The differential forms on $\underline{S}_\L$ form a complex of sheaves
\begin{equation*} \Omega^0_{\underline{S}_\L}\xrightarrow{\d_\mathrm{dR}} \Omega^1_{\underline{S}_\L} \xrightarrow{\d_\mathrm{dR}} ... 
\end{equation*} with differential given by the de Rham differential of forms on $\underline{S}_\L^\mathrm{prin}$, which we call the \textbf{de Rham complex} on $\underline{S}_\L$. Let $H^\bullet_\mathrm{dR}(\underline{S}_\L)$ denote the cohomology of the associated complex of global sections. 
\begin{ex}\label{ex:cohomclass:sympform:redsp} The symplectic form $\omega$ induces a differential form $\omega_{\underline{S}_\L}$ on $\underline{S}_\L$, determined by the property that 
\begin{equation*} i_\L^*\omega-(J\vert_{ S_\L^\mathrm{prin}})^*\omega_\L=q_\L^*\omega_{\underline{S}_\L},
\end{equation*} with $\omega_\L$ the KKS symplectic form on $\L$. So, $\omega$ defines a cohomology class $\varpi_\L:=[\omega_{\underline{S}_\L}]\in H^2_\mathrm{dR}(\underline{S}_\L)$.
\end{ex}

We now turn to the aforementioned de Rham theorem.
\begin{thm}[\cite{Sj2}]\label{thm:sjamaar:deRham} The de Rham complex on $\underline{S}_\L$ yields an acyclic resolution
\begin{equation*} 0\to \underline{\R}\hookrightarrow \Omega^0_{\underline{S}_\L}\xrightarrow{\mathrm{d}_\mathrm{dR}} \Omega^1_{\underline{S}_\L} \xrightarrow{\mathrm{d}_\mathrm{dR}} ...
\end{equation*} of the sheaf of locally constant functions with values in $\R$. Hence, there are natural isomorphisms
\begin{equation}\label{isomorphism:sjamaar:deRham}
H^\bullet_\mathrm{sing}(\underline{S}_\L,\R)\cong H^\bullet(\underline{S}_\L,\underline{\R})\cong H^\bullet_\mathrm{dR}(\underline{S}_\L).
\end{equation}
\end{thm}
Besides this, Sjamaar also showed the following.
\begin{thm}[\cite{Sj2}]\label{thm:sjamaar:finvol} Suppose that $\underline{S}_\L$ is compact. Then the following hold. 
\begin{itemize} 
\item[(i)] Its symplectic volume $\mathrm{Vol}(\underline{S}_\L,\omega_{\underline{S}_\L}):=\mathrm{Vol}(\underline{S}_\L^\mathrm{prin},\omega_{\underline{S}_\L})$ is finite. 
\item[(ii)] The integral over $\underline{S}_\L^\mathrm{prin}$ of any differential top-form on $\underline{S}_\L$ is finite.
\item[(iii)] Stokes' theorem: the integral over $\underline{S}_\L^\mathrm{prin}$ of an exact differential top-form on $\underline{S}_\L$ is zero. 
\end{itemize}
\end{thm}
Note here that, in \cite{Sj2}, it seems to be assumed for simplicity that all connected components of $\underline{S}^\mathrm{prin}_\L$ have the same dimension. By the same arguments as in that paper, the above remains true without this assumption, with the added note that, in Theorem \ref{thm:sjamaar:finvol}, by a \textit{top}-form we mean a finite collection of differential forms $\alpha_1,...,\alpha_k$ on $\underline{S}_\L$, one for each connected component $C^i_\L$ of $\underline{S}_\L^\mathrm{prin}$, such that $\alpha_i$ is a form of degree $n_i:=\dim(C^i_\L)$ and $\alpha_i\vert_{C^j_\L}=0$ for each $j\neq i$. Moreover, we call such a top-form exact if each $\alpha_i\in \Omega^{n_i}(\underline{S}_\L)$ is exact in the de Rham complex of $\underline{S}_\L$. Accordingly, $\mathrm{Vol}(\underline{S}_\L,\omega_{\underline{S}_\L})$ is the sum of the usual symplectic volumes of the connected components of $\underline{S}_\L$. 
\begin{rem}\label{rem:redspaceatpoints:comparison} The reduced space $\underline{S}_\L$ at a coadjoint orbit $\L$ can also be realized via reduction at a point $x\in \L$ and there is an analogous version of the above theory that uses the latter. To be more precise, given $x\in \L$, consider 
\begin{equation*} S_x:=J^{-1}(x),\quad \underline{S}_x:=J^{-1}(x)/G_x.
\end{equation*} By definition, the space $\underline{S}_x$ is the reduced space at $x$. The Lerman-Sjamaar stratification and differential forms on $\underline{S}_x$ can be defined as on $\underline{S}_\L$, with the role of the $G$-space $S_\L$ played by the $G_x$-space $S_x$. The symplectic form $\omega$ induces a differential form $\omega_{\underline{S}_x}$ on $\underline{S}_x$, determined by the simpler property that 
\begin{equation*} i_x^*\omega=q_x^*\omega_{\underline{S}_x},
\end{equation*} with $i_x:S^\mathrm{prin}_x\hookrightarrow S$ and $q_x:S_x^\mathrm{prin}\to \underline{S}_x^\mathrm{prin}$ the inclusion and quotient maps of the principal parts. Moreover, the inclusion $S_x\hookrightarrow S_\L$ then induces an isomorphism of stratified spaces
\begin{equation*} \underline{S}_x\cong \underline{S}_\L,
\end{equation*} the pull-back along which gives an isomorphism of de Rham complexes which takes the cohomology class of the reduced symplectic form on $\underline{S}_\L$ to that on $\underline{S}_x$. 
\end{rem}
\begin{rem}\label{rem:deRhamforms:non-principalstrata} For each stratum $\underline{\Sigma}$ of $\underline{S}_\L$, there is a canonical restriction map 
\begin{equation*} \Omega(\underline{S}_\L)\to \Omega(\underline{\Sigma}),\quad \alpha\mapsto \alpha\vert_{\underline{\Sigma}}.
\end{equation*} with $\alpha\vert_{\underline{\Sigma}}$ defined as follows. Given $\alpha\in \Omega(\underline{S}_\L)$, choose an extension $\widehat{\alpha}$ of $q_\L^*\alpha$ to an open in $S$ around $S_\L$. The restriction of $\widehat{\alpha}$ to the invariant submanifold $\Sigma$ of $S$ corresponding to $\underline{\Sigma}$ is $G$-basic and independent of the choice of extension, and $\alpha\vert_{\underline{\Sigma}}$ is defined to be the induced form. The form $\alpha\vert_{\underline{\Sigma}}$ can also be obtained using an extension of $q_x^*\alpha$ instead. Note here that the extension $\widehat{\alpha}$ is taken to be $G$-invariant in \cite{Sj2}, but that this is not necessary.
\end{rem}
For the purpose of studying the variation of the reduced spaces it is more natural to work with the reduced spaces $\underline{S}_\L$, because this does not require the choice of a point $x\in \L$ for each coadjoint orbit $\L$. Nonetheless, we will sometimes use their description via reduction at points instead, as this simplifies some arguments because the description of the reduced symplectic form on these does not involve the symplectic form on the corresponding coadjoint orbit. 
\subsubsection{Smooth singular chains and an explicit de Rham isomorphism} In preparation for the proof of the linear variation theorem, we will now show that the isomorphism \eqref{isomorphism:sjamaar:deRham} can be realized via a non-degenerate pairing between singular homology and de Rham cohomology which, like in the classical de Rham theorem, is defined explicitly on the level of complexes by integrating differential forms over smooth simplices. Here, by a smooth simplex we mean the following.
\begin{defi}\label{defi:smoothsimplices:redsp} We call an $n$-simplex $\sigma:\Delta_n\to \underline{S}_\L$ \textbf{smooth} if
\begin{itemize}
\item[(i)] it lifts to a map $\Delta_n\to S_x$ that is smooth as map into $S$, for some (equivalently, any) $x\in \L$,
\item[(ii)] it is a map of stratified spaces with respect the stratification of $\Delta_n$ by open faces and the Lerman-Sjamaar stratification of $\underline{S}_\L$. 
\end{itemize}
\end{defi}
The free abelian groups on the sets of smooth $n$-simplices form a subcomplex of the singular chain complex, which we denote as $C^{\infty-\mathrm{sing}}(\underline{S}_\L)$. As we will show below, this subcomplex still computes the singular homology. More precisely: 
\begin{prop}\label{prop:smsingchains:isohomology} The inclusion of the subcomplex of smooth singular chains gives an isomorphism
\begin{equation*} H_\bullet(C^{\infty-\mathrm{sing}}(\underline{S}_\L))\cong H^\mathrm{sing}_\bullet(\underline{S}_\L,\Z).
\end{equation*} 
\end{prop} In view of this, for each $n$, we can define the pairing
\begin{equation}\label{eqn:pairing:sjamaar:deRham} H^n_\mathrm{dR}(\underline{S}_\L)\times H^\mathrm{sing}_n(\underline{S}_\L,\Z) \to \R, \quad \big\langle \alpha, \big[\sum_{i=1}^k\lambda_i\sigma_i\big] \big\rangle=\sum_{i=1}^k\lambda_i\int_{\Delta_n}\widehat{\sigma}_i^*\widehat{\alpha},
\end{equation} where each $\sigma_i:\Delta_n\to \underline{S}_\L$ is a smooth simplex, $\widehat{\sigma}_i:\Delta_n\to S_x$ is a smooth lift of $\sigma_i$ as in Definition \ref{defi:smoothsimplices:redsp} and $\widehat{\alpha}$ is a smooth extension of $q_x^*\alpha$ to an open neighbourhood of $S_x$ in $S$. Note here that $\widehat{\sigma}_i^*\widehat{\alpha}$ does not depend on the choice of lift $\widehat{\sigma}_i$ or the choice of extension $\widehat{\alpha}$, since $\sigma_i$ restricts to a smooth map from the interior $\mathring{\Delta}_n$ of $\Delta_n$ (which is dense in $\Delta_n$) into a stratum $\underline{\Sigma}$ in $\underline{S}_\L$, and $(\widehat{\sigma}_i^*\widehat{\alpha})\vert_{\mathring{\Delta}_n}=(\sigma_i\vert_{\mathring{\Delta}_n})^*(\alpha\vert_{\underline{\Sigma}})$ (with $\alpha\vert_{\underline{\Sigma}}$ as in Remark \ref{rem:deRhamforms:non-principalstrata}). 
From the pairing \eqref{eqn:pairing:sjamaar:deRham} we obtain a map
\begin{equation}\label{eqn:explicitiso:sjamaar:deRham} H^\bullet_\mathrm{dR}(\underline{S}_\L)\to H^\bullet_\mathrm{sing}(\underline{S}_\L,\R).
\end{equation}
\begin{thm}\label{thm:explicitiso:sjamaar:deRham} The maps \eqref{isomorphism:sjamaar:deRham} and \eqref{eqn:explicitiso:sjamaar:deRham} coincide. In particular, \eqref{eqn:explicitiso:sjamaar:deRham} is an isomorphism. 
\end{thm}
In the remainder of this section we will prove Proposition \ref{prop:smsingchains:isohomology} and Theorem \ref{thm:explicitiso:sjamaar:deRham}]. The argument that we will give for Proposition \ref{prop:smsingchains:isohomology} is a variation of that in \cite{Bre}. It will use the following.
\begin{lemma}\label{lemma:smsingchains:isohomology} For each $\O\in \underline{S}_\L$ and each open neighbourhood $\underline{U}$ of $\O$ in $\underline{S}_\L$, there is open neighbourhood $\underline{W}$ of $\O$ in $\underline{U}$ together with a deformation retract $\underline{f}:\underline{W}\times [0,1]\to \underline{W}$ onto the point $\O$, with the following properties. 
\begin{itemize}
\item[(i)] For some (or equivalently any) $x\in \L$, it admits a lift
\begin{center} \begin{tikzcd} W\times [0,1]\arrow[r,"f"] \arrow[d,"q\times \mathrm{id}_{[0,1]}"'] & W \arrow[d,"q"] \\
\underline{W}\times [0,1] \arrow[r,"\underline{f}"] & \underline{W} 
\end{tikzcd} 
\end{center} by a deformation retract $f$ of an open $W$ in $S_x$ onto the $G_x$-orbit in $S_x$ corresponding to $\O$, that is smooth in the sense that there are an open $\widehat{W}$ in $S$ around $W$ and a smooth map $\widehat{f}:\widehat{W}\times [0,1]\to S$ that extends $f$.
\item[(ii)] It respects the Lerman-Sjamaar stratification of $\underline{S}_\L$, in the sense that, for each $\O'\in \underline{W}$, $\underline{f}(\{\O'\}\times [0,1[)$ is contained in the stratum through $\O'$.
\end{itemize} 
\end{lemma}
\begin{proof}[Proof of Lemma \ref{lemma:smsingchains:isohomology}] After restricting to a symplectic cross-section at a slice as in Theorem \ref{thm:symplectic:cross-section}, we can assume that $\L$ is the origin in $\g^*$. Moreover, appealing to the Marle-Guillemin-Sternberg normal form the proof reduces to the case in which $J$ is the Marle-Guillemin-Sternberg local model and $\O$ is the orbit through $[1,0,0]\in G\times _H(\h^0\oplus V)$ (using the same notation as in the proof of Proposition \ref{prop:images:canhamstrata:affineopen}). In this model, the fiber of $J$ over the origin in $\g^*$ is a $G$-invariant open neighbourhood of $G\times_HJ_V^{-1}(0)$ around $[1,0,0]$. Since the origin in $J_V^{-1}(0)$ admits arbitrarily small $H$-invariant open neighbourhoods that are invariant under scaling by $[0,1]$ and since $H_{tv}=H_v$ for all $v\in V$ and $t>0$, the desired deformation retract $\underline{f}$ can be obtained by restricting the $G$-equivariant deformation retract 
\begin{equation*} G\times _H(\h^0\oplus V)\times [0,1]\to G\times _H(\h^0\oplus V),\quad ([g,\alpha,v],t)\mapsto [g,(1-t)\alpha,(1-t)v],
\end{equation*} and descending to the level of orbit spaces.
\end{proof}
\begin{proof}[Proof of Proposition \ref{prop:smsingchains:isohomology} and Theorem \ref{thm:explicitiso:sjamaar:deRham}] First, we give a proof of Proposition \ref{prop:smsingchains:isohomology}. It is enough to show that the map $i^*$ in cohomology (rather than homology) induced by the inclusion map $i:C^{\infty-\mathrm{sing}}(\underline{S}_\L)\hookrightarrow C^\mathrm{sing}(\underline{S}_\L)$ is an isomorphism. Indeed, since $C(i):=C^\mathrm{sing}(\underline{S}_\L)/C^\mathrm{\infty-sing}(\underline{S}_\L)$ is a chain complex of free abelian groups, it would follow from this and the universal coefficient theorem that, for each $n$, $\mathrm{Hom}(H_n(C(i)),\Z)=0$ and $\mathrm{Ext}^1(H_n(C(i)),\Z)=0$, and so $H_n(C(i))=0$. \\

To prove that $i^*$ is an isomorphism, for $\underline{U}$ open in $\underline{S}_\L$, let $C^n_{\infty-\mathrm{sing}}(\underline{U})$ denote the abelian group of $\Z$-valued functions on the set of smooth $n$-simplices in $\underline{S}_\L$ with image in $\underline{U}$. These form a complex of pre-sheaves $C_{\infty-\mathrm{sing}}$ on $\underline{S}_\L$, with the natural restriction maps, and with differential inherited from that of the complex of singular chains on $\underline{S}_\L$ by restricting and dualizing. Let $\mathcal{C}$ denote the complex of sheaves on $\underline{S}_\L$ obtained by sheafification. Since the composition of sheafification with taking the stalk at a point $\O\in \underline{S}_\L$ is an exact functor, it commutes with taking cohomology. Therefore, $H^n(\mathcal{C}_\O)$ is isomorphic to the stalk of the pre-sheaf $H^n(C_{\infty-\mathrm{sing}})$ at $\O$. Note that, for any $\underline{V}$ around $\O$ as in Lemma \ref{lemma:smsingchains:isohomology}, 
\begin{equation*} H^n(C_{\infty-\mathrm{sing}})(\underline{V})=H^n(C_{\infty-\mathrm{sing}}(\underline{V}))=\begin{cases} \Z&\textrm{ if }n=0,\\
0&\textrm{ if }n>0.
\end{cases}
\end{equation*} Indeed, this follows from the observation that, for $\underline{f}$ as in the lemma, $\underline{f}(-,0)=\mathrm{id}_{\underline{V}}$ and the constant map $\underline{f}(-,1)=\O$ induce the same map in cohomology, because the chain homotopy on singular chains that is induced by $\underline{f}$ preserves the subcomplex of smooth singular chains. So, $H^0(\mathcal{C}_\O)=\Z$ and $H^n(\mathcal{C}_\O)=0$ for all $n>0$. This means that $\mathcal{C}^\bullet$ gives a resolution 
\begin{equation}\label{eqn0:pf:prop:smsingchains:isohomology} 0\to \underline{\Z}\to \mathcal{C}^0\xrightarrow{\d} \mathcal{C}^1\xrightarrow{\d} ...
\end{equation} of the constant sheaf $\underline{\Z}$. To prove that $i^*$ is an isomorphism, we will now show that this resolution is acyclic and that the canonical map of complexes
\begin{equation}\label{eqn1:pf:prop:smsingchains:isohomology} C_{\infty-\mathrm{sing}}(\underline{S}_\L)\to \mathcal{C}(\underline{S}_\L)
\end{equation} is a quasi-isomorphism. This would imply that, from the canonical map of resolutions from \eqref{eqn0:pf:prop:smsingchains:isohomology} into the usual acyclic resolution of $\underline{\Z}$ by singular chains, we obtain a commutative triangle
\begin{center} \begin{tikzcd} 
H^\bullet(C_{\infty-\mathrm{sing}}(\underline{S}_\L)) \arrow[dr,"\sim"'] \arrow[rr,"i^*"] & & H^\bullet_\mathrm{sing}(X,\Z) \arrow[dl,"\sim"] \\
& H^\bullet(X,\underline{\Z}) &
\end{tikzcd} 
\end{center} in which the horizontal map is that induced by the inclusion and the other maps are isomorphisms, and so it would follow that $i^*$ is indeed an isomorphism. Since the pre-sheaves $\mathcal{C}^n\vert_{\underline{U}}$ satisfy the glueing axiom, it follows from \cite[Theorem 3.9.1]{Gode} that the canonical map $C_{\infty-\mathrm{sing}}(\underline{U})\to \mathcal{C}(\underline{U})$ is term-wise surjective, for every open $\underline{U}$ in $\underline{S}_\L$. So, since the pre-sheaves $C^n_{\infty-\mathrm{sing}}$ are flasque, it follows that the sheaves $\mathcal{C}^n$ are flasque as well and, hence, acyclic. Moreover, denoting by $K^\bullet(\underline{S}_\L)$ the kernel of \eqref{eqn1:pf:prop:smsingchains:isohomology}, it follows that we have a short exact sequence of complexes
\begin{equation*} 
0\to K(\underline{S}_\L)\to C_{\infty-\mathrm{sing}}(\underline{S}_\L)\to \mathcal{C}(\underline{S}_\L)\to 0.
\end{equation*} So, it is remains to show that $K(\underline{S}_\L)$ is acyclic. Given an open cover $\U$ of $\underline{S}_\L$, write $C^\U(\underline{S}_\L)$ for the subcomplex of $C^{\infty-\mathrm{sing}}(\underline{S}_\L)$ generated by the smooth simplices with image contained in opens belonging to $\U$, let $C_\U(\underline{S}_\L)$ be the dual complex and $K_\U(\underline{S}_\L)$ the kernel of the canonical map 
\begin{equation}\label{eqn2:pf:prop:smsingchains:isohomology}
 C_{\infty-\mathrm{sing}}(\underline{S}_\L)\to C_\U(\underline{S}_\L).
\end{equation} The small simplices theorem holds for smooth singular chains, since the barycentric subdivision operator on singular chains preserves the subcomplex of smooth singular chains and the same goes for the usual chain homotopy between this operator and the identity map. So, for each $\U$, \eqref{eqn2:pf:prop:smsingchains:isohomology} is a quasi-isomorphism and, hence, $H^\bullet(K_\U(\underline{S}_\L),\d)=0$. Since $K(\underline{S}_\L)$ is the direct limit of the complexes $K_\U(\underline{S}_\L)$ and taking cohomology commutes with direct limits, it follows that $K(\underline{S}_\L)$ is acyclic. So, indeed, \eqref{eqn0:pf:prop:smsingchains:isohomology} is acyclic and \eqref{eqn1:pf:prop:smsingchains:isohomology} is a quasi-isomorphism, which proves Proposition \ref{prop:smsingchains:isohomology}. \\

For the proof of Theorem \ref{thm:explicitiso:sjamaar:deRham}, note that the same formula as in \eqref{eqn:pairing:sjamaar:deRham} defines a map of acyclic resolutions $\varphi:\Omega_{\underline{S}_\L}\to \mathcal{C}\otimes\underline{\R}$ covering the identity on $\underline{\R}$. From this and the canonical map of acyclic resolutions from $\mathcal{C}\otimes\underline{\R}$ into the usual acyclic resolution of $\underline{\R}$ by singular chains, we obtain the commutative diagram
\begin{center} 
\begin{tikzcd} 
H^\bullet_\mathrm{dR}(\underline{S}_\L) \arrow[dr,"\sim"'] \arrow[r,"(\ref{eqn1:pf:prop:smsingchains:isohomology})_*^{-1}\circ \varphi_*"] & H^\bullet(C_{\infty-\mathrm{sing}}(\underline{S}_\L),\R) \arrow[d,"\sim" {anchor=south, rotate=90}]\arrow[r,"i^*"] & H^\bullet_\mathrm{sing}(X,\R) \arrow[dl,"\sim"]\\
& H^\bullet(X,\underline{\R}) &
\end{tikzcd} 
\end{center} in which the non-horizontal maps are the isomorphisms induced by the three acyclic resolutions. Since the map \eqref{eqn:explicitiso:sjamaar:deRham} is the composition of the two horizontal arrows, Theorem \ref{thm:explicitiso:sjamaar:deRham} follows.
\end{proof}
\subsection{The linear variation theorem}
\subsubsection{On the meaning of linear variation}\label{subsec:meaning:linvar} 
Before turning to its proof, we elaborate on the meaning of Theorem \ref{thm:linvar:main} (the linear variation theorem). Consider $G$ and $J:(S,\omega)\to \g^*$ as in the introduction. Let $\sigma_+$ be a stratum of the stratification ${\sS}_\mathrm{Ham}(\Delta_+)$ in Theorem \ref{thm:stratification:image:momentum:map}. As argued in the proof of Corollary \ref{cor:constructibility:pushforwards}, the reduced spaces $\underline{S}_\L$ have finite-dimensional cohomology groups. So, since for $x\in \sigma_+$ the reduced spaces $\underline{S}_{\L_x}$ are the fibers of the topologically locally trivial fibration \eqref{eqn:restrictionstratum:transversemommap}, for each degree $n$ their real cohomology groups naturally form a flat vector bundle
\begin{equation}\label{eqn:vectorbundle:cohomgps}
\bigsqcup_{x\in \sigma_+} H^n(\underline{S}_{\L_x})\to \sigma_+.
\end{equation} The flat connection $\nabla$ on this vector bundle is the so-called Gauss-Manin connection. Theorem \ref{thm:linvar:main} states that the cohomology classes $\varpi_{x}:=\varpi_{\L_x}\in H^2(\underline{S}_{\L_x})$ of the reduced symplectic forms (as in Example \ref{ex:cohomclass:sympform:redsp}) vary linearly with $x\in \sigma_+$, in the sense that if $x_0\in \sigma_+$ and $W_+$ is a convex open in $\sigma_+$ around $x_0$, then the difference
\begin{equation*} P^\mathrm{GM}_{[x,x_0]}(\varpi_{x})-\varpi_{x_0}\in H^2(\underline{S}_{\L_{x_0}})
\end{equation*} depends linearly on $x-x_0$, where 
\begin{equation*} P^\mathrm{GM}_{[x,x_0]}:H^2(\underline{S}_{\L_x})\to H^2(\underline{S}_{\L_{x_0}})
\end{equation*} denotes the parallel transport of the Gauss-Manin connection along any path in $W_+$ from $x$ to $x_0$.
\subsubsection{Proof of linear variation} For the proof of Theorem \ref{thm:linvar:main} an alternative characterization of linear variation will be useful. To state this, let ($\H_{\sigma_+},\nabla)$ denote the flat vector bundle over $\sigma_+$ formed by the degree $2$ cohomology groups of the reduced spaces, as in  \eqref{eqn:vectorbundle:cohomgps}, and let $\varpi$ denote the section of $\H_{\sigma_+}$ given by the cohomology classes $\varpi_{x}\in H^2(\underline{S}_{\L_x})$.  
\begin{prop}\label{prop:redsympforms:smoothsection} The section $\varpi$ of $\H_{\sigma_+}$ is smooth. 
\end{prop}
\begin{proof} Let $x\in \sigma_+$ and consider a local trivialization $(\Phi,\varphi)$ as in Theorem \ref{thm:localtriviality} over an open $W_+$ in $\sigma_+$ around $x$. Fix a basis $\mathrm{c}_1,...,\mathrm{c}_n$ of $H_2(\underline{S}_{\L_x})$. Via the local trivialization of $\mathcal{H}_{\sigma_+}$ induced by $\Phi$, the basis $\check{\mathrm{c}}_1,...,\check{\mathrm{c}}_n$ of $H^2(\underline{S}_{\L_x})$ dual to this extends to a flat local frame $e$ over $W_+$ given by
\begin{equation*} e_i(y)=(\underline{\Phi}\vert_{\underline{S}_{\L_y}})^*(\check{\mathrm{c}}_i)\in H^2(\underline{S}_{\L_y}), \quad y\in W_+, \quad i=1,...,n.
\end{equation*} So, to prove smoothness of $\varpi$, we ought to show that $\big\langle \varpi_{y}, (\underline{\Phi}\vert_{\underline{S}_{\L_y}}^{-1})_*(\mathrm{c}_i) \big\rangle$ depends smoothly on $y$. By Proposition \ref{prop:smsingchains:isohomology} we can write 
\begin{equation*} \mathrm{c}_i=\big[\sum_{j=1}^{k_i}\lambda_{i,j}\sigma_{i,j}\big],
\end{equation*} for $\sigma_{i,j}$ smooth simplices in $\underline{S}_{\L_x}$, and in view of Theorem \ref{thm:explicitiso:sjamaar:deRham} it holds, that for each $y\in W_+$,
\begin{equation}\label{eqn:coefficients:redspforms:locfr} \big\langle \varpi_{y}, (\underline{\Phi}\vert_{\underline{S}_{\L_y}}^{-1})_*(\mathrm{c}_i) \big\rangle=\sum_{j=1}^{k_i}\lambda_{i,j}\int_{\Delta_2} (\widehat{\sigma}_{i,j})^*(i_y)^*(\widehat{\Phi}^{-1})^*\omega, 
\end{equation} where $\widehat{\Phi}^{-1}$ is a smooth extension of $\Phi^{-1}$ to an open of the form $V\times W_+$, with $V$ open in $S$ around $S_{\L_x}$, $\widehat{\sigma}_{i,j}:\Delta_2\to V$ is a smooth lift of $\sigma_{i,j}$ as in Definition \ref{defi:smoothsimplices:redsp} and $i_y:V\to V\times W_+$ is given by $p\mapsto (p,y)$. This indeed depends smoothly on $y$, because the integrands do so.
\end{proof}
Having established this, we can state the aforementioned characterization of linear variation.
\begin{lemma}\label{lemma:altcharlinvar} The cohomology classes $\varpi_{x}\in H^2(\underline{S}_{\L_x})$ vary linearly with $x\in \sigma_+$ in the sense explained in Subsection \ref{subsec:meaning:linvar} if and only if for every constant vector field $X$ on $\sigma_+$ the section $\nabla_X\varpi$ of $\H_{\sigma_+}$ is flat with respect to $\nabla$. 
\end{lemma}
\begin{proof} If $W_+$ is a convex open in $\sigma_+$, then there is a flat local frame $e$ of $\H_{\sigma_+}$ defined on $W_+$ and, writing 
\begin{equation*} \varpi\vert_{W_+}=\sum_{i=1}^nf_ie_i,\quad f_1,...,f_n\in C^\infty(W_+),
\end{equation*} it holds that
\begin{equation}\label{eqn:charoflinvarprop:1}
P^\mathrm{GM}_{[x,x_0]}(\varpi_{x})=\sum_{i=1}^nf_i(x)e_i(x_0)
\end{equation} for any $x,x_0\in W_+$. Moreover, for any vector field $X$ on $\sigma_+$,
\begin{equation}\label{eqn:charoflinvarprop:2}
\nabla_X\varpi\vert_{W_+}=\sum_{i=1}^n(\L_Xf_i)\,e_i.
\end{equation}
In view of \eqref{eqn:charoflinvarprop:2}, $\nabla_X\varpi$ is flat if and only if for any $W_+$ and $e$ as above the functions $\L_Xf_1,...,\L_Xf_n$ are constant. So, this holds for every  constant vector field $X$ on $\sigma_+$ if and only if for any such $W_+$ and $e$ the functions $f_1,...,f_n$ are restrictions of affine functions on the affine hull of $\sigma_+$. In view of \eqref{eqn:charoflinvarprop:1}, the lemma follows from this. 
\end{proof}

\begin{proof}[Proof of Theorem \ref{thm:linvar:main}] In the setting of Theorem \ref{thm:stratification:image:momentum:map}, let $\sigma_+ \in \sS_\mathrm{Ham}(\Delta_+)$. By Lemma \ref{lemma:altcharlinvar}, it is enough to show that $\nabla_{X}\varpi$ is flat for any constant vector field $X$ on $\sigma_+$. Suppose we are given such an $X$ and let $x\in \sigma_+$. Fix a local frame of $X_1,...,X_m$ of $\sigma_+$ consisting of constant vector fields, with $X_m:=X$. For each $l$, let $(\widehat{X}_l,\widehat{X}^\sigma_l)$ be a pair of lifts of $X_l$ constructed as in the proof of Lemma \ref{lemma:liftingvecfields:stratumwise}, and set $(\widehat{X},\widehat{X}^\sigma):=(\widehat{X}_m,\widehat{X}^\sigma_m)$. Let $(\Phi,\varphi)$ be the associated trivialization over a convex open $W_+$ in $\sigma_+$ around $x$ and let $\widehat{\Phi}^{-1}$ be the associated equivariant smooth extension of $\Phi^{-1}$, all constructed as in the proof of Theorem \ref{thm:localtriviality} using the flows of pairs of lifts $(\widehat{X}_1,\widehat{X}_1^\sigma)$,...,$(\widehat{X}_m,\widehat{X}^\sigma_m)$. Consider $\mathrm{c}_1,...,\mathrm{c}_n$ and $e$ as in the proof of Proposition \ref{prop:redsympforms:smoothsection}. Then
\begin{equation*} \varpi\vert_{W_+}=\sum_{i=1}^nf_ie_i,
\end{equation*} with $f_i$ given by
\begin{equation*} f_i(y)=\big\langle \varpi_{y}, ({\underline{\Phi}\vert_{\underline{S}_{\L_y}}^{-1}})_*(\mathrm{c}_i) \big\rangle, \quad y\in W_+. 
\end{equation*} 
For any $y\in W_+$, the $2$-form $\L_{\widehat{X}}\omega$ descends to a closed differential $2$-form on the reduced space $\underline{S}_{y}$ at the point $y$, that we denote by $\L_{\widehat{X}}\omega_{\underline{S}_{y}}$ and is determined by the property that 
\begin{equation*} i^*_{y}\L_{\widehat{X}}\omega=q_{y}^*(\L_{\widehat{X}}\omega_{\underline{S}_{y}}).
\end{equation*} Let $\L_{\widehat{X}}\omega_{\underline{S}_{\L_y}}$ denote the corresponding $2$-form on $\underline{S}_{\L_y}$. Taking the Lie derivative of the expression \eqref{eqn:coefficients:redspforms:locfr} we find that, for any $i=1,...,n$ and $y\in W_+$,
\begin{align*} (\L_{X}f_i)(y)&=\sum_{j=1}^{k_i}\lambda_{i,j}\int_{\Delta_2} (\widehat{\sigma}_{i,j})^*(i_y)^*\L_{(0,X)}(\widehat{\Phi}^{-1})^*\omega\\
&=\sum_{j=1}^{k_i}\lambda_{i,j}\int_{\Delta_2} (\widehat{\sigma}_{i,j})^*(i_y)^*(\widehat{\Phi}^{-1})^*\L_{\widehat{X}}\omega=\big\langle [\L_{\widehat{X}}\omega_{\underline{S}_{\L_y}}], ({\underline{\Phi}\vert_{\underline{S}_{\L_y}}^{-1}})_*(\mathrm{c}_i) \big\rangle
\end{align*} and, hence,
\begin{equation*} (\nabla_{X}\varpi)(y)=\sum_{i=1}^n(\L_{X}f_i)(y)\,e_i(y)=[\L_{\widehat{X}}\omega_{\underline{S}_{\L_y}}].
\end{equation*} For $y\in W_+$, let $\Phi_{\L_y}:S_{\L_y}\to S_{\L_x}$ be the restriction of $\Phi$. Note that pull-back along $(\underline{\Phi}^{-1}_{\L_y})^*$ takes the de Rham complex of $\underline{S}_{\L_y}$ into that of $\underline{S}_{\L_x}$, and parallel transport of the Gauss-Manin connection along a path in $W_+$ from $x$ to $y$ takes $[\L_{\widehat{X}}\omega_{\underline{S}_{\L_y}}]$ to $(\underline{\Phi}^{-1}_y)^*[\L_{\widehat{X}}\omega_{\underline{S}_{\L_y}}]$. To prove that $\nabla_{X}\varpi\vert_{W_+}$ is flat, it therefore suffices to show that 
\begin{equation}\label{eqn:pfthm:linvar:main} (\underline{\Phi}^{-1}_{\L_y})^*[\L_{\widehat{X}}\omega_{\underline{S}_{\L_y}}]=[\L_{\widehat{X}}\omega_{\underline{S}_{\L_x}}]
\end{equation} for all $y\in W_+$. Given such $y$, consider the $1$-form
\begin{equation*} \alpha:=(\Phi_y^{-1})^*i_{y}^*(\iota_{\widehat{X}}\omega)-i_{x}^*(\iota_{\widehat{X}}\omega)
\end{equation*} on $S_{x}^\mathrm{prin}$, where $\Phi_y:S^\textrm{prin}_y\to S^\textrm{prin}_x$ is the restriction of $\Phi_{\L_y}$. Note here that, by construction, $\widehat{X}^\sigma_l$ restricts to $X_l$ on $\sigma_+$ for each $l$ and, hence, $\varphi$ restricts to the map $y\mapsto (x,y)$ on $W_+$, which implies that $\Phi_{\L_y}$ indeed maps $S_y$ onto $S_x$. The $1$-form $\alpha$ extends smoothly to an open neighbourhood of $S_{x}$ in $S$ because $\Phi^{-1}$ extends smoothly and it is $G_x$-invariant since $\Phi$ is equivariant and since $\omega$ and $\widehat{X}$ are $G$-invariant. Moreover, it is $G_x$-horizontal, for if $\xi\in \g_x$ and $p\in S_{x}^\mathrm{prin}$, then using the momentum map condition we find that 
\begin{align*} (\iota_{a(\xi)}\alpha)_p&=\omega(a(\xi),\widehat{X})_p-\omega(a(\xi),\widehat{X})_{\Phi_y^{-1}(p)}\\
&=\langle \d J(\widehat{X})_p,\xi \rangle-\langle \d J(\widehat{X})_{\Phi_y^{-1}(p)},\xi\rangle\\
&=\langle \widehat{X}^\sigma_{x},\xi \rangle-\langle \widehat{X}^\sigma_{y},\xi\rangle,
\end{align*} which is zero because $\widehat{X}^\sigma$ restricts to a constant vector field on $\sigma_+$. So, $\alpha$ defines a $1$-form on $\underline{S}_x$. The corresponding $1$-form on $\underline{S}_{\L_x}$ is a primitive of $(\underline{\Phi}_{\L_y}^{-1})^*(\L_{\widehat{X}}\underline{\omega}_{S_{\L_y}})-(\L_{\widehat{X}}\underline{\omega}_{S_{\L_x}})$ and, hence, \eqref{eqn:pfthm:linvar:main} indeed holds. This proves the theorem. 
\end{proof}
\subsubsection{Polynomiality of symplectic volumes} Finally, we deduce the remaining corollary.

\begin{proof}[Proof of Corollary \ref{cor:sympvolpolynomial}] By definition, the dimension of $\underline{S}_{\L}$ is the maximum of the dimensions $n_1,...,n_k$ of the connected components $C^1_\L,...,C^k_\L$ of $\underline{S}_\L^\mathrm{prin}$. These dimensions are constant for $\L\in \underline{\sigma}$, since $\underline{J}$ is locally trivial (in the stratified sense of Remark \ref{rem:trivstratpreserving}) over the connected space $\underline{\sigma}$. To prove polynomiality of the symplectic volumes, let $x_0\in \sigma_+$, $\L_0:=\L_{x_0}$ and consider a trivialization $\Phi$ as in Theorem \ref{thm:localtriviality}, defined on a convex open $W_+$ around $x_0$ in $\sigma_+$. In view of Remark \ref{rem:trivstratpreserving}, the restriction $\underline{\Phi}_{\L}:\underline{S}_{\L}\to \underline{S}_{\L_{0}}$ of $\underline{\Phi}$ restricts to a diffeomorphism between $\underline{S}_{\L}^\mathrm{prin}$ and $\underline{S}_{\L_0}^\mathrm{prin}$. The resulting diffeomorphism $\underline{\Phi}_\L$ is orientation preserving with respect to the orientations on $\underline{S}_{\L_{0}}^\mathrm{prin}$ and $\underline{S}_{\L}^\mathrm{prin}$ induced by $\underline{\omega}_{S_{\L_{0}}}$ and $\omega_{\underline{S}_\L}$, since $(\underline{\Phi}^{-1}_{\L_0})^*\underline{\omega}_{S_{\L_{0}}}=\underline{\omega}_{S_{\L_{0}}}$, since $W_+$ is connected and since the family of forms $(\Phi^{-1}_{\L_x})^*\underline{\omega}_{S_{\L_x}}$ on $\underline{S}_{\L_{0}}^\mathrm{prin}$ depends smoothly on $x\in W_+$. Therefore,
\begin{equation*}
\mathrm{Vol}(\underline{S}_{\L},\underline{\omega}_{S_{\L}})=\sum_{i=1}^k\int_{C^i_{\L}} \frac{1}{n_i!}\underline{\omega}_{S_{\L}}^{n_i}=\sum_{i=1}^k\int_{C^i_{\L_{0}}}\frac{1}{n_i!}(\Phi_\L^{-1})^*\underline{\omega}_{S_{\L}}^{n_i}.
\end{equation*} Using parts (ii) and (iii) of Theorem \ref{thm:sjamaar:finvol} and the fact that, by the linear variation theorem,
\begin{equation*} (\underline{\Phi}_{\L_x}^{-1})^*\varpi_{x}=\sum_{i=1}^{\dim(\sigma_+)}(x^i-x_0^i)\mathrm{c}_i+\varpi_{x_0}
\end{equation*} for some $\mathrm{c}_i\in H_{\mathrm{dR}}^2(\underline{S}_{\L_0})$, the corollary follows.
\end{proof} 

\section{Beyond actions of connected compact Lie groups}
\subsection{Extension to actions of arbitrary compact Lie groups}\label{sec:extension:disconnectedgps} 
\subsubsection{Extension of the main results}\label{sec:extensionpf:disconnectedgps} Next, we extend our main results to actions of possibly disconnected compact Lie groups. In a nutshell, the main difference is that in this generality the coadjoint orbit space $\underline{\g^*}$ need no longer be parametrized by $\t^*_+$. Nevertheless, $\underline{\g^*}$ still has a natural integral affine stratification (the stratification by connected components of $G$-orbit types) that corresponds to the stratification of $\t^*_+$ by its open faces when $G$ is connected and allows to extend our results. \\

Turning to the details, let $G$ be a compact Lie group and let $J:(S,\omega)\to \g^*$ be a Hamiltonian $G$-space with the property that the momentum map is proper as map into its image $\Delta$. Further, fix a maximal torus $T$ in $G$ and a choice of closed Weyl chamber $\t^*_+\subset \t^*$. Given $x \in \g^*$, we denote by $G_x$ the isotropy group of the coadjoint $G$-action and by $\g_x$ its Lie algebra. 
As when $G$ is connected (see \S \ref{sec:background:Lie}), we identify $\t^*$ and $\g^*_x$ with $(\g^*)^T$ and its $G_x$-saturation in $\g^*$, for $x\in \t^*$. We let $G^0$ denote the identity component of $G$ and $\Gamma:=G/G^0$ its group of connected components. Moreover, we denote by 
\begin{equation*} \pi:\g^*\to \underline{\g^*},\quad \pi_0:\g^*/G^0\to\underline{\g^*}, \quad\quad (\underline{\g^*}:= \g^*/G)
\end{equation*} the quotient maps for the coadjoint $G$-action and the residual $\Gamma$-action on the orbit space of the coadjoint $G^0$-action.
The inclusion $\t^*_+\subset \g^*$ descends to a homeomorphism between $\t^*_+$ and $\g^*/G^0$, whose composition with $\pi_0$ we denote as
\begin{equation}\label{eqn:coadorbproj:restrictionWeylchamber} \pi_+:\t^*_+\to \underline{\g^*}.
\end{equation} By the transverse momentum map we mean the map 
\begin{equation*} \underline{J}:\underline{S}\to \underline{\g^*}
\end{equation*} from the orbit space $\underline{S}:=S/G$ into $\underline{\g^*}$, induced by 
\begin{equation*} \pi\circ J:S\to \underline{\g^*}.
\end{equation*} These both have image $\underline{\Delta}$. Given a subset $\underline{X}\subset \underline{\g^*}$ with corresponding $G$-invariant subset $X\subset \g^*$, we use the notation
\begin{equation*} X_+:=(\pi_+)^{-1}(\underline{X})=X\cap \t^*_+.
\end{equation*} Finally, we let $\sS_G(\underline{\g^*})$ denote the stratification of $\underline{\g^*}$ by connected components of $G$-orbit types. Theorem \ref{thm:stratification:image:momentum:map} extends as follows.
\begin{thm}\label{thm:stratification:image:momentum:map:disconnectedgps} The subspace $\underline{\Delta}$ of $\underline{\g^*}$ admits a natural 
integral affine stratification ${\sS}_\mathrm{Ham}(\underline{\Delta})$ with the property that, for each $\underline{\sigma}\in {\sS}_\mathrm{Ham}(\underline{\Delta})$, the restriction
\begin{equation}\label{eqn:restrictionstratum:mommap:disconnectedgps} 
J:J^{-1}(\sigma)\to \sigma
\end{equation} is a $G$-equivariantly locally trivial fibration (in the sense of Theorem \ref{thm:localtriviality}, but with the role of $\t^*_+$ taken over by $\underline{\g^*}$, meaning that $\Delta_+$ is replaced by $\underline{\Delta}$, $\sigma$ by $\underline{\sigma}$, $x$ and $\L_x$ by $\L$, and $W_+$ by $\underline{W}$). This is uniquely 
determined by the property that each stratum is an integral affine submanifold of a stratum of $\sS_G(\underline{\g^*})$ and that for any $x\in \Delta_+$ the tangent space at $\pi_+(x)$ to the stratum $\underline{\sigma}$ through that point is given by
\begin{equation}\label{eqn:definingproperty:stratofimage:disconnectedgps}
T_{\pi_+(x)}\underline{\sigma}=\mathrm{d}(\pi_+\vert_{\Sigma_+})_x\left(\bigcap_{p\in J^{-1}(x)} (\g_p^0)^{G_x}\right)\subset T_{\pi_+(x)}\underline{\Sigma}
\end{equation}  where 
\begin{itemize} 
\item $\underline{\Sigma}\in \sS_G(\underline{\g^*})$ is the stratum through $\pi_+(x)$,
\item $\g_p$ is the isotropy Lie algebra at $p\in S$ of the $G$-action,
\item $G_x$ is the isotropy group at $x\in \g^*$ of the coadjoint $G$-action and $\g_x$ is its Lie algebra,
\item $\g_p^0\subset \g_x^*$ is the annihilator of $\g_p$ in $\g_x$, 
\item $(\g_p^0)^{G_x}:=\g_p^0\cap (\g_x^*)^{G_x}$ is the linear subspace of $\t^*$ given by the intersection of $\g_p^0$ with the fixed point set of the coadjoint $G_x$-action.
\end{itemize}  
\end{thm}
In the part on the defining property of ${\sS}_\mathrm{Ham}(\underline{\Delta})$, the statement now makes reference to an integral affine structure on each stratum of $\sS_G(\underline{\g^*})$. These structures arise as follows. 
\begin{prop}\label{prop:affstr:coadjoint-orbit-types} For each $\underline{\Sigma}\in \sS_G(\underline{\g^*})$, $\Sigma_+$ is an integral affine submanifold of $(\t^*,\Lambda_T)$, $\pi_+$ restricts to a surjective local diffeomorphism $\Sigma_+\to \underline{\Sigma}$ and the integral affine structure on $\Sigma_+$ descends to $\underline{\Sigma}$, meaning that $\underline{\Sigma}$ admits a unique integral affine structure making $\pi_+:\Sigma_+\to \underline{\Sigma}$ a local integral affine isomorphism. 
\end{prop}
When $G$ is connected, these integral affine structures are inherited from that of the integral affine vector space $(\t^*,\Lambda_T)$ via the canonical identification of these strata with open faces of $\t^*_+$. Therefore, Theorem \ref{thm:stratification:image:momentum:map} indeed extends Theorem \ref{thm:stratification:image:momentum:map:disconnectedgps}. Proofs of Theorem \ref{thm:stratification:image:momentum:map:disconnectedgps} and Proposition \ref{prop:affstr:coadjoint-orbit-types}, as well as a more intrinsic description of the integral affine structures on the strata of $\sS_G(\underline{\g^*})$ will be given in the next sections. \\

Moving on to the remaining results: Corollary \ref{cor:thm:stratification:image:transverse:momentum:map} and Remark \ref{rem:trivstratpreserving} go through after replacing $J_+$ by $\pi\circ J$, $\underline{J}_+$ by $\underline{J}$, $\Delta_+$ by $\underline{\Delta}$, $\sigma_+$ by $\underline{\sigma}$, $x$ and $\L_x$ by $\L$, and $W_+$ by $\underline{W}$. Moreover, Corollary \ref{cor:constructibility:pushforwards} and its proof go through with the same replacements. Proposition \ref{prop:specializationmaps} extends as follows. 
\begin{prop}\label{prop:specializationmaps:disconngps} Let $x,y\in \Delta$ be regular values of $J$ with $\L_x$ in an open stratum of $\sS_\G(\underline{\g^*})$, and suppose that $\L_x$ and $\L_y$ are connected by a path $\gamma$ in $\underline{\Delta}$ that remains in the stratum $\underline{\sigma}\in \sS_\mathrm{Ham}(\underline{\Delta})$ through $\L_x$, except for at its end-point $\L_y$. Then the same conclusion as in Proposition \ref{prop:specializationmaps} holds. 
\end{prop}
The proof of Proposition \ref{prop:specializationmaps} goes through to prove the above extension, using that the following lemma, which we applied in the proof for connected $G$, holds for disconnected $G$ as well.
\begin{lemma} Let $G$ be a (possibly disconnected) compact Lie group and $T$ a maximal torus in $G$. Then $N_G(T)/T$ is finite and $G/N_G(T)$ is connected and has trivial cohomology in degrees $n>0$.
\end{lemma} 
\begin{proof} Both statements can be reduced to the case in which $G$ is connected. First, that $N_G(T)/T$ is finite follows from the observation that we have a (canonical) short exact sequence of groups
\begin{equation*} 1\to N_{G^0}(T)/T\to N_G(T)/T\to G/G^0\to 1,
\end{equation*} and the fact that both $N_{G^0}(T)/T$ and $G/G^0$ are finite (for surjectivity of $N_G(T)/T\to G/G^0$, recall that, in fact, for any choice of closed Weyl chamber $\t^*_+$ in $\t^*$ the inclusion $N_G(\t^*_+)\hookrightarrow N_G(T)\hookrightarrow G$ induces an isomorphism between $N_G(\t^*_+)/T$ and $G/G^0$; cf. \cite[Proposition 3.15.1]{DuKo}). Second, since any two maximal tori are conjugate in $G^0$, the map
\begin{equation}\label{eqn:lemma:G/N:trivcohom:disconnectedgps} G^0/N_{G^0}(T)\to G/N_G(T)
\end{equation} induced by the inclusion $G^0\hookrightarrow G$ is surjective. So, since \eqref{eqn:lemma:G/N:trivcohom:disconnectedgps} is also injective and a local diffeomorphism, it is in fact a diffeomorphism. Therefore, the statement about $G/N_G(T)$ reduces to the case in which $G$ is connected as well.
\end{proof}
Finally, with the same replacements as above, the statements of Theorem \ref{thm:linvar:main} and Corollary \ref{cor:sympvolpolynomial} remain true in the disconnected case, where the variation being linear and the symplectic volume function being polynomial is now meant locally, in some (or equivalently any) affine coordinate chart for the stratum $\underline{\sigma}\in \sS_\mathrm{Ham}(\underline{\Delta})$ under consideration. Proofs of these claims will also be given in the next sections.  
\subsubsection{Integral affine geometry of the coadjoint orbit type strata}\label{sec:iastr:orbittypes:coadjointact} We now turn to the proof of Proposition \ref{prop:affstr:coadjoint-orbit-types}, for which we will use the following.
\begin{lemma}\label{lemma:orbit-types:dis-con-gps} Consider a manifold $M$ equipped with a proper, smooth action of a Lie group $G$. Let $G^0$ be the identity component of $G$ and $\Gamma:=G/G^0$ its group of connected components. The partition of $\sS_G(\underline{M})$ of $\underline{M}:=M/G$ by connected components of $G$-orbit types coincides with the partition consisting of, for each $G^0$-orbit type $\Sigma_0$ in $M/G^0$, the connected components of the $\Gamma$-orbit types in $\Sigma_0/\Gamma$ of the residual $\Gamma$-action on $\Sigma_0$ (viewed canonically as subsets of $\underline{M}$).
\end{lemma}
\begin{proof} By \cite[Lemma 6]{CrMe} it suffices to show that for any $x\in M$ the members of the two partitions through the $G$-orbit $\O_x$ containing $x$ have the same germ in $\underline{M}$ at $\O_x$. Since the $G$-action is proper, there is a slice $S$ at $x$ for the $G$-action, meaning $S$ is a $G_x$-invariant submanifold of $M$ for which the map $G\times_{G_x}S\to M$, $[g,y]\mapsto g\cdot y$, is an embedding onto a $G$-invariant open $U$ around $\O_x$. This map canonically induces identifications
\begin{equation*} \underline{U}\cong S/G_x,\quad U/G^0\cong \Gamma\times_{\Gamma_x}(S/G^0_x),\quad\mathrm{where}\quad \Gamma_x:=G_x/G^0_x, \quad G^0_x:=G^0\cap G_x. 
\end{equation*}
Using that for a compact Lie group any closed subgroup isomorphic to it must be equal to it, it follows that under these identifications the $G$-orbit type through $\O_x$ becomes the $G_x$-fixed point set $S^{G_x}\subset S/{G_x}$, whereas the $G^0$-orbit type $\Sigma_0$ through the $G^0$-orbit containing $x$ becomes $\Gamma\times_{\Gamma_x}S^{G_x^0}\subset \Gamma\times_{\Gamma_x}(S/G_x^0)$, with $S^{G_x^0}\subset S/G_x^0$ the $G_x^0$-fixed point set, and the $\Gamma$-orbit type through $\O_x$ of the residual $\Gamma$-action on $\Sigma_0$ becomes the $\Gamma_x$-fixed point set $(S^{G_x^0})^{\Gamma_x}=S^{G_x}\subset S/{G_x}$. So, the germs at $\O_x$ of the members of the two partitions through $\O_x$ indeed coincide. 
\end{proof}
\begin{lemma}\label{lemma:orbittypes:affineaction:finitegp} Let $M$ be an integral affine manifold equipped with an action of a finite group $\Gamma$ by integral affine isomorphisms, and let $\underline{\Sigma}$ be a connected component of a $\Gamma$-orbit type in $\underline{M}:=M/\Gamma$. Then the corresponding invariant subset $\Sigma$ of $M$ is an integral affine submanifold and $\underline{\Sigma}$ admits a (necessarily unique) integral affine structure making the quotient map $\pi:\Sigma\to \underline{\Sigma}$ (which is a surjective local diffeomorphism) a local integral affine isomorphism.
\end{lemma}
\begin{proof} By the standard theory of proper Lie group actions, $\Sigma$ is a submanifold of $M$ and $\underline{\Sigma}$ admits a (necessarily unique) smooth structure for which $\pi:\Sigma\to \underline{\Sigma}$ is a submersion. Since $\Gamma$ is finite, $\pi:\Sigma\to \underline{\Sigma}$ in fact a local diffeomorphism and, given $x\in \Sigma$, there is a $\Gamma_x$-invariant connected open $U$ in $M$ around $x$ which is the domain of an integral affine chart for $M$ and is such that $(\gamma\cdot U)\cap U\neq \emptyset$ if and only if $\gamma\in \Gamma_x$. Note that, for such $x$ and $U$, there is a unique integral affine open embedding $\iota:U\to T_xM$ that maps $x$ to the origin and has derivative at $x$ equal to the identity map of $T_xM$. This can be constructed, for instance, by choosing an integral affine chart $\chi$ with domain $U$ and setting $\iota:=(\d \chi_x)^{-1}\circ \big(\chi-\chi(x)\big)$. The map $\iota$ is $\Gamma_x$-equivariant with respect to the induced linear $\Gamma_x$-action on $T_xM$, since for any $\gamma\in \Gamma_x$ the map $(\d m_{\gamma})^{-1}_x\circ \iota\circ m_{\gamma}:U\to T_xM$ is again an integral affine open embedding $U\to T_xM$ that maps $x$ to the origin and has derivative at $x$ equal to the identity map of $T_xM$, hence it must equal $\iota$. So, $\iota$ linearizes both the $\Gamma_x$-action on $U$ and the integral affine structure. Since $U$ is a slice at $x$ for the $\Gamma$-action it holds that $U\cap \Sigma=U^{\Gamma_x}$ and, hence, $\iota(U\cap \Sigma)=\iota(U)\cap T_xM^{\Gamma_x}$. This shows that $\Sigma$ is an integral affine submanifold of $M$, because $T_xM^{\Gamma_x}$ is a rational linear subspace of the $(T_xM,\Lambda_x^*)$ (as follows from the fact that the linear map $T_xM\to T_xM^{\Gamma_x}$, $v\mapsto \sum_{\gamma\in \Gamma_x}\gamma\cdot v$ is surjective and preserves $\Lambda_x^*$). Since $\pi:\Sigma\to \underline{\Sigma}$ is a quotient map for the induced $\Gamma$-action on $\Sigma$ and since this action is by integral affine automorphisms of $\Sigma$, the integral affine structure descends to $\underline{\Sigma}$.
\end{proof}

\begin{proof}[Proof of Proposition \ref{prop:affstr:coadjoint-orbit-types}] In view of Lemma \ref{lemma:orbit-types:dis-con-gps} there is a unique $G^0$-orbit type $\Sigma_0$ in $\g^*/G^0$ for which $\underline{\Sigma}$ is a connected component of a $\Gamma$-orbit type in $\Sigma_0/\Gamma$ of the residual $\Gamma$-action on $\Sigma_0$. For the coadjoint $G^0$-action, the $G^0$-orbit types in $\g^*/G^0$ have equidimensional connected components, because the dimension of its connected component through the orbit containing a point $x\in \g^*$ equals the dimension of $(\g_x^*)^{G^0_x}$. So, the subset $M$ of $\t^*_+$ corresponding to $\Sigma_0$ via the canonical homeomorphism between $\t^*_+$ and $\g^*/G^0$ is a disjoint union of equidimensional open faces of $\t^*_+$ and, hence, $M$ is an integral affine submanifold of $(\t^*,\Lambda_T^*)$. The group $\Gamma$ is isomorphic to $\Gamma_+:=N_G(\t^*_+)/T$, with $N_G(\t^*_+)$ the normalizer of $\t^*_+$ in $G$, via the group homomorphism induced by the inclusion of $N_G(\t^*_+)$ into $G$ (cf. \cite[Proposition 3.15.1]{DuKo}). Furthermore, the canonical homeomorphism between $\t^*_+$ and $\g^*/G^0$ intertwines the $\Gamma_+$-action on $\t^*_+$ induced by the coadjoint $G$-action with the residual $\Gamma$-action on $\g^*/G^0$. So, since $\Sigma_0$ is $\Gamma$-invariant, $M$ is $\Gamma_+$-invariant. Because of this and the fact that $\Gamma_+$ acts by linear integral affine automorphisms of $(\t^*,\Lambda^*_T)$, $M$ inherits a $\Gamma_+$-action by integral affine isomorphisms. Applying Lemma \ref{lemma:orbittypes:affineaction:finitegp} to $M$ and $\underline{\Sigma}$, the proposition follows. 
\end{proof}
The integral affine structures on the strata of $\sS_G(\underline{\g^*})$ have a more intrinsic description. For this recall that, as for any proper Lie group action, the tangent spaces to the strata $\underline{\Sigma}\in \sS_G(\underline{\g^*})$ are determined by the property that the differential of the quotient map $\Sigma\to \underline{\Sigma}$ induces an isomorphism
\begin{equation*} T_{\L_x}\underline{\Sigma}\cong\No_x^{G_x}, \quad x\in \Sigma,
\end{equation*} where $\No_x$ is the fiber of the normal bundle to the $G$-orbit $\L_x$, equipped with the linear $G_x$-action obtained by differentiating the action. For the coadjoint action, the tangent space $T_x\L$ is the annihilator $\g_x^0\subset \g^*$. Therefore, $\No_x$ is canonically isomorphic to $\g_x^*$. Since this isomorphism identifies the $G_x$-action on $\No_x$ with the coadjoint action on $\g^*_x$, we obtain an isomorphism
\begin{equation}\label{eqn:tgntsp:stratum:orbsp:coadjointact} T_{\L_x}\underline{\Sigma}\cong (\g_x^*)^{G_x}.
\end{equation} 
For any finite-dimensional real representation $V$ of a compact Lie group $G$, the fixed point set $V^G$ has a canonical $G$-invariant linear complement $\c_V$ given by the linear span of the set
\begin{equation*} \{g\cdot v- v\mid v\in V,\text{ }g\in G\}.
\end{equation*}
Alternatively, $\c_V$ is the orthogonal complement of $V^G$ with respect to any $G$-invariant inner product on $V$. It is also characterized by the property that the annihilator of $\c_V$ is the fixed point set $(V^*)^G$ of the dual representation. Using the latter and dualizing \eqref{eqn:tgntsp:stratum:orbsp:coadjointact}, we get an identification 
\begin{equation}\label{eqn:cotangentspace:orbittypestratum:coadjointact} T^*_{\L_x}\underline{\Sigma}\cong \g_x^{G_x}
\end{equation} with the Lie algebra of the center $Z(G_x)$ of $G_x$. Since the canonical complement $\c_{\g_x}$ to $\g_x^{G_x}$ is a Lie subalgebra of $\g_x$ (in fact, it is a Lie algebra ideal), we can consider the connected Lie subgroup $C_{\g_x}$ of $G_x$ that integrates it. With this we can explicitly describe the lattice subbundle $\Lambda_{\underline{\Sigma}}$ of $T^*\underline{\Sigma}$ that encodes the integral affine structure (in the sense of \cite[Proposition 3.1.2]{CrFeTo}). 
\begin{prop}\label{prop:intrinsicdescr:IAstr:orbittypestrata:coadjointact} Let $\underline{\Sigma}\in \sS_G(\underline{\g^*})$ and $x\in \Sigma$. Under the identification \eqref{eqn:cotangentspace:orbittypestratum:coadjointact}, the fiber of $\Lambda_{\underline{\Sigma}}$ at $\L_x$ becomes the lattice $\exp_{Z(G_x)}^{-1}(C_{\g_x})$ in $\g_x^{G_x}$.
\end{prop}
\begin{proof} We can assume that $x\in \Sigma_+$. The proofs of Proposition \ref{prop:affstr:coadjoint-orbit-types} and Lemma \ref{lemma:orbittypes:affineaction:finitegp} show that 
$T_x\Sigma_+=(T_xF)^{(\Gamma_+)_x}$, with $F$ the open fact of $\t^*_+$ through $x$. Using that the inclusion map from $N_{G_x}(\t^*_+)$ into $G_x$ induces an isomorphism between $(\Gamma_+)_x$ and $G_x/G_x^0$, and using that $T_xF=(\g_x^*)^{G^0_x}$, it follows that $T_x\Sigma_+=(\g_x^*)^{G_x}\subset \t^*$. Therefore, the fiber at $x$ of the subbundle $\Lambda_{\Sigma_+}$ of $T^*\Sigma_+$ that encodes the integral affine structure on $\Sigma_+$ is identified with the image of $\Lambda_T$ under the projection $\t=\g_x^{G_x}\oplus (\c_{\g_x}\cap \t)\to \g_x^{G_x}$. If $G_x$ were connected, it would hold that $Z(\g_x)=\g_x^{G_x}$ and $\c_{\g_x}=\g_x^\mathrm{ss}$, so the proof of Lemma \ref{lemma:rationality:Lietheoreticsubspaces} would show that this image is the lattice $\exp_{Z(G_x)}^{-1}(C_{\g_x})$. The arguments in that proof readily extend to show that the same holds when $G_x$ is not connected. \end{proof}
\subsubsection{The remaining proofs} We now turn to the remaining proofs of the statements in \S\ref{sec:extensionpf:disconnectedgps}. For the proof of Theorem \ref{thm:stratification:image:momentum:map:disconnectedgps}, denote by 
\begin{equation*} q^0:S\to S/G^0,\quad q_0:S/G^0\to S/G,
\end{equation*} the respective quotient maps for $G^0$-action on $S$ and the residual $\Gamma$-action on $S/G^0$. Note that, since $\sS_\mathrm{Ham}(\underline{S})$ refines the orbit type stratification, each stratum $\underline{\Sigma}\in \sS_\mathrm{Ham}(\underline{S})$ is a submanifold of a connected component $\underline{\Sigma}_G$ of a $G$-orbit type in $\underline{S}$. Moreover, in view of Lemma \ref{lemma:orbit-types:dis-con-gps}, $q_0^{-1}(\underline{\Sigma}_G)$ is a submanifold of a unique $G^0$-orbit type in $S/G^0$ and $q_0$ restricts to a submersion $q_0^{-1}(\underline{\Sigma}_G)\to \underline{\Sigma}_G$. Therefore, $q_0^{-1}(\underline{\Sigma})$ is a submanifold of the same $G^0$-orbit type. Let $\sS_\mathrm{Ham}(S/G^0)$ be the partition of $S/G^0$ consisting of the connected components of the members of $\{q_0^{-1}(\underline{\Sigma})\mid \underline{\Sigma}\in\sS_\mathrm{Ham}(\underline{S})\}$, equipped with the induced manifold structure. We will use the following fact. 
\begin{prop}\label{prop:canhamstratG^0=strat} $\sS_\mathrm{Ham}(S/G^0)$ is a stratification. 
\end{prop}
\begin{proof} Since the members of $\sS_\mathrm{Ham}(\underline{S})$ are locally closed in $\underline{S}$, those of $\sS_\mathrm{Ham}(S/G^0)$ are locally closed in $S/G^0$. For local finiteness note that, given $\underline{\Sigma}\in\sS_\mathrm{Ham}(\underline{S})$, the $\Gamma$-action on $q_0^{-1}(\underline{\Sigma})$ permutes connected components. It follows from this, the fact that $q_0$ is an open map, and connectedness of $\underline{\Sigma}$, that the image of each connected component of $q_0^{-1}(\underline{\Sigma})$ under $q_0$ is all of $\underline{\Sigma}$. Consequently, for each such $\underline{\Sigma}$ the action of $\Gamma$ on the set of connected components of $q_0^{-1}(\underline{\Sigma})$ is transitive and, hence,  $q_0^{-1}(\underline{\Sigma})$ has finitely many connected components. This and local finiteness of $\sS_\mathrm{Ham}(\underline{S})$ imply that  $\sS_\mathrm{Ham}(S/G^0)$ is locally finite as well. Finally, using the symplectic cross-section theorem (Theorem \ref{thm:symplectic:cross-section}) and the MGS normal form theorem, it follows from an argument along the same lines as that at the end of the proof of \cite[Theorem 2.1]{LeSj} that $\sS_\mathrm{Ham}(S/G^0)$ satisfies the frontier condition.  
\end{proof}
\begin{proof}[Proof of Theorem \ref{thm:stratification:image:momentum:map:disconnectedgps}] As in the case of connected $G$, we will first show that there is a unique stratification $\sS_\mathrm{Ham}(\underline{\Delta})$ of $\underline{\Delta}$ for which each stratum is an affine submanifold of a stratum of $\sS_G(\underline{\g^*})$ with tangent space given by \eqref{eqn:definingproperty:stratofimage:disconnectedgps}. 
Let $J_0:S/G^0\to \g^*/G^0$ be the map induced by $J$ and let $J_{0,+}:S/G^0\to \t^*_+$ be the composition of $J_0$ with the canonical homeomorphism between $\g^*/G^0$ and $\t^*_+$. Consider the cover $\A_\mathrm{Ham}$ of $\Delta_+=J_{0,+}(S/G^0)$ given by
\begin{equation*} \A_\mathrm{Ham}:=\{J_{0,+}(\Sigma_0)\mid \Sigma_0\in \sS_\mathrm{Ham}(S/G^0)\}.
\end{equation*} As in the proof of Proposition \ref{prop:images:canhamstrata:affineopen} one can show that each member $J_{0,+}(\Sigma_0)$ of $\A_\mathrm{Ham}$ is affine-open in $\t^*$, that $J_{0,+}$ restricts to a submersion from $\Sigma_0$ into the affine hull of $J_{0,+}(\Sigma_0)$, and that this affine hull given by
\begin{equation*} \mathrm{Af}(J_{0,+}(\Sigma_0))=x+(\g_p^0)^{G_x}\subset \t^*
\end{equation*} for any $x\in J_{0,+}(\Sigma_0)$ and $p\in (q^0)^{-1}(\Sigma_0)$ such that $J(p)=x$. Furthermore, it follows as in the proof of Theorem \ref{thm:stratification:image:momentum:map} (using Proposition \ref{prop:canhamstratG^0=strat}) that $\A_\mathrm{Ham}$ is a piecewise-affine cover of $\Delta_+$. So, appealing to Proposition \ref{prop:piecewise-affine:stratification}, we find that there is a (necessarily unique) stratification $\sS_\mathrm{Ham}(\Delta_+)$ of $\Delta_+$ into affine submanifolds of $\t^*$ with the property that 
\begin{equation*} T_x\widehat{\sigma}=\bigcap_{p\in J^{-1}(x)}(\g_p^0)^{G_x}\subset \t^*
\end{equation*} for all $\widehat{\sigma}\in \sS_\mathrm{Ham}(\Delta_+)$ and $x\in \widehat{\sigma}$. 
Moreover, by Lemma \ref{lemma:affinestrat:refinementcrittangentspaces}, each stratum of $\sS_\mathrm{Ham}(\Delta_+)$ is contained in a connected component of the subset $\Sigma_+$ of $\t^*_+$ corresponding to some $\underline{\Sigma}\in \sS_G(\underline{\g^*})$. Since the finite group $\Gamma_+$ (as in the proof of Proposition \ref{prop:affstr:coadjoint-orbit-types}) acts by linear automorphisms of $\t^*$ and permutes the members of $\A_\mathrm{Ham}$, it also permutes the strata of $\sS_\mathrm{Ham}(\Delta_+)$. Therefore, the saturation $\sigma_+:=\Gamma_+\cdot \widehat{\sigma}$ is a finite disjoint union of such strata and, hence, is an affine submanifold of $\Sigma_+$ as well. Quotienting each such $\sigma_+$ by the $\Gamma_+$-action, we obtain a partition $\sS_\mathrm{Ham}(\underline{\Delta})$ of $\underline{\Delta}$ with the property that each member is a connected affine submanifold of a stratum of $\sS_G(\underline{\g^*})$, with tangent spaces given by \eqref{eqn:definingproperty:stratofimage:disconnectedgps}. To see that $\sS_\mathrm{Ham}(\underline{\Delta})$ is a stratification, we ought to verify that its members are locally closed, that it is locally finite and that it satisfies the frontier condition. Each of its members is locally closed in $\underline{\Delta}$, being a submanifold of a stratum of $\sS_G(\underline{\g^*})$. Moreover, local finiteness of $\sS_\mathrm{Ham}(\underline{\Delta})$ readily follows from that of $\sS_\mathrm{Ham}(\Delta_+)$. Finally, to see that $\sS_\mathrm{Ham}(\underline{\Delta})$ satisfies the frontier condition suppose that $\underline{\sigma}_1,\underline{\sigma}_2\in \sS_\mathrm{Ham}(\underline{\Delta})$ are such that $\underline{\sigma}_1$ intersects the closure $\mathrm{Cl}_{\underline{\Delta}}(\underline{\sigma}_2)$ of $\underline{\sigma}_2$ in $\underline{\Delta}$. By construction of $\sS_\mathrm{Ham}(\underline{\Delta})$, there are $\widehat{\sigma}_1,\widehat{\sigma}_2\in \sS_\mathrm{Ham}(\Delta_+)$ such that $\pi_+(\widehat{\sigma}_1)=\underline{\sigma}_1$ and $\pi_+(\widehat{\sigma}_2)=\underline{\sigma}_2$. Since $\pi_+:\Delta_+\to \underline{\Delta}$ is continuous and closed, it holds that 
\begin{equation*} 
\pi_+(\mathrm{Cl}_{\Delta_+}(\widehat{\sigma}_2))=\mathrm{Cl}_{\underline{\Delta}}(\pi_+(\widehat{\sigma}_2))=\mathrm{Cl}_{\underline{\Delta}}(\underline{\sigma}_2).
\end{equation*}
 In view of this and the fact that $\Gamma_+$ permutes the strata of $\sS_\mathrm{Ham}(\Delta_+)$, there is in fact a $\widehat{\sigma}_1$ as above such that $\widehat{\sigma}_1$ intersects $\mathrm{Cl}_{\Delta_+}(\widehat{\sigma}_2)$. Since $\sS_\mathrm{Ham}(\Delta_+)$ satisfies the frontier condition, it follows that $\widehat{\sigma}_1\subset \mathrm{Cl}_{\Delta_+}(\widehat{\sigma}_2)$ and $\dim(\widehat{\sigma}_1)<\dim(\widehat{\sigma}_2)$ and, hence, $\underline{\sigma}_1\subset \mathrm{Cl}_{\underline{\Delta}}(\underline{\sigma}_2)$ and $\dim(\underline{\sigma}_1)<\dim(\underline{\sigma}_2)$. This shows that $\sS_\mathrm{Ham}(\underline{\Delta})$ also satisfies the frontier condition. So, it is a stratification of $\underline{\Delta}$ for which each stratum is an affine submanifold of a stratum of $\sS_G(\underline{\g^*})$ with tangent spaces given by \eqref{eqn:definingproperty:stratofimage:disconnectedgps}. Such a stratification is unique for the same reason that the stratification in Lemma \ref{lemma:localtoglobal:abstractaffstrat} is unique.\\
 
To prove that the strata of this stratification are integral affine, we ought to show that $\bigcap_{p\in J^{-1}(x)} (\g_p^0)^{G_x}$ is a rational linear subspace of $(\t^*,\Lambda_T^*)$ for each $x\in \Delta_+$. Since any intersection of rational linear subspaces is again rational, this boils down to showing that Lemma \ref{lemma:rationality:Lietheoreticsubspaces} in fact holds for arbitrary (possibly disconnected) compact Lie groups $G$, which can be done using the same arguments as in the connected case, with the role of $Z(\g)$ taken by the Lie algebra of $Z(G)$, which still equals $\g^G$, and the role of $\g^\mathrm{ss}$ now taken by the linear complement $\c_{\g}$ (see \S\ref{sec:iastr:orbittypes:coadjointact}), the annihilator of which still equals $(\g^*)^G$. \\

To complete the proof of the theorem it remains to extend the arguments for the proof of Theorem \ref{thm:localtriviality} to prove the equivariant local triviality in Theorem \ref{thm:stratification:image:momentum:map:disconnectedgps}. It suffices to extend the statement of Lemma \ref{lemma:liftingvecfields:stratumwise} for a vector field $X$ defined on the domain of an affine chart for $\underline{\sigma}$ in which it is a constant vector field. To prove the existence of lifts $(\widehat{X},\widehat{X}^\sigma)$ as in that lemma for a given such $X$, notice that, for $\underline{\sigma} \in \sS_\mathrm{Ham}(\underline{\Delta})$, $x\in \sigma_+$ and $\mathfrak{S}$ a slice through $x$ that is open in $x+\g^*_x$ and such that $\mathfrak{S}\cap \sigma$ is open in $\sigma_+$, the maps $q:\sigma\to \underline{\sigma}$ and $\pi_+:\sigma_+\to \underline{\sigma}$ restrict to the same map $\mathfrak{S}\cap\sigma \to \underline{\sigma}$, which is an affine open embedding into the domain of $X$ when $\mathfrak{S}$ is chosen small enough. So, we can extend the definition of the local lift $\overline{X}^\sigma$ by taking it to be the pull-back of $X$ along this open embedding. Then $\overline{X}^\sigma$ is a constant vector field on the affine-open $\mathfrak{S}\cap \sigma$ in $\t^*$ and the rest of the proof of Lemma \ref{lemma:liftingvecfields:stratumwise} goes through. It is readily verified that with this extension of Lemma \ref{lemma:liftingvecfields:stratumwise} the other arguments in the proof of Theorem \ref{thm:localtriviality} go through as well, by choosing $\varepsilon>0$ such that 
\begin{equation*} [-\varepsilon,\varepsilon]^n\to \underline{\sigma},\quad t\mapsto (\Phi_{X_1}^{t_1}\circ ...\circ \Phi_{X_n}^{t_n})(\O),
\end{equation*} is defined (and, hence, restricts to an affine diffeomorphism from $]-\varepsilon,\varepsilon[^n$ onto an open $\underline{W}$ in $\underline{\sigma}$), and by letting the roles of $q_+$, $W_+$ and $\tau$ be played by the quotient map $q:\g^*\to \underline{\g^*}$, $\underline{W}$ and the inverse of this affine diffeomorphism.
\end{proof}

Next, we turn to the proofs of the remaining two statements.
\begin{proof}[Proof of the extensions of Theorem \ref{thm:linvar:main} and Corollary \ref{cor:sympvolpolynomial}] The proof of Theorem \ref{thm:linvar:main} extends as follows. First, Proposition \ref{prop:redsympforms:smoothsection} readily goes through and Lemma \ref{lemma:altcharlinvar} extends by replacing constant vector fields on $\sigma_+$ by vector fields on $\underline{\sigma}$ that are parallel with respect to the flat torsion-free connection defining the affine structure on $\underline{\sigma}$, which is equivalent to the vector fields being constant in any affine coordinate chart with connected domain. Given such a vector field $X$ and an $\L\in \underline{\sigma}$, we can complete $X$ to a commuting frame $X_1,...,X_m:=X$ of $\underline{\sigma}$ defined on the domain of an affine chart for $\underline{\sigma}$ around $\L$, such that each of the vector fields $X_l$ is constant in this affine chart. Following the proof of equivariant local triviality in Theorem \ref{thm:stratification:image:momentum:map:disconnectedgps} outlined at the end of the proof of that theorem, we obtain pairs of lifts $(\widehat{X}_1,\widehat{X}^\sigma_1),...,(\widehat{X}_m,\widehat{X}^\sigma_m)$ with the property that $\widehat{X}^\sigma_l$ coincides with $(\pi_+\vert_{\sigma_+})^*X_l$ on $\sigma_+$ for each $l$, as well as an associated local trivialization $(\Phi,\varphi)$ over $\underline{W}$ and a smooth equivariant extension of $\Phi^{-1}$. Then, for each $x\in \L\cap \t^*_+$, $\varphi^{-1}(\{x\}\times \underline{W})$ is contained (and open) in $W_+$, $\pi_+$ restricts to an affine diffeomorphism between $\varphi^{-1}(\{x\}\times \underline{W})$ and $\underline{W}$, and $\widehat{X}^\sigma:=\widehat{X}^\sigma_m$ restricts to a constant vector field on $\varphi^{-1}(\{x\}\times \underline{W})$. Therefore, by fixing such an $x$ and parametrizing $\underline{W}$ by $y\in \varphi^{-1}(\{x\}\times \underline{W})$, the rest of the proof of Theorem \ref{thm:linvar:main} goes through. Finally, the proof of Corollary \ref{cor:sympvolpolynomial} readily extends by working locally in affine coordinates for $\underline{\sigma}$. 
\end{proof}

\subsection{Extension to actions of proper quasi-symplectic groupoids}\label{sec:extension:quasisympgpoids}
\subsubsection{The canonical stratification of a proper quasi-symplectic groupoid and its integral affine structure}\label{sec:canstrat:quasisympgpoids} In the remainder of the paper we extend the theory developed thus far to Hamiltonian actions of proper quasi-symplectic groupoids (as in \cite{BuCrWeZh,Xu1}). As shown in \cite{CrFeTo2}, the canonical stratification of the orbit space of such a groupoid (as in \cite{PfPoTa,CrMe}) comes with a natural integral affine structure on each stratum. This stratification generalizes the orbit-type stratification $\sS_G(\underline{\g^*})$ of the coadjoint orbit space $\underline{\g^*}$ and, in the extended theory, it will take the role played by $\sS_G(\underline{\g^*})$ in the versions of our main theorems for compact Lie group actions. Before stating the extended theorems, in this subsection we recall the definition of the canonical stratification and the integral affine structures on its strata. \\

Let $(\G,\Omega,\phi)$ be a proper quasi-symplectic groupoid over $M$, where $\Omega\in \Omega^2(\G)$ denotes the multiplicative quasi-symplectic $2$-form on $\G$ and $\phi \in \Omega^3(M)$ a $3$-form such that $\d\Omega=s^*\phi-t^*\phi$ (for basics on quasi-symplectic groupoids we refer the reader to \cite{BuCrWeZh,Xu1}). As for any proper Lie groupoid, its orbit space $\underline{M}:=M/\G$ comes with a natural stratification $\sS_\G(\underline{M})$, which consists of the connected components of the members of the partition $\P$ of $\underline{M}$ defined by the equivalence relation that declares the orbits through two points $x,y\in M$ to be equivalent when the isotropy groups $\G_x$ and $\G_y$ of $\G$ at these points are isomorphic (see \cite{PfPoTa,CrMe}). The usual proof of the fact that the members of $\sS_\G(\underline{M})$ are smooth further shows that 
\begin{itemize} \item for any $\underline{\Sigma}\in \sS_\G(\underline{M})$ the corresponding invariant subspace $\Sigma\subset M$ is a submanifold with tangent space determined by the fact that, for any $\L\in \underline{M}$ and $x\in \L$,
\begin{equation*} \frac{T_x\Sigma}{T_x\L}=\No_x^{\G_x}\subset \No_x:=\frac{T_xM}{T_x\L},
\end{equation*} where $\No_x^{\G_x}$ is the fixed-point set of the so-called normal representation $\No_x$ of $\G_x$,
\item the orbit space projection $q:M\to \underline{M}$ restricts to a smooth submersion $\Sigma \to \underline{\Sigma}$, the differential of which induces an isomorphism
\begin{equation}\label{eqn:tgntsp:stratum:orbsp:Liegpoid} \No_x^{\G_x}\cong T_{\L_x}\underline{\Sigma}.  
\end{equation}
\end{itemize}
The theory developed in \cite{CrFeTo,CrFeTo2} shows that the orbit space of a proper quasi-symplectic groupoid has a much richer geometric structure than that of a general proper Lie groupoid. Here we only need part of this richer structure: the natural integral affine structures on the strata of $\sS_\G(\underline{M})$. Given such a stratum $\underline{\Sigma}$, the lattice subbundle $\Lambda_{\underline{\Sigma}}$ of $T^*\underline{\Sigma}$ encoding its integral affine structure can be described explicitly, as follows. Recall from \cite{BuCrWeZh} that the base $M$ of the quasi-symplectic groupoid $(\G,\Omega,\phi)$ has an induced $\phi$-twisted Dirac structure $L\subset TM\oplus T^*M$, uniquely determined by the property that the target map $t:(\G,\Omega)\to (M,L)$ is forward Dirac. The Lie algebroid of $\G$ is isomorphic to the Lie algebroid $(L,[-,-]_\phi,\mathrm{pr}_{TM})$ of the Dirac structure $L$, via the isomorphism
\begin{equation*} \rho_\Omega:L\to \mathrm{Lie}(\G)
\end{equation*} that sends $(v,\alpha)\in L_x$ to the element $\rho_\Omega(v,\alpha)\in \mathrm{Lie}(\G)_x$ determined by the equations
\begin{equation*} \iota_{\rho_\Omega(v,a)}\Omega=\d t^*_{1_x}\alpha,\quad \d t(\rho_\Omega(v,a))=v. 
\end{equation*} This induces an isomorphism of isotropy Lie algebras
\begin{equation}\label{eqn:isotropyLiealg:Dirac:gpoid:iso} (\rho_\Omega)_x:\g(L)_x\to \g_x,\quad x\in M.
\end{equation} The vector space underlying the isotropy Lie algebra $\g(L)_x$ of $L$ at $x$ is the annihilator of $T_x\L$ (the tangent space to the orbit $\L$ of $\G$ through $x$, which equals the tangent space to the pre-symplectic leaf of $L$ through $x$) in $T_x^*M$ and so the dual of \eqref{eqn:isotropyLiealg:Dirac:gpoid:iso} induces a linear isomorphism
\begin{equation}\label{eqn:isotropyLiealg:Dirac:gpoid:iso:dual} (\rho_\Omega)_x^*:\g_x^*\to \No_x.
\end{equation} In fact, this is an isomorphism of representations between the coadjoint $\G_x$-representation and the aforementioned normal representation of $\G_x$. So, via \eqref{eqn:tgntsp:stratum:orbsp:Liegpoid} and \eqref{eqn:isotropyLiealg:Dirac:gpoid:iso:dual} we obtain an isomorphism
\begin{equation}\label{eqn:tgntsp:stratum:orbsp:quasisympgpoid} T_{\L_x}\underline{\Sigma}\cong (\g_x^*)^{\G_x}. 
\end{equation} Using that $(\g_x^*)^{\G_x}$ is the annihilator of $\c_{\g_x}$ (see \S\ref{sec:iastr:orbittypes:coadjointact}) and dualizing \eqref{eqn:tgntsp:stratum:orbsp:quasisympgpoid}, we obtain an isomorphism
\begin{equation*} T_{\L_x}^*\underline{\Sigma}\cong \g_x^{\G_x}. 
\end{equation*}
With this we can explicitly describe the lattice bundle $\Lambda_{\underline{\Sigma}}$: its fiber at $\L_x$ is the lattice in $T^*_{\L_x}\underline{\Sigma}$ corresponding to the lattice
\begin{equation}\label{eqn:lattice:IAstr:stratumquasisympgpoid} \exp_{Z(\G_x)}^{-1}(C_{\g_x})\subset \g_x^{\G_x},
\end{equation} where $C_{\g_x}$ is the connected Lie subgroup of $\G_x$ integrating $\c_{\g_x}$. To see that this lattice subbundle is indeed smooth and Lagrangian, one can verify that it coincides with the lattice subbundle of the integral affine structure induced (as in \cite[Section 3.3]{CrFeTo}) by a regular proper quasi-symplectic groupoid over $\Sigma$ with smooth orbit space equal to $\underline{\Sigma}$, which is constructed out of $(\G,\Omega,\phi)$ via reduction along $\Sigma$ as in the proof of \cite[Theorem 3.29]{CrFeTo2} (see also \cite[Proposition 2.106]{Mol1} for a more direct exposition of this reduction procedure). 
\begin{ex} Let $G$ be a compact Lie group and consider the cotangent symplectic groupoid $(G\ltimes\g^*,\Omega_\textrm{can})\rightrightarrows\g^*$. Since its underlying Lie groupoid is the action groupoid of the coadjoint $G$-action, it is proper, its orbit space is $\underline{\g^*}:=\g^*/G$ and the corresponding canonical stratification coincides with $\sS_G(\underline{\g^*})$. In view of Proposition \ref{prop:intrinsicdescr:IAstr:orbittypestrata:coadjointact}, the induced integral affine structures on the strata of $\sS_G(\underline{\g^*})$ coincide with those in Proposition \ref{prop:affstr:coadjoint-orbit-types}.
\end{ex}
For further details on the above integral affine structures and their relation to the geometry of the quasi-symplectic groupoid we refer to \cite{CrFeTo2}. 

\subsubsection{Extension of the statements} Turning to the extension of our theorems, let $(\G,\Omega,\phi)$ be a proper quasi-symplectic groupoid over $M$ and $J:(S,\omega)\to M$ a quasi-Hamiltonian $(\G,\Omega,\phi)$-space (as in \cite{BuCrWeZh,Xu1}). We let $\underline{S}:=S/\G$ and $\underline{M}:=M/\G$ denote the orbit spaces of the $\G$-action and of $\G$ itself, and we denote by
\begin{equation*} \underline{J}:\underline{S}\to \underline{M}
\end{equation*} the map induced by $J$. We call this the transverse momentum map and denote its image by $\underline{\Delta}$.
\begin{thm}\label{thm:stratification:image:momentum:map:groupoids} The subspace $\underline{\Delta}$ of $\underline{M}$ admits a natural 
integral affine stratification ${\sS}_\mathrm{Ham}(\underline{\Delta})$ with the property that, for each $\underline{\sigma}\in {\sS}_\mathrm{Ham}(\underline{\Delta})$, the restriction
\begin{equation}\label{eqn:restrictionstratum:mommap:gpoids} 
J:J^{-1}(\sigma)\to \sigma
\end{equation} is a $\G$-equivariantly locally trivial fibration. This integral affine stratification is uniquely 
determined by the property that each stratum is an integral affine submanifold of a stratum of $\sS_\G(\underline{M})$ and that for any $x\in \Delta$ the tangent space at $\L_x$ to the stratum $\underline{\sigma}$ through that orbit is given by
\begin{equation}\label{eqn:definingproperty:stratofimage:gpoids}
T_{\L_x}\underline{\sigma}=\bigcap_{p\in J^{-1}(x)} (\g_p^0)^{\G_x}\subset T_{\L_x}\underline{\Sigma}, 
\end{equation}  where 
\begin{itemize} 
\item $\g_p$ is the isotropy Lie algebra at $p\in S$ of the $\G$-action,
\item $\G_x$ is the isotropy group at $x$ of the Lie groupoid $\G$ and $\g_x$ is its Lie algebra,
\item $\g_p^0\subset \g_x^*$ is the annihilator of $\g_p$ in $\g_x$, 
\item $(\g_p^0)^{G_x}:=\g_p^0\cap (\g_x^*)^{G_x}$ is viewed as linear subspace of the tangent space at $x$ of the stratum $\underline{\Sigma}\in \sS_\G(\underline{M})$ containing $\underline{\sigma}$ via the identification \eqref{eqn:tgntsp:stratum:orbsp:quasisympgpoid}. 
\end{itemize}  
\end{thm}
By equivariant local triviality we now mean that, given $\underline{\sigma}\in \sS_\mathrm{Ham}(\underline{\Delta})$ and $\L \in \underline{\sigma}$, there is an open $W$ in $\sigma$ around $\L$ (that is not necessarily $\G$-invariant), an open $O$ in $\L\times \underline{W}$ around $\L\times \{0\}$ (with $\underline{W}$  the image of $W$ under the quotient map $\sigma\to \underline{\sigma}$), an isomorphism of Lie groupoids
\begin{equation*}
\widehat{\varphi}:\G\vert_W\xrightarrow{\sim} (\G_{\L}\times \underline{W})\vert_{O}
\end{equation*}




that restricts to the identity on $\G_\L:=\G\vert_\L$ (where $\G_{\L}\times \underline{W}$ denotes the product of $\G_\L$ with the unit groupoid of $\underline{W}$) and a $\widehat{\varphi}$-equivariant homeomorphism $\Phi$ fitting into a commutative diagram
\begin{center} 
\begin{tikzcd} 
J^{-1}(W)\arrow[rr,"\Phi", "\sim"'] \arrow[d,"J"'] & & (J\times\textrm{id}_{\underline{W}})^{-1}(O) \arrow[hookrightarrow,r] & J^{-1}(\L)\times \underline{W} \arrow[d,"J\times\mathrm{id}_{\underline{W}}"] \\
W \arrow[rd,"\pi"']\arrow[rr,"\varphi", "\sim"'] & & O \arrow[hookrightarrow,r]& \L\times \underline{W}\arrow[lld,"\mathrm{pr}_2"] \\
  & \underline{W} & &
\end{tikzcd} 
\end{center} with $\varphi$ the map covered by $\widehat{\varphi}$. Moreover, these can be chosen so that $\Phi$ and $\Phi^{-1}$ extend to smooth maps from an open in $S$ into $S\times \underline{W}$ and from an open in $S\times \underline{W}$ into $S$, respectively. 
\begin{rem}\label{rem:choiceofgpoidiso} In the proof of Theorem \ref{thm:stratification:image:momentum:map:groupoids} the isomorphism $\widehat{\varphi}$ will be obtained from an isomorphism $\widetilde{\varphi}$ between the restriction of the quasi-symplectic groupoid to a neighbourhood $V$ of $\L$ and its linear local model \cite[Section 2]{CrFeTo2}, by restricting $\widetilde{\varphi}$ to an open $W$ in $V\cap \sigma$ around $\L$. The proof will show that, for any $\widehat{\varphi}$ obtained by restricting such a linearization $\widetilde{\varphi}$ to a small enough such $W$, there is a $\Phi$ as above. 
\end{rem}


\begin{ex} Various proper quasi-symplectic groupoids of specific interest, such as the cotangent groupoid or the AMM groupoid of a compact Lie group (see, e.g., \cite{Xu1}), are the action groupoid of a Lie group action and, around any orbit $\L$, can be linearized via an isomorphism of quasi-symplectic groupoids that is induced by an equivariant tubular neighbourhood embedding around $\L$. For instance, for the cotangent groupoid $(G\ltimes \g^*,\Omega_\textrm{can})$ of a compact Lie group $G$ and a coadjoint orbit $\L$ through $x\in \g^*$, consider the $G$-equivariant tubular neighbourhood embedding 
\begin{equation*}
G\times_{G_x}\mathfrak{S} \to \g^*
\end{equation*} associated to a slice $\mathfrak{S}$ at $x$ that is open in $\g^*_x$. This  lifts to an isomorphism of Lie groupoids
\begin{equation*} \widetilde{\varphi}:G\ltimes (G\times_{G_x}\mathfrak{S}) \xrightarrow{\sim} (G\ltimes \g^*)\vert_{G\cdot \mathfrak{S}}
\end{equation*} that identifies the symplectic groupoid $(T^*G,\omega_\textrm{can})$ with its linear local model around $\L$ (see, e.g., \cite[\S 1.3.4]{Mol1} for more details). From the existence of $\Phi$ as in Remark \ref{rem:choiceofgpoidiso} for this specific $\widetilde{\varphi}$, we recover Theorem \ref{thm:stratification:image:momentum:map:disconnectedgps}. Using a similar choice of $\widetilde{\varphi}$ for the AMM groupoid, it follows that the analogue of Theorem \ref{thm:stratification:image:momentum:map:disconnectedgps} for the group-valued momentum maps of \cite{AlMaMe} holds as well.
\end{ex} 

All of the other results and corollaries stated in the introduction (as well as Remark \ref{rem:trivstratpreserving}) continue to hold for quasi-Hamiltonian $(\G,\Omega,\phi)$-spaces after making the same substitutions as for their extension to possibly disconnected compact Lie groups, but with the role of $\underline{\g^*}$ now taken by $\underline{M}$ and the role of $EG$ in Corollary \ref{cor:constructibility:pushforwards} taken over by the total space $E\G$ of the universal $\G$-bundle $E\G\to B\G$ (as in \cite{Haef}). More details on this will be given in the coming sections. To be more precise, in \S\ref{sec:pf:constructibility:quasisymgpoids} we prove Theorem \ref{thm:stratification:image:momentum:map:groupoids}. The extensions of Corollary \ref{cor:thm:stratification:image:transverse:momentum:map} and Corollary \ref{cor:constructibility:pushforwards}(b) are immediate from that. Those of Corollary \ref{cor:constructibility:pushforwards}(a) and Proposition \ref{prop:specializationmaps} are proved in \S\ref{sec:pushforwardsheaves:gpoidactions} and, finally, in \S\ref{sec:Sjamaar:deRham:gpoidactions} we elaborate on Theorem \ref{thm:linvar:main} and Corollary \ref{cor:sympvolpolynomial}.  
\subsubsection{Morita equivalence of Hamiltonian actions and proof of constructibility}\label{sec:pf:constructibility:quasisymgpoids} We now turn to the proof of Theorem \ref{thm:stratification:image:momentum:map:groupoids}, the idea of which is to reduce to the case of Hamiltonian actions of compact Lie groups using the linearization theorem for proper quasi-symplectic groupoids and basic properties of Morita equivalence between Hamiltonian actions. This variant of Morita equivalence was introduced in \cite{Mol1} for Hamiltonian actions of symplectic groupoids. The definition given there extends to the quasi-symplectic setting simply by using quasi-symplectic Morita equivalences instead of just symplectic ones. More explicitly, this means that given two quasi-symplectic groupoids $(\G_1,\Omega_1,\phi_1)$ and $(\G_2,\Omega_2,\phi_2)$, a Morita equivalence between a Hamiltonian $(\G_1,\Omega_1,\phi_1)$-space $J_1:(S_1,\omega_1)\to M_1$ and a Hamiltonian $(\G_2,\Omega_2,\phi_2)$-space $J_2:(S_2,\omega_2)\to M_2$ consists of
\begin{itemize}
\item a quasi-symplectic Morita equivalence $(P,\omega_P,\alpha_1,\alpha_2)$ between $(\G_1,\Omega_1,\phi_1)$ and $(\G_2,\Omega_2,\phi_2)$,
\item a Morita equivalence $(Q,\beta_1,\beta_2)$ between $\G_1\ltimes S_1$ and $\G_2\ltimes S_2$,
\item a smooth map $j:Q\to P$ making
\begin{center}
\begin{tikzpicture} 
\node (H1) at (-0.8,0) {$\G_1\ltimes S_1$};
\node (S1) at (-0.8,-1.3) {$S_1$};
\node (Q) at (1.35,0) {$Q$};
\node (S2) at (3.5,-1.3) {$S_2$};
\node (H2) at (3.5,0) {$\G_2\ltimes S_2$};

\node (G1) at (0,-3) {$\G_1$};
\node (M1) at (0,-4.3) {$M_1$};
\node (P) at (1.35,-3) {$P$};
\node (M2) at (2.7,-4.3) {$M_2$};
\node (G2) at (2.7,-3) {$\G_2$};

\draw[->, bend right=50](H1) to node[pos=0.45,below] {$\mathrm{pr}_{\G_1}\text{ }\text{ }\text{ }\text{ }\text{ }$} (G1);
\draw[->, bend right=20](S1) to node[pos=0.45,below] {$J_1\text{ }\text{ }\text{ }\text{ }$} (M1);
\draw[->, bend left=50](H2) to node[pos=0.45,below] {$\text{ }\text{ }\text{ }\text{ }\text{}\text{ }\text{ }\text{ }\mathrm{pr}_{\G_2}$} (G2);
\draw[->, bend left=20](S2) to node[pos=0.45,below] {$\text{ }\text{ }\text{ }\text{ }J_2$} (M2);
\draw[->](Q) to node[pos=0.45,below] {$j\text{ }\text{ }\text{ }\text{ }$} (P);
 
\draw[->,transform canvas={xshift=-\shift}](H1) to node[midway,left] {}(S1);
\draw[->,transform canvas={xshift=\shift}](H1) to node[midway,right] {}(S1);
\draw[->,transform canvas={xshift=-\shift}](H2) to node[midway,left] {}(S2);
\draw[->,transform canvas={xshift=\shift}](H2) to node[midway,right] {}(S2);
\draw[->](Q) to node[pos=0.25, below] {$\text{ }\text{ }\beta_1$} (S1);
\draw[->] (0.3,-0.15) arc (315:30:0.25cm);
\draw[<-] (2.3,0.15) arc (145:-145:0.25cm);
\draw[->](Q) to node[pos=0.25, below] {$\beta_2$\text{ }} (S2); 
 
\draw[->,transform canvas={xshift=-\shift}](G1) to node[midway,left] {}(M1);
\draw[->,transform canvas={xshift=\shift}](G1) to node[midway,right] {}(M1);
\draw[->,transform canvas={xshift=-\shift}](G2) to node[midway,left] {}(M2);
\draw[->,transform canvas={xshift=\shift}](G2) to node[midway,right] {}(M2);
\draw[->](P) to node[pos=0.25, below] {$\text{ }\text{ }\alpha_1$} (M1);
\draw[->] (0.8,-3.15) arc (315:30:0.25cm);
\draw[<-] (1.9,-2.85) arc (145:-145:0.25cm);
\draw[->](P) to node[pos=0.25, below] {$\alpha_2$\text{ }} (M2);
\end{tikzpicture} 
\end{center} 
a map of bibundles and such that 
\begin{equation*} j^*\omega_P=\beta_1^*\omega_1-\beta_2^*\omega_2.
\end{equation*}
\end{itemize}  Recall here that a quasi-symplectic Morita equivalence
\begin{center}
\begin{tikzpicture} \node (G1) at (-0.7,0) {$(\G_1,\Omega_1,\phi_1)$};
\node (M1) at (-0.7,-1.3) {$M_1$};
\node (S) at (1.4,0) {$(P,\omega_P)$};
\node (M2) at (3.5,-1.3) {$M_2$};
\node (G2) at (3.5,0) {$(\G_2,\Omega_2,\phi_2)$};
 
\draw[->,transform canvas={xshift=-\shift}](G1) to node[midway,left] {}(M1);
\draw[->,transform canvas={xshift=\shift}](G1) to node[midway,right] {}(M1);
\draw[->,transform canvas={xshift=-\shift}](G2) to node[midway,left] {}(M2);
\draw[->,transform canvas={xshift=\shift}](G2) to node[midway,right] {}(M2);
\draw[->](S) to node[pos=0.25, below] {$\text{ }\text{ }\alpha_1$} (M1);
\draw[->] (0.7,-0.15) arc (315:30:0.25cm);
\draw[<-] (2.1,0.15) arc (145:-145:0.25cm);
\draw[->](S) to node[pos=0.25, below] {$\alpha_2$\text{ }} (M2);
\end{tikzpicture}
\end{center}
is a Morita equivalence between the Lie groupoids underlying two quasi-symplectic groupoids, together with a two-form $\omega_P\in \Omega^2(P)$, such that
\begin{equation*} m_{\G_1}^*\omega_P=\mathrm{pr}_P^*\omega_P+\mathrm{pr}_{\G_1}^*\Omega_{\G_1},\quad m_{\G_2}^*\omega_P=\mathrm{pr}_P^*\omega_P+\mathrm{pr}_{\G_2}^*\Omega_{\G_1},
\end{equation*} and
\begin{equation*} \alpha_2^*\phi-\alpha_1^*\phi=\d\omega_P.
\end{equation*} Here $m_{\G_1}:\G_1\times_{M_1}P\to P$ and $m_{\G_2}:P\times_{M_2}\G_2\to P$ denote the action maps and the other maps are the projections from $\G_1\times_{M_1}P$ and $P\times_{M_2}\G_2$, respectively. 
\begin{ex}\label{ex:Moreq:associatedHamsp} By \cite{Xu1}, given a quasi-symplectic Morita equivalence as above, any Hamiltonian $(\G_1,\Omega_1,\phi_1)$-space $J:(S,\omega_S)\to M_1$ has a corresponding Hamiltonian $(\G_2,\Omega_2,\phi_2)$-space 
\begin{equation*} P_*(J):(P\ast_{\G_1}S,\omega_{PS})\to M_1,
\end{equation*} with
\begin{equation*} P\ast_{\G_1}S:=(P\times_{M_1}S)/\G_1
\end{equation*} the quotient manifold obtained via the diagonal action, $\omega_{PS}$ the $2$-form that $\omega_P\oplus \omega_S\in \Omega^2(P\times_{M_1}S)$ descends to, $P_*(J)$ the map that sends $[p_P,p_S]\mapsto \alpha_2(p_P)$ and with $\G_2$-action along $P_*(J)$ given by 
\begin{equation*} g\cdot[p_P,p_S]=[p_P\cdot g^{-1},p_S].
\end{equation*} The given quasi-symplectic Morita equivalence canonically extends to a Morita equivalence 
\begin{center}
\begin{tikzpicture} 
\node (H1) at (-1,0) {$\G_1\ltimes S$};
\node (S1) at (-1,-1.3) {$S$};
\node (Q) at (1.35,0) {$P\times_{M_1} S$};
\node (S2) at (3.9,-1.3) {$P\ast_{\G_1}S$};
\node (H2) at (3.9,0) {$\G_2\ltimes \left(P\ast_{\G_1}S\right)$};

\node (G1) at (0,-3) {$\G_1$};
\node (M1) at (0,-4.3) {$M_1$};
\node (P) at (1.35,-3) {$P$};
\node (M2) at (2.7,-4.3) {$M_2$};
\node (G2) at (2.7,-3) {$\G_2$};

\draw[->, bend right=50](H1) to node[pos=0.45,below] {$\textrm{pr}_{\G_1}\quad\quad$} (G1);
\draw[->, bend right=20](S1) to node[pos=0.45,below] {$J\text{ }\text{ }\text{ }\text{ }$} (M1);
\draw[->, bend left=50](H2) to node[pos=0.45,below] {$\quad\quad\quad\textrm{pr}_{\G_2}$} (G2);
\draw[->, bend left=20](S2) to node[pos=0.45,below] {$\quad\quad\text{ }\text{ }P_*(J)$} (M2);
\draw[->](Q) to node[pos=0.45,below] {$\textrm{pr}_P\quad\quad$} (P);
 
\draw[->,transform canvas={xshift=-\shift}](H1) to node[midway,left] {}(S1);
\draw[->,transform canvas={xshift=\shift}](H1) to node[midway,right] {}(S1);
\draw[->,transform canvas={xshift=-\shift}](H2) to node[midway,left] {}(S2);
\draw[->,transform canvas={xshift=\shift}](H2) to node[midway,right] {}(S2);
\draw[->](Q) to node[pos=0.25, below] {$\quad\text{ }\textrm{pr}_{S}$} (S1);
\draw[->] (0.5,-0.15) arc (315:30:0.25cm);
\draw[<-] (2.2,0.15) arc (145:-145:0.25cm);
\draw[->](Q) to node[pos=0.25, below] {$\textrm{pr}_{PS}$\text{ }} (S2); 
 
\draw[->,transform canvas={xshift=-\shift}](G1) to node[midway,left] {}(M1);
\draw[->,transform canvas={xshift=\shift}](G1) to node[midway,right] {}(M1);
\draw[->,transform canvas={xshift=-\shift}](G2) to node[midway,left] {}(M2);
\draw[->,transform canvas={xshift=\shift}](G2) to node[midway,right] {}(M2);
\draw[->](P) to node[pos=0.25, below] {$\text{ }\text{ }\alpha_1$} (M1);
\draw[->] (0.8,-3.15) arc (315:30:0.25cm);
\draw[<-] (1.9,-2.85) arc (145:-145:0.25cm);
\draw[->](P) to node[pos=0.25, below] {$\alpha_2$\text{ }} (M2);
\end{tikzpicture}
\end{center}
between these Hamiltonian spaces.

\end{ex} 
As in the symplectic case \cite{Mol1}, any Morita equivalence between Hamiltonian actions is, in fact, of this form (up to isomorphism). So, this notion of Morita equivalence gives a conceptual and symmetric way of capturing the relationship between the Hamiltonian spaces $J$ and $P_*(J)$. Heuristically speaking, this notion of equivalence preserves the transverse part of the geometry of Hamiltonian spaces, in the same way that Morita equivalences of Lie or quasi-symplectic groupoids preserve their transverse geometry. The following lemma gives a precise instance of this, which will be key in the proofs in this and the next section. Recall that a Morita equivalence $(P,\alpha_1,\alpha_2)$ between Lie groupoids $\G_1$ and $\G_2$ induces a homeomorphism $h_P$ between the orbit spaces of $\G_1$ and $\G_2$, which relates a $\G_1$-orbit $\L_1$ with the unique $\G_2$-orbit $\L_2$ for which $\alpha_1^{-1}(\L_1)=\alpha_2^{-1}(\L_2)$.
\begin{lemma}\label{lemma:Moritainv:iacanstrat} Consider a Morita equivalence between two Hamiltonian spaces, denoted as above. 
\begin{itemize}
\item[(a)] The homeomorphisms $h_Q:\underline{S}_1\to \underline{S}_2$ and $h_P:\underline{M}_1\to\underline{M}_2$ fit into a commutative square
\begin{center}
\begin{tikzcd} 
\underline{S}_1\arrow[d,"\underline{J}_1"] \arrow[r,"h_Q"]  & \underline{S}_2 \arrow[d,"\underline{J}_2"]\\
\underline{M}_1 \arrow[r, "h_P"] & \underline{M}_2
\end{tikzcd}
\end{center} 
\item[(b)] For each $p\in P$, we have an associated isomorphism of Lie groups
\begin{equation}\label{eqn:indisoisotrgps:moreq} \Phi_p:\G_{x_1}\xrightarrow{\sim} \G_{x_2},\quad x_1:=\alpha_1(p),\quad x_2:=\alpha_2(p),
\end{equation} defined by the relation
\begin{equation*} g\cdot p=p\cdot \Phi_p(g).
\end{equation*} Given $q\in j^{-1}(p)$, \eqref{eqn:indisoisotrgps:moreq} identifies $\G_{p_1}$ with $\G_{p_2}$ for $p_1:=\beta_1(q)$ and $p_2:=\beta_2(q)$.  
\item[(c)] The homeomorphism $h_P$ identifies $\sS_{\G_1}(\underline{M}_1)$ with $\sS_{\G_2}(\underline{M}_2)$ and restricts to an integral affine isomorphism between strata.  
\item[(d)] Let $\underline{\Sigma}_1\in \sS_{\G_1}(\underline{M}_1)$ and $\underline{\Sigma}_2\in \sS_{\G_2}(\underline{M}_2)$ such that $h_P(\underline{\Sigma}_1)=\underline{\Sigma}_2$, $p\in \alpha_1^{-1}(\Sigma_1)=\alpha_2^{-1}(\Sigma_2)$ and $q\in j^{-1}(p)$. The differential of $h_P:\underline{\Sigma}_1\to \underline{\Sigma}_2$ at $\L_{x_1}$ identifies $(\g_{p_1}^0)^{\G_{x_1}}$ with $(\g_{p_2}^0)^{\G_{x_2}}$ (viewed as subspaces of $T_{\L_{x_1}}\underline{\Sigma}_1$ and $T_{\L_{x_2}}\underline{\Sigma}_2$, respectively, via \eqref{eqn:tgntsp:stratum:orbsp:quasisympgpoid}), with $x_1,x_2$ and $p_1,p_2$ as in part $\mathrm{(b)}$.  
\end{itemize}
\end{lemma}
\begin{proof} Parts (a) and (b) are particular cases of the first two parts of \cite[Proposition 1.52]{Mol1} (and are readily verified). The fact that $h_P$ identifies $\sS_{\G_1}(\underline{M}_1)$ with $\sS_{\G_2}(\underline{M}_2)$ is immediate from the fact that isotropy groups at points through $P$-related orbits are isomorphic. Suppose that $\underline{\Sigma}_1\in \sS_{\G_1}(\underline{M}_1)$ and $\underline{\Sigma}_2\in \sS_{\G_2}(\underline{M}_2)$ such that $h_P(\underline{\Sigma}_1)=\underline{\Sigma}_2$. Since the diagram
\begin{center}
\begin{tikzcd} & \alpha^{-1}(\Sigma_1)=\alpha_2^{-1}(\Sigma_2)\arrow[ld,"\alpha_1"']\arrow[rd,"\alpha_2"] & \\
\Sigma_1\arrow[d,"\pi_1"'] & & \Sigma_2\arrow[d,"\pi_2"]  \\
\underline{\Sigma}_1\arrow[rr,"h_P"'] & & \underline{\Sigma}_2   
\end{tikzcd}
\end{center}
commutes and all vertical and skew arrows are surjective submersions, the map $h_P:\underline{\Sigma}_1\to \underline{\Sigma}_2$ is indeed a diffeomorphism and, for any $p\in P$, the diagram
\begin{center}
\begin{tikzcd} T_{\L_{x_1}}\underline{\Sigma}_1 \arrow[r,"\d(h_P)"]\arrow[d,"\eqref{eqn:tgntsp:stratum:orbsp:Liegpoid}"'] & T_{\L_{x_2}}\underline{\Sigma}_2\arrow[d,"\eqref{eqn:tgntsp:stratum:orbsp:Liegpoid}"] \\
\No_{x_1}^{\G_{x_1}}\arrow[r,"\Psi_p"'] & \No_{x_2}^{\G_{x_2}}
\end{tikzcd}
\end{center} commutes as well, where $x_1:=\alpha_1(p)$, $x_2:=\alpha_2(p)$ and $\Psi_p:\No_{x_1}\to \No_{x_2}$ is the $\Phi_p$-equivariant linear isomorphism that maps $[v]$ to $[\d \alpha_2(\widehat{v})]$, with $\widehat{v}\in T_pP$ such that $\d\alpha_1(\widehat{v})=v$. The diagram
\begin{center}
\begin{tikzcd} 
\No_{x_1}\arrow[r,"\Psi_p"] & \No_{x_2} \\
\g^*_{x_1}\arrow[r,"(\Phi_p)^*"']\arrow[u,"(\rho_{\Omega_1})_{x_1}^*"] & \g^*_{x_2} \arrow[u,"(\rho_{\Omega_2})_{x_2}^*"']   
\end{tikzcd}
\end{center} 
commutes as well (this is \cite[Lemma 1.37]{Mol2} in the symplectic case, and the proof given there readily extends to the quasi-symplectic case). Hence, so does the diagram
\begin{center}
\begin{tikzcd} T_{\L_{x_1}}\underline{\Sigma}_1 \arrow[r,"\d(h_P)"]\arrow[d,"\eqref{eqn:tgntsp:stratum:orbsp:quasisympgpoid}"'] & T_{\L_{x_2}}\underline{\Sigma}_2\arrow[d,"\eqref{eqn:tgntsp:stratum:orbsp:quasisympgpoid}"] \\
(\g^*_{x_1})^{\G_{x_1}}\arrow[r,"(\Phi_p)^*"'] & (\g^*_{x_2})^{\G_{x_2}}
\end{tikzcd}
\end{center} The fact that $h_P:\underline{\Sigma}_1\to \underline{\Sigma}_2$ is an integral affine isomorphism (for which we ought to show that $h_P^*(\Lambda_{\underline{\Sigma}_2})=\Lambda_{\underline{\Sigma}_1}$) follows from commutativity of this diagram, because for all $p,x_1$ and $x_2$ as above
\begin{equation*} (\Phi_p)_*\big(\exp_{Z(\G_{x_1})}^{-1}(C_{\g_{x_1}})\big)=\exp_{Z(\G_{x_2})}^{-1}(C_{\g_{x_2}}).
\end{equation*} Moreover, part (d) follows from commutativity of this diagram and part (b).
\end{proof}
With this at hand, we turn to the extension of the constructibility theorem. 
\begin{proof}[Proof of Theorem \ref{thm:stratification:image:momentum:map:groupoids}] To prove the theorem, we will first show that Lemma \ref{lemma:localtoglobal:abstractaffstrat} applies to $Y:=\underline{M}$, $\P:=\sS_\G(\underline{M})$, $X:=\underline{\Delta}$ and $\D$ given by 
\begin{equation*} \D_{\L_x}:=\bigcap_{p\in J^{-1}(x)} (\g_p^0)^{\G_x}.
\end{equation*}
To this end, let $x\in \Delta$ and $G:=\G_x$. Due to the linearization theorem for proper symplectic groupoids, there is an invariant open $V$ in $M$ around $x$, an invariant open $V_{\g^*}$ in $\g^*$ around the origin and a quasi-symplectic Morita equivalence 
\begin{center}
\begin{tikzpicture} \node (G1) at (0,0) {$(\G,\Omega,\phi)\vert_V$};
\node (M1) at (0,-1.3) {$V$};
\node (S) at (3.3,0) {$(P,\omega_P)$};
\node (M2) at (6.6,-1.3) {$V_{\g^*}$};
\node (G2) at (6.6,0) {$(G\ltimes \g^*,-\d \lambda_{\textrm{can}})\vert_{V_{\g^*}}$};
 
\draw[->,transform canvas={xshift=-\shift}](G1) to node[midway,left] {}(M1);
\draw[->,transform canvas={xshift=\shift}](G1) to node[midway,right] {}(M1);
\draw[->,transform canvas={xshift=-\shift}](G2) to node[midway,left] {}(M2);
\draw[->,transform canvas={xshift=\shift}](G2) to node[midway,right] {}(M2);
\draw[->](S) to node[pos=0.25, below] {$\quad\quad\alpha_1$} (M1);
\draw[->] (2.25,-0.15) arc (315:30:0.25cm);
\draw[<-] (4.35,0.15) arc (145:-145:0.25cm);
\draw[->](S) to node[pos=0.25, below] {$\alpha_2$\text{ }} (M2);
\end{tikzpicture}
\end{center}
that relates $\L_x$ to the origin in $\g^*$ (this follows from \cite[Corollary 2.6]{Zu} by restricting to a slice for $\G$ through $x$, like in the proof of \cite[Theorem 2.1]{CrFeTo2}). As in Example \ref{ex:Moreq:associatedHamsp}, the restriction of $J$ to $V$ has an associated Hamiltonian $(G\ltimes \g^*)\vert_{V_{\g^*}}$-space. The data of a Hamiltonian $(G\ltimes \g^*)\vert_{V_{\g^*}}$-space is the same as that of a Hamiltonian $G$-space with momentum map taking values in $V_{\g^*}$ (see, e.g., \cite[Example 1.34, Example 1.35]{Mol1}). By Lemma \ref{lemma:Moritainv:iacanstrat}, the stratification obtained by applying Theorem \ref{thm:stratification:image:momentum:map:disconnectedgps} to this Hamiltonian $G$-space transports via $h_P$ to a stratification $\sS_{\underline{\Delta}\cap \underline{V}}$ of $\underline{\Delta}\cap \underline{V}$ into integral affine submanifolds of the members of $\sS_\G(\underline{M})$, with the property that $T_{\L}\underline{\Sigma}=\D_{\L}$ for each $\underline{\Sigma}\in \sS_{\underline{\Delta}\cap \underline{V}}$ and $\L\in \underline{\Sigma}$. So, Lemma \ref{lemma:localtoglobal:abstractaffstrat} indeed applies and, hence, there is a (necessarily unique) stratification $\sS_\mathrm{Ham}(\underline{\Delta})$ of $\underline{\Delta}$ into integral affine submanifolds of the members of $\sS_\G(\underline{M})$ with tangent spaces given by \eqref{eqn:definingproperty:stratofimage:gpoids}. \\

To address the equivariant local triviality, suppose that $\underline{\sigma}\in \mathcal{S}_\textrm{Ham}(\underline{\Delta})$ and $x\in \sigma$. Let $G:=G_x$ and $Q\to \L$ the principal $G$-bundle over the $\G$-orbit through $x$ with total space the source-fiber of $\G$ over $x$. We denote by $(\G_\textrm{lin},\Omega_\textrm{lin})$ the linear local model of $(\G,\Omega,\phi)$ around $\L$ as in \cite[Section 2]{CrFeTo2}. By the linearization theorem for quasi-symplectic groupoids (see \cite{CrFeTo2}) there is an open $V$ around $\L$ in $M$, an open $V_\textrm{lin}$ around $\L$ (viewed as $Q\times_G\{0\}$) in $Q\times_G\g^*$, a $2$-form $\beta$ on $V_\textrm{lin}$ and an isomorphism of quasi-symplectic groupoids $\widetilde{\varphi}$ between $(\G,\Omega,\phi)\vert_V$ and $(\G_\textrm{lin},\Omega_\textrm{lin}+t^*\beta-s^*\beta,\textrm{pr}_\L^*\phi-\d\beta)\vert_{V_\textrm{lin}}$ that restricts to the identity map on $\G_\L$. By construction of $(\G_\textrm{lin},\Omega_\textrm{lin})$, there is a canonical quasi-symplectic Morita equivalence 
\begin{center}
\begin{tikzpicture} \node (G1) at (-1.5,0) {$(\G_\textrm{lin},\Omega_\textrm{lin}+t^*\beta-s^*\beta,\textrm{pr}_\L^*\phi-\d\beta)\vert_{V_\textrm{lin}}$};
\node (M1) at (-1.5,-1.3) {$V_\textrm{lin}$};
\node (S) at (3.3,0) {$(P,\omega_P)$};
\node (M2) at (8.1,-1.3) {$V_{\g^*}$};
\node (G2) at (8.1,0) {$(G\ltimes \g^*,-\d \lambda_{\textrm{can}})\vert_{V_{\g^*}}$};
 
\draw[->,transform canvas={xshift=-\shift}](G1) to node[midway,left] {}(M1);
\draw[->,transform canvas={xshift=\shift}](G1) to node[midway,right] {}(M1);
\draw[->,transform canvas={xshift=-\shift}](G2) to node[midway,left] {}(M2);
\draw[->,transform canvas={xshift=\shift}](G2) to node[midway,right] {}(M2);
\draw[->](S) to node[pos=0.25, below] {$\quad\quad\textrm{pr}_{V_\textrm{lin}}$} (M1);
\draw[->] (2.25,-0.15) arc (315:30:0.25cm);
\draw[<-] (4.35,0.15) arc (145:-145:0.25cm);
\draw[->](S) to node[pos=0.25, below] {$\textrm{pr}_{\g^*}$\text{ }} (M2);
\end{tikzpicture}
\end{center}
in which
\begin{itemize}
\item $V_{\mathfrak{g}^*}$ is a $G$-invariant open in $\g^*$ around the origin, 
\item $P$ is the $G$-invariant open in $Q\times \g^*$ around $Q\times \{0\}$ corresponding to $V_\textrm{lin}$, 
\item the action of $G\ltimes\g^*$ along $\textrm{pr}_{\g^*}$ is that induced by the diagonal $G$-action,
\item $\omega_P:=\omega_\theta+(\textrm{pr}_{V_\textrm{lin}})^*\beta$, with $\theta\in \Omega^1(Q,\g)$ some connection $1$-form and 
\begin{equation*}
\omega_\theta:=(\textrm{pr}_{\L})^*\omega_\L-\d\widehat{\theta}\in \Omega^2(Q\times \g^*),
\end{equation*} where $\widehat{\theta}\in \Omega^1(Q\times \g^*)$ is given by $\widehat{\theta}_{(p,\alpha)}=\langle\alpha,\theta_p\rangle$.
\end{itemize} Composing this with the isomorphism $\widetilde{\varphi}$ we obtain a quasi-symplectic Morita equivalence
\begin{center}
\begin{tikzpicture} \node (G1) at (0,0) {$(\G,\Omega,\phi)\vert_V$};
\node (M1) at (0,-1.3) {$V$};
\node (S) at (3.3,0) {$(P,\omega_P)$};
\node (M2) at (6.6,-1.3) {$V_{\g^*}$};
\node (G2) at (6.6,0) {$(G\ltimes \g^*,-\d \lambda_{\textrm{can}})\vert_{V_{\g^*}}$};
 
\draw[->,transform canvas={xshift=-\shift}](G1) to node[midway,left] {}(M1);
\draw[->,transform canvas={xshift=\shift}](G1) to node[midway,right] {}(M1);
\draw[->,transform canvas={xshift=-\shift}](G2) to node[midway,left] {}(M2);
\draw[->,transform canvas={xshift=\shift}](G2) to node[midway,right] {}(M2);
\draw[->](S) to node[pos=0.25, below] {$\quad\alpha$} (M1);
\draw[->] (2.25,-0.15) arc (315:30:0.25cm);
\draw[<-] (4.35,0.15) arc (145:-145:0.25cm);
\draw[->](S) to node[pos=0.25, below] {$\textrm{pr}_{\g^*}$\text{ }} (M2);
\end{tikzpicture}
\end{center} that relates $\L$ to the origin in $\g^*$. Let $J_{\g^*}:(S_{\g^*},\omega_{\g^*})\to \g^*$ be the associated Hamiltonian $G$-space and $\Delta_{\g^*}$ its image. In view of Lemma \ref{lemma:Moritainv:iacanstrat}, $h_P$ maps an open in $\underline{\sigma}\cap \underline{V}$ around $\L$ onto the stratum $\underline{\sigma}_{\g^*}\in\sS_\mathrm{Ham}(\underline{\Delta}_{\g^*})$ through the origin. Because the corresponding invariant subset of $\g^*$ is contained in the $G$-fixed point set, from Theorem \ref{thm:stratification:image:momentum:map:disconnectedgps} we obtain a $G$-equivariant topological trivialization 
\begin{center} \begin{tikzcd} J_{\g^*}^{-1}(W_{\g^*})\arrow[rr,"\Phi_{\g^*}", "\sim"'] \arrow[rd,"J_{\g^*}"'] & & J_{\g^*}^{-1}(0)\times W_{\g^*} \arrow[ld,"\mathrm{pr}_2"] \\
  & W_{\g^*} &
\end{tikzcd} 
\end{center} with $W_{\g^*}$ an invariant open neighbourhood of the origin in $\sigma_{\g^*}$. This transports to a $\widehat{\varphi}$-equivariant topological trivialization
\begin{center} 
\begin{tikzcd} 
J^{-1}(W)\arrow[rr,"\Phi", "\sim"'] \arrow[d,"J"'] & & (J\times\textrm{id}_{\underline{W}})^{-1}(O) \arrow[hookrightarrow,r] & J^{-1}(\L)\times \underline{W} \arrow[d,"J\times\mathrm{id}_{\underline{W}}"] \\
W \arrow[rd,"\pi"']\arrow[rr,"\varphi", "\sim"'] & & O \arrow[hookrightarrow,r]& \L\times \underline{W}\arrow[lld,"\mathrm{pr}_2"] \\
  & \underline{W} & &
\end{tikzcd} 
\end{center} 
over the open $\underline{W}:=h_P^{-1}(\underline{W}_\g^*)$ in $\underline{\sigma}$ around $\L$. Here $W$ is the $\G\vert_V$-invariant open in $V$ corresponding to $\underline{W}$ and $\widehat{\varphi}:\G\vert_W\to \G_\L\times \underline{W}$ is the restriction of $\widetilde{\varphi}$ (where we use that $\G_\L\times \underline{W}$ canonically embeds into $\G_\textrm{lin}$, since $W_{\g^*}$ is contained in the $G$-fixed point set of $\g^*$). Moreover, $\varphi$ is the diffeomorphism covered by $\widehat{\varphi}$ and $O$ is its image. To define $\Phi$, we use that (by \cite{Xu1}) the Hamiltonian $(\G,\Omega,\phi)\vert_V$-space $J:(J^{-1}(V),\omega)\to V$ is canonically isomorphic to the Hamiltonian $(\G,\Omega,\phi)\vert_V$-space 
 \begin{equation*} P_*(J_{\g^*}):(P\ast_GS_{\g^*},\omega_{PS_{\g^*}})\to V,\quad [p_P,p_{S}]\mapsto\alpha(p_P),
 \end{equation*} constructed from $J_{\g^*}$ as in Example \ref{ex:Moreq:associatedHamsp}, but with the roles of left and right reversed. The trivialization $\Phi_{\g^*}$ transports to a $\widehat{\varphi}$-equivariant homeomorphism 
\begin{equation*} P_*(J_{\g^*})^{-1}(W)
\xrightarrow{\sim} (P(J_{\g^*})\times \textrm{id}_{\underline{W}})^{-1}(O)\subset P(J_{\g^*})^{-1}(\L)\times \underline{W},
\end{equation*} given by the composition of 
\begin{equation*} P_*(\Phi_{\g^*}):P_*(J_{\g^*})^{-1}(W)=P\ast_GJ_{\g^*}^{-1}(W_{\g^*})\xrightarrow{\sim} P\ast_G(J_{\g^*}^{-1}(0)\times W_{\g^*}),\quad [p_P,p_S]\mapsto [p_P,\Phi_{\g^*}(p_S)],
\end{equation*} with 
\begin{equation*} P\ast_G(J_{\g^*}^{-1}(0)\times W_{\g^*})\xrightarrow{\sim} 
(P(J_{\g^*})\times \textrm{id}_{\underline{W}})^{-1}(O),\quad [q,\eta,p_S,\eta]\mapsto \big([q,0,p_S],(\pi\circ\alpha)(q,\eta)\big),
\end{equation*} 
and the corresponding $\widehat{\varphi}$-equivariant homeomorphism 
\begin{equation*} \Phi: J^{-1}(W)\xrightarrow{\sim} (J\times\textrm{id}_{\underline{W}})^{-1}(O)
\end{equation*} fits into the desired commutative diagram. Following this construction, it is readily verified that smooth extensions of $\Phi_{\g^*}$ and $\Phi_{\g^*}^{-1}$ like in Theorem \ref{thm:localtriviality} can also be transported to obtain smooth extensions of the same type for $\Phi$ and $\Phi^{-1}$ if the extensions of $\Phi_{\g^*}$ and $\Phi_{\g^*}^{-1}$ are $G$-equivariant. From the construction of the local trivializations for Hamiltonian $G$-spaces it is clear that $\Phi_{\g^*}$ can be chosen so that it admits such extensions and, hence, the local trivializations of $J:J^{-1}(\sigma)\to \sigma$ can be chosen so that they admit the desired smooth extensions.  
 \end{proof}
Finally, let us point out the following, which is immediate from Lemma \ref{lemma:Moritainv:iacanstrat}.
\begin{cor}\label{cor:morinv:canstratim} For a Morita equivalence as in Lemma \ref{lemma:Moritainv:iacanstrat}, the homeomorphism $h_P:\underline{M}_1\to \underline{M}_2$ restricts to a homeomorphism between the respective images $\underline{\Delta}_1$ and $\underline{\Delta}_2$ of $\underline{J}_1$ and $\underline{J}_2$, which identifies $\sS_\mathrm{Ham}(\underline{\Delta}_1)$ with $\sS_\mathrm{Ham}(\underline{\Delta}_2)$ and restricts to integral affine isomorphisms between the strata. 
\end{cor}
\subsubsection{Pushforwards of constant sheaves}\label{sec:pushforwardsheaves:gpoidactions} Moving on to the extension of Corollary \ref{cor:constructibility:pushforwards}(a), as mentioned before, the role of $EG$ will now be taken by $E\G$, by which we mean the total space of the universal principal $\G$-bundle $E\G\to B\G$ obtained by extending Milnor's construction for topological groups to topological groupoids (as in \cite{Haef}). Explicitly, $E\G$ consists of equivalence classes of sequences $(t_0,g_0,t_1,g_1,...)$ with $t_i\in [0,1]$ and $g_i\in \G$ such that 
\begin{itemize} \item[(i)] $t_i\neq 0$ for only finitely many $i$ and $\sum_{i\geq 0} t_i=1$,
\item[(ii)] $t(g_i)=t(g_j)$ for all $i,j$,
\end{itemize}
where two such sequences $(t_0,g_0,t_1,g_1,...)$ and $(s_0,h_0,s_1,h_1,...)$ are considered equivalent if $t_i=s_i$ for all $i$ and $g_i=h_i$ whenever $t_i\neq 0$. The topology on $E\G$ is the weakest one making the maps 
\begin{equation*} f_i:E\G\to [0,1],\quad [t_0,g_0,t_1,g_1,...]\mapsto t_i,\quad\quad f_i^{-1}(]0,1])\to \G,\quad [t_0,g_0,t_1,g_1,...]\mapsto g_i,
\end{equation*} continuous for every $i$. The $\G$-action is along the continuous map
\begin{equation*} \mu:E\G\to M,\quad [t_0,g_0,t_1,g_1,...]\mapsto t(g_0),
\end{equation*} and is given by
\begin{equation*} g\cdot [t_0,g_0,t_1,g_1,...]=[t_0,gg_0,t_1,gg_1,...].
\end{equation*} 
Let $B\G:=E\G/\G$ denote the orbit space of this action. Its cohomology is that of the stack presented by $\G$ (see, e.g., \cite[Section 2]{BeGiNoXu}).
\begin{ex} For the action groupoid $G\ltimes S$ of a $G$-space $S$, there is a canonical homeomorphism
\begin{equation*} B(G\ltimes S)\cong EG\times _GS:=(EG\times S)/G.
\end{equation*} So, the cohomology of $B(G\ltimes S)$ is the equivariant cohomology of the $G$-space. More generally, for the action groupoid $\G\ltimes S$ of a $\G$-space $S\to M$, the homeomorphism
\begin{equation*} E(\G\ltimes S)\to E\G\times_MS, \quad [t_0,(h_0,p_0),t_1,(h_1,p_1),...]\mapsto ([t_0,h_0,t_1,h_1,...],h_0\cdot p_0),
\end{equation*}  descends to a homeomorphism
\begin{equation}\label{eqn:classfsp:borelmodel:gpoid} B(\G\ltimes S)\cong E\G\ast _\G S:=(E\G\times_MS)/\G.
\end{equation} 
\end{ex}
Being an invariant of the stack presented by $\G$, the cohomology of $B\G$ only depends on the Morita equivalence class of $\G$. We will use the following version of this fact.
\begin{prop}\label{prop:moreq:indisoshfcoh} A Morita equivalence of Lie groupoids
\begin{center}
\begin{tikzpicture} \node (G1) at (0,0) {$\G_1$};
\node (M1) at (0,-1.3) {$M_1$};
\node (S) at (1.4,0) {$P$};
\node (M2) at (2.7,-1.3) {$M_2$};
\node (G2) at (2.7,0) {$\G_2$};
 
\draw[->,transform canvas={xshift=-\shift}](G1) to node[midway,left] {}(M1);
\draw[->,transform canvas={xshift=\shift}](G1) to node[midway,right] {}(M1);
\draw[->,transform canvas={xshift=-\shift}](G2) to node[midway,left] {}(M2);
\draw[->,transform canvas={xshift=\shift}](G2) to node[midway,right] {}(M2);
\draw[->](S) to node[pos=0.25, below] {$\text{ }\text{ }\alpha_1$} (M1);
\draw[->] (0.8,-0.15) arc (315:30:0.25cm);
\draw[<-] (1.9,0.15) arc (145:-145:0.25cm);
\draw[->](S) to node[pos=0.25, below] {$\alpha_2$\text{ }} (M2);
\end{tikzpicture}
\end{center}
induces an isomorphism in sheaf cohomology 
\begin{equation*} H^\bullet(B\G_1,\R)\cong H^\bullet(B\G_2,\R). 
\end{equation*}
\end{prop}
\begin{proof}
The Morita equivalence induces an isomorphism between pull-back groupoids
\begin{equation*} \alpha_1^*\G_1\cong \alpha_2^*\G_2, 
\end{equation*} which identifies the respective arrows $(p,g_1,q)\in P\times_{M_1} \G_1\times_{M_1}P$ and $(p,g_2,q)\in  P\times_{M_2} \G_2\times_{M_2}P$ of $\alpha_1^*\G_1$ and $\alpha_2^*\G_2$ for which $g_1\cdot q=p\cdot g_2$. The two projections in the diagram
\begin{center} \begin{tikzcd} &\arrow[ld,"\textrm{pr}_{\G_1}"'] \alpha_1^*\G_1\cong \alpha_2^*\G_2\arrow[rd,"\textrm{pr}_{\G_2}"]  & \\
\G_1 & & \G_2
\end{tikzcd} 
\end{center} are weak equivalences. It is well-known (although we could not find a reference) that, if $\Phi:\G\to \H$ be a weak equivalence of Lie groupoids, then the continuous map 
\begin{equation*}
\Phi_*:B\G\to B\H
\end{equation*}
induced by
\begin{equation*} \Phi_*:E\G\to E\H,\quad [t_0,g_0,t_1,g_1,...]\mapsto [t_0,\Phi(g_0),t_1,\Phi(g_1),...]
\end{equation*} is a homotopy equivalence. Applying this to the above two projections, the proposition follows.
\end{proof}
With the above background at hand, we are now ready to prove the extension of Corollary \ref{cor:constructibility:pushforwards}(a), starting with the following lemma.
\begin{lemma}\label{lemma:moreq:indisoderpushfwd} Let $i\in \N$. A Morita equivalence as in Lemma \ref{lemma:Moritainv:iacanstrat} induces an isomorphism of sheaves 
\begin{equation}\label{eqn:moreq:indisoderpushfwd} R^i(\underline{J}_1^{E\G_1})_*(\underline{\R})\cong R^i(\underline{J}_2^{E\G_2})_*(\underline{\R})
\end{equation} covering the induced homeomorphism $h_P:\underline{\Delta}_1\to \underline{\Delta}_2$ between the images of $\underline{J}_1$ and $\underline{J}_2$.
\end{lemma}
\begin{proof} Given $P$-related opens $\underline{V}_1$ and $\underline{V}_2$, the Morita equivalence can be restricted to $V_1$ and $V_2$ to obtain a Morita equivalence between the restrictions of the given Hamiltonian spaces to $V_1$ and $V_2$, respectively. Applying Proposition \ref{prop:moreq:indisoshfcoh} to the upper part of this restricted Morita equivalence and using the homeomorphism \eqref{eqn:classfsp:borelmodel:gpoid} for both $\G_1\ltimes J^{-1}_1(V_1)$ and $\G_2\ltimes J_2^{-1}(V_2)$, we obtain an isomorphism
\begin{equation*} H^i(E\G_1\ast_{\G_1}J_1^{-1}(V_1),\R)\cong H^i(E\G_2\ast_{\G_2}J_2^{-1}(V_2),\R),
\end{equation*} for each such pair of $P$-related opens. Since these isomorphisms are compatible with restriction to smaller opens, this defines an isomorphism of pre-sheaves
\begin{equation*} H^i\big((\underline{J}_1^{E\G_1})^{-1}(-),\R\big)\cong H^i\big((\underline{J}_2^{E\G_2})^{-1}(-),\R\big)
\end{equation*} covering $h_P:\underline{\Delta}_1\to \underline{\Delta}_2$. After applying sheafification, we obtain the desired isomorphism \eqref{eqn:moreq:indisoderpushfwd}.  
\end{proof}
\begin{proof}[Proof of the extension of Corollary \ref{cor:constructibility:pushforwards}(a)] Let $x\in \Delta$. Consider a quasi-symplectic Morita equivalence as in the proof of Theorem \ref{thm:stratification:image:momentum:map:groupoids} and the induced Hamiltonian Morita equivalence as in Example \ref{lemma:Moritainv:iacanstrat}. Applying Lemma \ref{lemma:moreq:indisoderpushfwd} to this, we obtain an isomorphism of sheaves (covering $h_P$) between the restriction of $R^i(\underline{J}^{E\G})_*(\underline{\R})$ to $\underline{\Delta}\cap \underline{V}$ and $R^i(\underline{J}_{\g^*}^{EG})_*(\underline{\R})$, where $J_{\g^*}$ denotes the associated Hamiltonian $G$-space (as in the proof of Theorem \ref{thm:stratification:image:momentum:map:groupoids}). By the extension of Corollary \ref{cor:constructibility:pushforwards}(a) to actions of possibly disconnected compact Lie groups, $R^i(\underline{J}_{\g^*}^{EG})_*(\underline{\R})$ is constructible with respect to $\sS_\mathrm{Ham}(\underline{\Delta}_{\g^*})$, meaning that its stalks are finite-dimensional and its restriction to any stratum is locally constant. Therefore, the stalks of $R^i(\underline{J}^{E\G})_*(\underline{\R})$ near $\L_x$ are finitely-dimensional as well and its restriction to the stratum of $\sS_\mathrm{Ham}(\underline{\Delta})$ through $\L_x$ is constant near $\L_x$. This being true for all $x\in \Delta$, we conclude that $R^i(\underline{J}^{E\G})_*(\underline{\R})$ is indeed constructible with respect to $\sS_\mathrm{Ham}(\underline{\Delta})$.
\end{proof}
Next, we turn to the extension of Proposition \ref{prop:specializationmaps}. 
\begin{proof}[Proof of the extension of Proposition \ref{prop:specializationmaps}] Replacing $\gamma$ by the path $\gamma^\varepsilon$ given by $\gamma^\varepsilon(t)=\gamma((1-t)\varepsilon+t)$, for small enough $\varepsilon >0$, we can assume without loss of generality that the path $\gamma$ is contained in an open $\underline{V}$ around $\L_y$ for which there is a quasi-symplectic Morita equivalence as in the proof of Theorem \ref{thm:stratification:image:momentum:map:groupoids}. Consider the induced Hamiltonian Morita equivalence as in Example \ref{ex:Moreq:associatedHamsp}. For this Morita equivalence, the homeomorphism $h_Q$ as in Lemma \ref{lemma:Moritainv:iacanstrat}(a) restricts to a homeomorphism between reduced space of $J$ at $\L_y$ and the reduced space of the associated Hamiltonian $G$-space $J_{\g^*}$ at $0$, and between those at $\L_x$ and the coadjoint $G$-orbit that is $P$-related to $\L_x$. Moreover, the induced maps between the cohomology groups of these reduced spaces intertwines the map between the stalks of $R^i(\underline{J})_*(\underline{\R})$ associated to $\gamma$ with that of $R^i(\underline{J}_{\g^*})_*(\underline{\R})$ associated to $h_P(\gamma)$. In view of Lemma \ref{lemma:Moritainv:iacanstrat}(c) and Corollary \ref{cor:morinv:canstratim}, Proposition \ref{prop:specializationmaps:disconngps} applies to $J_{\g^*}$ and the path $h_P\circ \gamma$. Transferring the Weyl group action obtained from this via a Lie group isomorphism as in Lemma \ref{lemma:Moritainv:iacanstrat}(c), we obtain a $\mathbb{W}_y$-action on the cohomology of the reduced space of $J$ at $\L_y$ and the map associated to $\gamma$ is an isomorphism onto the $\mathbb{W}_y$-fixed point set. 
\end{proof}
\subsubsection{Sjamaar's complex, linear variation and the symplectic volumes}\label{sec:Sjamaar:deRham:gpoidactions} In this section we elaborate on the extensions of Theorem \ref{thm:linvar:main} and Corollary \ref{cor:sympvolpolynomial}. The statements of these extensions are obtained by making the same substitutions as for their extension to possibly disconnected compact Lie groups, but with the role of $\underline{\g^*}$ now taken by $\underline{M}$. The resulting statements make sense because, as we will now explain, for Hamiltonian actions of proper quasi-symplectic groupoids the reduced spaces are still symplectic, come with an associated cohomology class and have finite symplectic volume when they are compact. Consider such an action of a proper quasi-symplectic groupoid $(\G,\Omega,\phi)$ along $J:(S,\omega)\to M$. As before, for $\L\in \underline{M}$, we use the notation
\begin{equation*} S_\L:=J^{-1}(\L), \quad \underline{S}_\L:=J^{-1}(\L)/\G.
\end{equation*} We let $\omega_\L$ denote the unique $2$-form on $\L$ such that $t^*\omega_\L-s^*\omega_\L=\Omega\vert_{\G_\L}$, where $t,s:\G_\L\to \L$ are the target and source map of the restriction of $\G$ to $\L$. Moreover, for $x\in M$ we denote
\begin{equation*} S_x:=J^{-1}(x), \quad \underline{S}_x:=J^{-1}(x)/\G_x.
\end{equation*}
To start with, the following still holds.
\begin{prop}\label{prop:sympstrredsp} For any $\L\in \underline{M}$ (respectively $x\in M$), the reduced space $\underline{S}_\L$ (respectively $\underline{S}_x$) has a natural stratification. For each stratum $\underline{\Sigma}$ of this stratification, the corresponding $\G$-invariant (respectively $\G_x$-invariant) subset $\Sigma$ of $S_\L$ (respectively $S_x$) is a submanifold of $S$ and the quotient map $q_\Sigma:\Sigma\to \underline{\Sigma}$ is a smooth submersion. Moreover, each such $\underline{\Sigma}$ comes with a symplectic form $\omega_{\underline{\Sigma}}$ that is uniquely determined by the fact that its pull-back along $q_\Sigma$ is the restriction of $\omega{\vert_\Sigma}-(J\vert_\Sigma)^*\omega_\L$ (respectively $\omega\vert_\Sigma$). When $x\in \L$, these stratifications and symplectic forms are identified with each other by the canonical homeomorphism between $\underline{S}_x$ and $\underline{S}_\L$ (induced by the inclusion $S_x\hookrightarrow S_\L$). 
\end{prop}
\begin{proof} The partition underlying the stratification of $\underline{S}_x$ consists of the connected components of the members of the partition defined by declaring the $\G_x$-orbits $\O_p$ and $\O_q$ through $p,q$ in $S_x$ to be equivalent when $\G_p$ and $\G_q$ are isomorphic Lie groups. For the reduced space at zero of a Hamiltonian action of a compact Lie group this partition coincides with that underlying the stratification of \cite{LeSj}, as follows from standard arguments similar to those in the proof of \cite[Proposition 4.12]{CrMe}. Therefore, the fact that this is a stratification and the claims about it in the proposition follow by combining the analogous results in \cite{LeSj} with \cite[Proposition 3.4]{Zu}. For the stratification on $\underline{S}_\L$, first consider the stratification $\mathcal{S}_{\G\ltimes S}(\underline{S})$ of $\underline{S}$ associated to the action groupoid $\G\ltimes S$ (as in \cite{PfPoTa,CrMe}). Explicitly, this consists of the connected components of the members of the partition $\mathcal{P}_{\cong}$ defined by declaring the $\G$-orbits $\O_p$ and $\O_q$ through $p,q$ in $S$ to be equivalent when $\G_p$ and $\G_q$ are isomorphic Lie groups. By restricting the partition $\mathcal{P}_{\cong}$ to the subspace $\underline{S}_\L$ and taking connected components of its members, we obtain the desired stratification of $\underline{S}_\L$. One way to see that this is indeed a stratification is that the canonical homeomorphism between $\underline{S}_x$ and $\underline{S}_\L$ takes it to the above stratification on $\underline{S}_x$. The proof of the fact that the strata in $\underline{S}_x$ are indeed smooth manifolds along the lines indicated above shows that, in fact, each stratum in $\underline{S}_x$ is a submanifold of the stratum of $\mathcal{S}_{\G\ltimes S}(\underline{S})$ that it is contained in. Hence, the same holds for the strata in $\underline{S}_\L$. Since for each stratum $\underline{\Sigma}_S$ of $\mathcal{S}_{\G\ltimes S}(\underline{S})$ the corresponding $\G$-invariant subset $\Sigma_S$ is a submanifold in $S$ and the quotient map $\Sigma_S\to \underline{\Sigma}_S$ is a smooth submersion (cf. \cite[Section 4.6]{CrMe}), the same holds for the strata in $\underline{S}_\L$. To conclude the proof, note that for each such stratum $\underline{\Sigma}$, the symplectic form $\omega_{\underline{\Sigma}}$ coming from the identification with the corresponding stratum in $\underline{S}_x$ pulls-back to the form $\omega\vert_{\Sigma}-(J\vert_\Sigma)^*\omega_\L$ along $q_\Sigma$, because $\omega\vert_{\Sigma}-(J\vert_\Sigma)^*\omega_\L$ restricts to $\omega$ on the $\G_x$-invariant subset of $S_x$ corresponding to $\underline{\Sigma}$ and because it is a basic form for the action groupoid $\G_\L\ltimes \Sigma$ of the induced $\G_\L$-action on $\Sigma$ along $J\vert_\Sigma$ (meaning that the pull-backs of $\omega\vert_{\Sigma}-(J\vert_\Sigma)^*\omega_\L$ along the source and target of $\G_\L\ltimes \Sigma$ coincide \cite{Wa}). 
\end{proof}
Given this proposition, Sjamaar's de Rham complex on $\underline{S}_\L$ can be defined just as in the case of compact Lie group actions (see \S\ref{subsec:deRhammodel:1}) and the reduced symplectic form on the union of open strata still defines a cohomology class $\varpi_\L\in H^2_\mathrm{dR}(\underline{S}_\L)$. The same goes for $\underline{S}_x$ and the caonical identification with $\underline{S}_\L$ still induces an isomorphism of de Rham complexes. Moreover, in view of \cite[Proposition 3.4]{Zu} both Theorem \ref{thm:sjamaar:deRham} and \ref{thm:sjamaar:finvol} hold as well for Hamiltonian actions of proper quasi-symplectic groupoids. For the interpretation of linear variation and polynomiality of the symplectic volume functions, the same comments apply as for the extension to actions of disconnected compact Lie groups (see the last paragraph of \S\ref{sec:extensionpf:disconnectedgps}). This explains the ingredients going into the statements of the extensions of Theorem \ref{thm:linvar:main} and Corollary \ref{cor:sympvolpolynomial}. For their proofs, we will use the following. 
\begin{lemma}\label{lemma:deRhamcomplex:Moritainv} Given a Morita equivalence between Hamiltonian spaces of proper quasi-symplectic groupoids, denoted as in \S\ref{sec:pf:constructibility:quasisymgpoids}, and $P$-related orbits $\L_1\subset M_1$ and $\L_2\subset M_2$, the homeomorphism between $\underline{S}_{\L_1}$ and $\underline{S}_{\L_2}$ obtained by restricting $h_Q$ (see Lemma \ref{lemma:Moritainv:iacanstrat}) is an isomorphism of stratified spaces and $h_Q^*:\Omega^\bullet(\underline{S}^\mathrm{prin}_{\L_2})\xrightarrow{\sim} \Omega^\bullet(\underline{S}^\mathrm{prin}_{\L_1})$ restricts to an isomorphism between the de Rham complexes. Moreover, $h_Q^*$ identifies the reduced symplectic forms on $\underline{S}_{\L_1}$ and $\underline{S}_{\L_2}$. 
\end{lemma}
We leave the straightforward proof of this to the reader and conclude with the following.  
\begin{proof}[Proof of the extensions of Theorem \ref{thm:linvar:main} and Corollary \ref{cor:sympvolpolynomial}] Let $\underline{\sigma}\in \sS_\textrm{Ham}(\underline{\Delta})$ and $\L\in \underline{\sigma}$. Consider a Morita equivalence as at the start of the proof of Theorem \ref{thm:stratification:image:momentum:map:groupoids} and let $J_{\g^*}:(S_{\g^*},\omega_{\g^*})\to \g^*$ be the associated Hamiltonian $G$-space. Denote by $\Delta_{\g^*}$ its image, by $\underline{\sigma}_{\g^*}$ the stratum of $\sS_{\mathrm{Ham}}(\underline{\Delta}_{\g^*})$ through the origin and by $\underline{W}$ the open $h_P^{-1}(\underline{\sigma}_{\g^*})$ in $\underline{\sigma}$. The homeomorphisms $h_P$ and $h_Q$ in Lemma \ref{lemma:Moritainv:iacanstrat} give an isomorphism of topological fiber bundles
\begin{center}
\begin{tikzcd} 
\underline{J}^{-1}(\underline{W})\arrow[d,"\underline{J}"] \arrow[r,"h_Q","\sim"']  & \underline{S}_{\g^*} \arrow[d,"\underline{J}_{\g^*}"]\\
\underline{W} \arrow[r, "h_P","\sim"'] & \underline{\sigma}_{\g^*}
\end{tikzcd}
\end{center} By Lemma \ref{lemma:deRhamcomplex:Moritainv} the induced isomorphism between the flat vector bundles given by the degree two cohomology of their fibers identifies the respective sections defined by the reduced symplectic forms and, by Corollary \ref{cor:morinv:canstratim}, $h_P$ is an isomorphism of integral affine manifolds. In view of this, the desired extensions of the linear variation theorem and the polynomiality of the symplectic volume function (which are both of a local nature) follow from the case of compact Lie group actions.  
\end{proof}
\bibliographystyle{plain}

\bibliography{ref}
\Addresses
\end{document}